\documentclass[a4paper,18pt,dvips]{article}  
\usepackage{amsmath, amssymb, amscd,amsthm,amsfonts,amscd}
\usepackage{graphics, graphicx}
\usepackage{float,epsfig,color}
\usepackage[latin1]{inputenc}
\usepackage[T1]{fontenc}
\newtheorem{theorem}{Theorem}[section]
\newtheorem{lemma}[theorem]{Lemma}
\newtheorem{proposition}[theorem]{Proposition}
\newtheorem{corollary}[theorem]{Corollary}
\theoremstyle{definition}
\newtheorem{definition}[theorem]{Definition}
\newtheorem{remark}[theorem]{Remark}
\newtheorem*{acknowledgements}{Acknowledgements}

\newcommand{\restrictto}[1]{\,\rule[-2.5mm]{0.125mm}{5.5mm}_{\rule[-1.0mm]{0mm}{4mm}\, #1}}
\newcommand{\M}{\mathcal{M}}

\newcommand{\HH}{\mathbb{H}}
\newcommand{\Z}{\mathbb{Z}}
\newcommand{\D}{\mathbb{D}}
\newcommand{\DD}{\mathcal{D}}
\newcommand{\LL}{\mathcal{L}}
\newcommand{\bs}{\mathcal{B}}
\newcommand{\rs}{\hat{\mathbb{C}}}
\newcommand{\C}{\mathbb{C}}
\newcommand{\R}{\mathbb{R}}

\newcommand{\hol}{holomorphic}
\newcommand{\eqcl}{equivalence class}
\newcommand{\eqre}{equivalence relation}
\newcommand{\multeq}{multiple equilibrium point}
\newcommand{\eqpt}{equilibrium point}
\newcommand{\cds}{combinatorial data set}

\newcommand{\Arcsin}{\textrm{Arcsin}}
\newcommand{\Res}{\textrm{Res}}
\newcommand{\id}{\textrm{id}}
\newcommand{\Arg}{\textrm{Arg}}
\newcommand{\Log}{\textrm{Log}}
\newcommand{\vf}{\frac{d}{dz}}
\newcommand{\vfw}{\frac{d}{dw}}

\begin{document}
\title{Classification of Complex Polynomial Vector Fields in One Complex Variable}  
\author{Bodil Branner and Kealey Dias}         
\date{\today}    
\maketitle
\begin{center}
Dedicated to Bob Devaney on his 60$^{th}$ birthday, \\
and dedicated to the memory of Adrien Douady.
\end{center}
\begin{abstract}
  The global characterization of structurally stable monic and centered complex polynomial vector fields presented 
in 
~\cite{Sent} is extended to include all monic and centered complex polynomial vector fields in $\C$.
\end{abstract}
\begin{figure}[b]%
\flushleft
\rule[0mm]{54.0mm}{0.15mm}\\ \hspace*{5.0mm}%
{\footnotesize  2000 Mathematics Subject Classification: 37F75, 37C15, 37F20.\par
 }
\end{figure}
\section{History/Motivation}
There are various motivations for the study of holomorphic vector fields. These motivations are largely interconnected, but there seem to be three schools of focus which we outline below. \par
The first school of focus we will mention is that which concentrates on the use of holomorphic vector fields to obtain results within iterative holomorphic dynamics, specifically in studying parabolic bifurcations \cite{Ben1993}, \cite{Shi}, \cite{Outhesis}, \cite{BT2007}, \cite{AEnotes} and in the proof that there exist Julia sets of positive Lebesgue measure \cite{BC2006}. \par
The second school and the one this paper belongs to is studying holomorphic vector fields in their own right (classification) \cite{BT76}, \cite{Sent}, and in the study of quadratic differentials \cite{JAJ},  \cite{Str}. Quadratic differentials and holomorphic vector fields can both be viewed as foliations with singularities.  However, the leaves are canonically oriented in the latter case and not in the first. \par
The last school of focus we will mention is the one that utilizes the integrability of holomorphic vector fields  to study real vector fields and continuous dynamical systems. More specifically, limit cycles and inverse problems relating to Hilbert's $16^{th}$ problem which (among other things) considers the upper bound of the number of limit cycles that can occur for real polynomial vector fields in the plane \cite{Sve1978}, \cite{Sve1981}, \cite{Ben1991}; Hamiltonian systems \cite{Sve1978}; and the existence of centers \cite{Nee1994} are studied. Although no limit cycles exist for holomorphic vector fields, an important strategy when dealing with Hilbert's $16^{th}$ problem is to study perturbations of holomorphic vector fields with centers \cite{AGP09}, \cite{LS04}.\par
In \cite{BT76}, the authors use ``two numerical invariants'' (one being the multiplicity of the zero, and the other coinciding with the \emph{dynamical residue} defined in Section \ref{analinvsec} of this paper) associated with the analytic function defining the vector field, which ``classify the induced flow up to conformal equivalence.'' However, this is a local, not global, result. In this vein, \cite{BT76} and \cite{GGJ2000} also give the local normal forms for such flows. \par
Partial results for the global classification of complex polynomial vector fields go back to the classification of quadratic differentials having poles of order $\geq 2$ \cite{JAJ},  \cite{Str}. The main case where classification of the global structure of quadratic differentials that applies to holomorphic vector fields is only a classification for given quadratic differentials.  When it comes to proving existence of a quadratic differential with prescribed combinatorial and analytic properties (for instance, the so-called ``moduli problem'' for quadradic differentials with closed trajectories  \cite{Str}), the full classification of vector fields is excluded since these theorems are only applicable to a special subclass of vector fields, or they assume a finite area condition, which holomorphic vector fields with at least one zero never satisfy. Muci$\tilde{n}$o-Raymundo \cite{JMR02} studies vector fields but also only covers a specific class. \par
The study of global classification for complex polynomial vector fields in $\C$ is recent \cite{Sent}, and is formulated as follows.\par
The global characterization of a monic, centered complex polynomial vector field corresponding to a fixed polynomial $P$ is determined by a \emph{combinatorial data set} (describing the topology) and an \emph{analytic data set} (describing the geometry).  These invariants uniquely describe such a vector field.  That is, two monic and centered polynomial vector fields with the same invariants must be identical, and  given a set of invariants, there exists a unique monic, centered polynomial vector field  having those invariants.\par
Douady, Estrada, and Sentenac \cite{Sent} classified the \emph{structurally stable}\footnote{Called \emph{generic} in \cite{Sent}} case, i.e. the vector fields  such that there are neither homoclinic separatrices nor \multeq s. It will be shown in \cite{TLD09} that these vector fields in fact make up an open and dense subset of full dimension in parameter space. \par
In this paper, we complete the classification of the global structure of complex polynomial vector fields in $\C$ by extending the result in \cite{Sent} to include the non-structurally stable vector fields.\par
The contents of this paper are as follows.\par
In Section \ref{prelims}, we review the basic definitions and results of monic
centered polynomial vector fields in $\mathbb C$ for polynomials of a
fixed degree $d\geq2$. Sections \ref{eqcl} and \ref{analinvsec} contain the definitions of the
combinatorial and the analytic invariants of such a given vector field
and synthesize their properties so that we are able to abstractly define
combinatorial and analytic data sets.   In Section \ref{floweq} we show that two polynomial vector fields with
the same invariants are identical, and also that any two polynomial
vector fields in the same combinatorial class have quasi-conformally
equivalent flows. In Sections \ref{strthm} - \ref{finalpf}, the invariants are proven to be realizable, formulated in Section \ref{strthm} as the  main theorem, the \emph{Structure Theorem}\footnote{Note that the \emph{Structure Theorem} in Jenkins \cite{JAJ} is a different statement.}, which is proven in several steps in the sections that follow.
 From any given
combinatorial and analytic data set we construct in Section \ref{atlassec} a Riemann
surface $\mathcal M$, called the \emph{rectified surface}. In Section \ref{MisoC} it is
shown to be isomorphic to the Riemann sphere $\hat{\mathbb C}$. The
manifold $\mathcal M$ is obtained from an open non-compact Riemann
surface $\mathcal M^*$ as the compactification of $\mathcal M^*$ by
adding a number of points determined by the combinatorial data
set. By construction, $\mathcal M^*$ has a canonical vector field
$\xi_{\mathcal M}$ assigned. In Section \ref{assvfs}, this vector field is shown to
extend continuously to the zero vector at all the added points,
identifying the equilibrium points of $\xi_{\mathcal M}$. In the final
Section \ref{finalpf}, we show the existence of a unique biholomorphism $\Phi:
\mathcal M \to \hat{\mathbb C}$, which makes the vector field
$\Phi_*(\xi_{\mathcal M})$ a polynomial vector field realizing the given
combinatorial and analytic invariants. The combinatorial and analytic
data sets therefore give a complete global classification of polynomial
vector fields.
\begin{acknowledgements}
We are indebted to Adrien Douady for suggesting the research topic. We would like to thank Douady and likewise Pierette Sentenac for several helpful conversations about their work on the topic. In addition, we want to thank Christian Henriksen for helpful comments and suggestions and  Xavier Buff for the suggestion of using Fatou coordinates in Subsection \ref{vfmult}.
\end{acknowledgements}

\section{Introduction/Preliminaries/Definitions}
\label{prelims}
We recall now some general concepts. Details for much of this section can be found in \cite{Sent} and \cite{BT76}. \par
Given $P$ in $\mathcal{P}_d$, the set of monic and centered polynomials of degree $d \geq 2$, there is an associated vector field $\xi_P \in \Xi_d$ where $\xi_P(z)=P(z) \vf$.  Such a vector field has a corresponding differential equation
\begin{equation}
\label{diffeq}
\dot{z}=P(z),
\end{equation}
where $\dot{z}=\frac{dz}{dt}$, $t \in \R$. The vector field $\xi_P$ has \eqpt s $\zeta$ at the roots of $P$, i.e. $P\left(\zeta\right)=0$.  There are four types of \eqpt s: three types correspond to simple roots of $P$, and one type corresponds to multiple roots of $P$.  
\begin{definition}
An \eqpt\ $\zeta$ is a \emph{sink}, \emph{source}, or \emph{center} if and only if $P\left(\zeta\right)=0$, $P'\left(\zeta\right) \neq 0$, and $\Re \left( P'\left(\zeta\right)\right)$ is negative, positive, or zero respectively.  An \eqpt\ $\zeta$ is a \emph{\multeq} of multiplicity $m$ if and only if $P\left(\zeta\right)=0$, $P'\left(\zeta\right)=\dots =P^{\left(m-1\right)}\left(\zeta\right)=0$, and $P^{\left(m\right)}\left(\zeta\right) \neq 0$.
\end{definition}
The \emph{maximal solution} $\gamma\left(t,z_0\right)$ of Equation \eqref{diffeq} through $z_0$ at $t=0$ satisfies $\gamma\left(0,z_0\right)=z_0$ and $\gamma'\left(t,z_0\right)=P\left(\gamma\left(t,z_0\right)\right)$ for each $t$ in the maximal interval $]t_{min},t_{max}[$ where $t_{min} \in \R_- \cup \{ -\infty \}$ and  $t_{max} \in \R_+ \cup \{ +\infty \}$.  The image $\gamma\left(]t_{min},t_{max}[,z_0 \right)$ is the \emph{trajectory} through $z_0$. If $z_0$ is not an \eqpt\ and if the maximal interval of $\gamma\left(\cdot, z_0\right)$ is mapped bijectively onto the trajectory through $z_0$, then the limit points of the trajectory are as follows. 
\begin{itemize}
\item
For $t_{min}=-\infty$:
$\lim \limits_{t \rightarrow -\infty} \gamma\left(t,z_0\right)=\zeta_{\alpha}$,\newline
a source or a \multeq. 
\item
For $t_{max}=+\infty$:
$\lim \limits_{t \rightarrow +\infty} \gamma\left(t,z_0\right)=\zeta_{\omega}$,\newline
a sink or a \multeq.
\item
For $t_{min}<0$:
$\lim \limits_{t \rightarrow t_{min}} \gamma\left(t,z_0\right)=\infty$, \newline
the point at infinity 
\item
For $t_{max}>0$:
$\lim \limits_{t \rightarrow t_{max}} \gamma\left(t,z_0\right)=\infty,$ \newline
the point at infinity. 
\end{itemize}
If $z_0$ is not an \eqpt\ and if the maximal interval of $\gamma\left(\cdot, z_0\right)$ is not mapped bijectively onto the trajectory through $z_0$, then the maximal interval is $\R$ and $\gamma\left(\cdot,z_0\right)$ is periodic of period $\tau$ where $\tau>0$ is minimal such that $\gamma\left(t+\tau, z_0\right)=\gamma\left(t,z_0\right)$, for all $t \in \R$.  The bounded component of the complement $\C \setminus \gamma\left(\R,z_0\right)$ of the maximal trajectory through $z_0$  contains a center $\zeta$. The period $\tau$ satisfies
\begin{equation}
\tau = \begin{cases}
\ 2 \pi i \frac{1}{P'\left(\zeta\right)} & \text{ if } \Im P'\left(\zeta\right)>0 \\
-2 \pi i \frac{1}{P'\left(\zeta\right)}& \text{ if } \Im P'\left(\zeta\right)<0, 
\end{cases}
\end{equation} 
where in the first case, $\zeta$ is on the left of the periodic trajectory, and in the second case, on the right.\par
The basins of \eqpt s are defined by the following:
\begin{itemize}
\item
$\zeta$ source:\newline
$\bs \left(\zeta\right)=\{ z \in \C \mid \gamma\left(t,z\right) \rightarrow \zeta \text{ for } t\rightarrow -\infty \}$
\item
$\zeta$ sink:\newline
$\bs \left(\zeta\right)=\{ z \in \C \mid \gamma\left(t,z\right) \rightarrow \zeta \text{ for } t\rightarrow +\infty \}$
\item
$\zeta$ center: \newline
$\bs \left(\zeta\right)= \{ \zeta \}\ \cup \ \{ z \in \C \mid \gamma\left(\cdot,z\right) \text{ periodic and } \zeta \text{ is in the bounded}$ \newline $\text{component of } \C \setminus \gamma\left(\R,z\right) \}$
\item
$\zeta$ \multeq :
\begin{align}
\bs\left(\zeta\right)&= \bs_{\alpha}\left(\zeta\right) \cup \bs_{\omega}\left(\zeta\right) \cup \{\zeta \}\quad \text{where} \nonumber \\
\bs_{\alpha}\left(\zeta\right)&=\{ z \neq \zeta \mid  \gamma\left(t,z\right) \rightarrow \zeta \text{ for } t\rightarrow -\infty \} \quad \text{is the repelling basin}\nonumber \\
\bs_{\omega}\left(\zeta\right)&=\{ z \neq \zeta \mid  \gamma\left(t,z\right) \rightarrow \zeta \text{ for } t\rightarrow +\infty \} \quad \text{is the attracting basin}. 
\end{align}
\end{itemize}
The connected components of $\bs_{\alpha}\left(\zeta\right)$ and $\bs_{\omega}\left(\zeta\right)$ are called \emph{repelling petals} and \emph{attracting petals} respectively. In all four cases, $\bs\left(\zeta\right)$ is an open, simply-connected domain. \par
If either $t_{min}$ or $t_{max}$ is finite, then the trajectory associated with this maximal interval is unbounded.  It is  interesting to study the behavior of a polynomial vector field in a neighborhood of infinity. We state a result from \cite{Sent}.
\begin{proposition}[From Chapter I in \cite{Sent}]
\label{2.1}
For every polynomial $P \in \mathcal{P}_d$, there exists a unique isomorphism, tangent to the identity at infinity, which conjugates the vector field $\xi_P$ to $\xi_0\left(z\right)=z^d \vf$ in a neighborhood of $\infty$.
\end{proposition}
Another result from \cite{Sent} states
\begin{proposition}[From Chapter I in \cite{Sent}]
There exist $2d-2$ solutions $\gamma_{\ell}$, $\ell \in \{ 0,1,\dots,2d-3 \}$ such that 
\begin{itemize}
\item
for $\ell$ odd, $\gamma_{\ell}$ is defined for $]-\alpha_{\ell},0]$ and $|\gamma_{\ell}\left(t\right)|\rightarrow \infty$ as $t \rightarrow -\alpha_{\ell}$, and
\item
for $\ell$ even, $\gamma_{\ell}$ is defined for $[0,\beta_{\ell}[$ and $|\gamma_{\ell}\left(t\right)|\rightarrow \infty$ as $t \rightarrow \beta_{\ell}$.
\end{itemize}
In a neighborhood $V_{\infty}$ of infinity, the trajectory $\gamma_{\ell}$ is asymptotic to the ray $t\delta_{\ell}$, $t\in \R_+$ determined by $\delta_{\ell}$ at $\infty$, where
\begin{equation}
\delta_{\ell}=\exp \left(2 \pi i \frac{\ell}{2\left(d-1\right)}\right), \quad \ell \in \{ 0,1,\dots,2d-3 \}.
\end{equation} 
\end{proposition}
\begin{definition}
The \emph{separatrices}  $s_{\ell}$, $\ell=0,\dots,2d-3$, of $\xi_P$ at infinity are the maximal trajectories of $\xi_P$ that have asymptotic directions $\delta_{\ell}$ at infinity. For $\ell$ odd, the separatrix is called \emph{outgoing}, and for $\ell$ even, the separatrix is called \emph{incoming}. 
\end{definition}
\begin{remark}
\label{seplabelremark}
Note the labelling of the separatrices which will be important, among other things, for the uniqueness in the Structure Theorem \ref{fundthm}.
\end{remark}
There are two possibilities for a separatrix $s_{\ell}$.  Either it is \emph{landing}, 
or it is \emph{homoclinic}. 
\begin{definition}
We say $s_{\ell}$ is \emph{landing} or $s_{\ell}$ \emph{lands} at $\zeta$ in $\C$ if  $\zeta=\bar{s}_{\ell} \setminus 
s_{\ell}$ is the limit point of $s_{\ell}$ in $\C$ as $t$ tends to $+\infty$ or $-\infty$, depending on whether $s_{\ell}$ is an outgoing or incoming separatrix to the point at infinity.
\end{definition} 
The limit point $\zeta$ is a sink, source, or \multeq. 
\begin{definition}
If 
$\bar{s}_{\ell} \setminus s_{\ell} = \emptyset$ in $\C$, then the separatrix is both outgoing from and 
incoming to infinity, and is called a \emph{homoclinic separatrix} of infinity.
\end{definition}
   \begin{remark}[Remark and Notation]
When we need to specify that $\ell$ is even or odd, we use $j$ and $k$ 
respectively, while $\ell$ is used to denote either even or odd. 
\end{remark}
Homoclinic separatrices will sometimes be notated by $s_{k,j}$ with two indices, one odd and one even to specify both the outgoing and incoming directions it has at infinity.\par
 We will often consider $\xi_P$ on the Riemann sphere $\rs$, where infinity is a pole of order $d-2$ (see \cite{Sent}).\par
\begin{definition}
The \emph{separatrix graph} of $\xi_P$ is
\begin{equation}
\Gamma_P=\bigcup \limits_{\ell=0}^{2d-3} \hat{s_{\ell}}, 
\end{equation}
where $\hat{s}_{\ell}$ is the closure in $\rs$.
\end{definition} 
\begin{remark}
The separatrix graph $\Gamma_P$ contains infinity and \emph{all} sinks, sources, and \multeq s, but not centers. That is, for every \eqpt \ which is not a center, there exists at least one separatrix which lands at it.
\end{remark}
The connected components $Z$ of $\rs\setminus \Gamma_P$ are called \emph{zones}.  Note that every separatrix is on the boundary of one or two zones (for proof, see \cite{Sent}).  Within the zones, all maximal solutions are defined on all of $\R$.  The following two propositions characterize their types, and it is important to note that the types of zones are determined by the types of their boundaries.
\begin{proposition}
\label{centerzoneprop}
Trajectories in a zone $Z$ that contains an \eqpt\ are periodic, and the \eqpt \ is a center. In this case $Z$ is called a \emph{center zone}. The boundary of a center zone consists of one or several homoclinic separatrices and the point at infinity.
\end{proposition} 
\begin{proof}
See \cite{Sent} and, for instance, \cite{Nee1994}.
\end{proof}
 Otherwise,  a zone does not contain an \eqpt\ in its interior and has at least one \eqpt\  on the boundary $\partial Z$. In such a zone, each solution has an $\alpha$-limit point $\zeta_{\alpha}$ and an $\omega$-limit point $\zeta_{\omega}$ \cite{Sent}.\par
\begin{proposition}[From Chapter I in \cite{Sent}]
\label{zones}
Trajectories in a zone that does not contain \eqpt s have a common $\alpha$-limit point $\zeta_{\alpha}$ and a common $\omega$-limit point $\zeta_{\omega}$. 
From this, one can deduce two types of zones not having an \eqpt\ in its interior.  
\begin{itemize}
\item[1.]
There are exactly two \eqpt s on $\partial Z$, i.e. $\zeta_{\alpha} \neq \zeta_{\omega}$. In this case, $Z$ is called an \emph{$\alpha \omega$-zone} and is of four subtypes:
\begin{itemize}
\item
$Z= \bs\left(\zeta_{\alpha}\right) \cap \bs\left(\zeta_{\omega}\right)$, where $\zeta_{\alpha}$ and $\zeta_{\omega}$ are a source and a sink respectively.
\item
$Z= \bs\left(\zeta_{\alpha}\right) \cap \bs_{\omega}\left(\zeta_{\omega}\right)$, where $\zeta_{\alpha}$ and $\zeta_{\omega}$ are a source and a \multeq\ respectively. In this case, $Z$ is called an \emph{attracting interpetal} for $\zeta_{\omega}$.
\item
$Z= \bs_{\alpha}\left(\zeta_{\alpha}\right) \cap \bs\left(\zeta_{\omega}\right)$, where $\zeta_{\alpha}$ and $\zeta_{\omega}$ are a \multeq\ and a sink respectively. In this case, $Z$ is called a \emph{repelling interpetal} for $\zeta_{\alpha}$.
\item
$Z= \bs_{\alpha}\left(\zeta_{\alpha}\right) \cap \bs_{\omega}\left(\zeta_{\omega}\right)$, where $\zeta_{\alpha}$ and $\zeta_{\omega}$ are both \multeq s. In this case, $Z$ is a repelling interpetal for $\zeta_{\alpha}$ and an attracting interpetal for $\zeta_{\omega}$.
\end{itemize}
The boundary $\partial Z$ contains one or two incoming separatrices and one or two outgoing separatrices, and possibly one or several homoclinic separatrices.
\item[2.]
There is exactly one \eqpt\ on $\partial Z$, i.e. $\zeta_{\alpha} = \zeta_{\omega}$. This case corresponds necessarily to a \multeq :
\begin{itemize}
\item
$Z= \bs_{\alpha}\left(\zeta\right) \cap \bs_{\omega}\left(\zeta\right)$ and is called a \emph{sepal zone} or just \emph{sepal}. There are exactly $2m-2$ sepals corresponding to a \multeq\ of multiplicity $m$ (see \cite{Sent} for details).
\end{itemize}
The boundary $\partial Z$ contains exactly one incoming and one outgoing separatrix and posssibly one or several homoclinic separatrices.
\end{itemize}
\end{proposition}
The point at infinity is a boundary point for each zone.  Within each zone, there are one or several accesses to infinity.  
\begin{definition}
The $2d-2$ \emph{ends at infinity} are the \emph{principal points} of the $\left(2d-2\right)$ prime ends at infinity, defined by the accesses to infinity, where \emph{prime end} is in the sense of Carath\'{e}odory.
\end{definition} 
The ends $e_{\ell}$ are numbered according to their accessibility; the end $e_{\ell}$ has access between the separatrices $s_{\ell-1}$ and $s_{\ell}$. \par
\begin{figure}%
   \centering
   \resizebox{!}{8cm}{\input{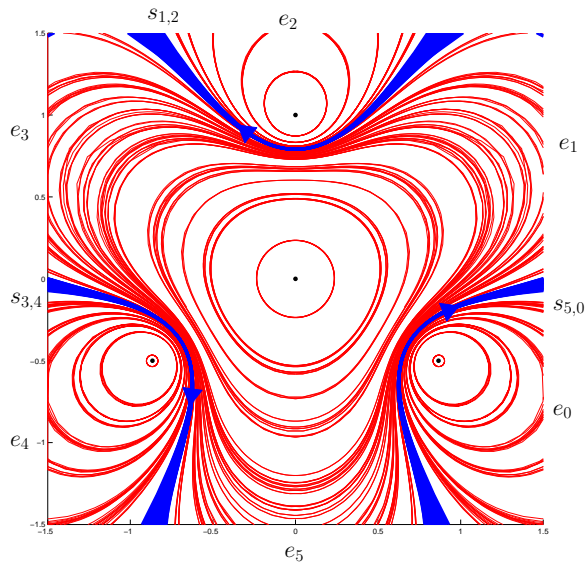}}%
   \caption{An example of the trajectories of a monic and centered polynomial vector field with a counter-clockwise center zone having three homoclinic separatrices $s_{5,0}$, $s_{1,2}$, and  $s_{3,4}$  and three odd ends $e_1$, $e_3$,  and $e_5$ on the boundary.}
  \label{centerzone}
\end{figure}
\begin{figure}%
   \centering
   \resizebox{!}{8cm}{\input{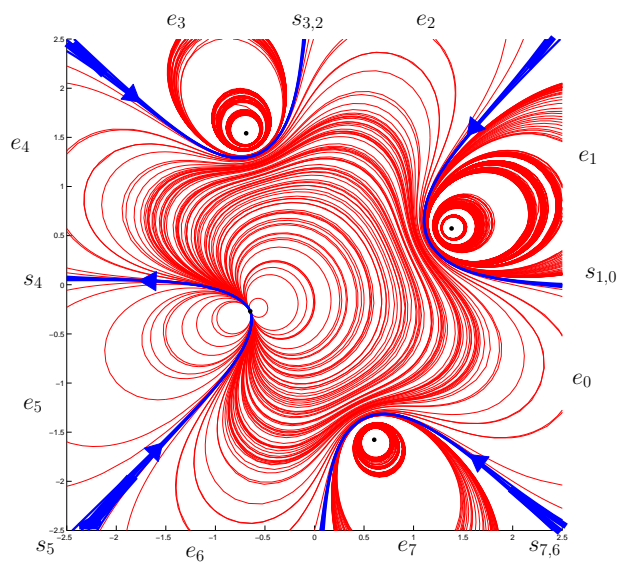}}%
   \caption{An example of the trajectories of a monic and centered polynomial vector field with one odd sepal zone with one end $e_5$ on the boundary and one even sepal zone with three  homoclinic separatrices $s_{1,0}$, $s_{3,2}$, and  $s_{7,6}$  and four even ends $e_0$, $e_2$, $e_4$, and $e_6$ on the boundary.}
  \label{sepalzone}
\end{figure}
\begin{figure}%
   \centering
   \resizebox{!}{8cm}{\input{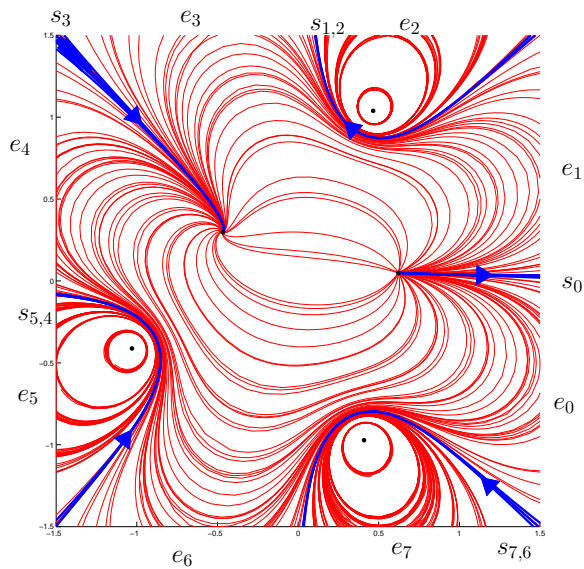}}%
   \caption{An example of the trajectories of a monic and centered polynomial vector field with an $\alpha \omega$-zone $Z^{\alpha \omega}_{3,0}$ having one clockwise homoclinic separatrix $s_{1,2}$, two counter-clockwise homoclinic separatrices $s_{5,4}$ and  $s_{7,6}$, two odd ends $e_1$ and $e_3$, and three even ends $e_0$, $e_4$, and $e_6$ on the boundary.}
  \label{alphomegzone}
\end{figure}

\begin{proposition}[Characterization of the ends on the boundary of zones]
\mbox{}
\begin{itemize}
\item
If a center zone $Z$ has $n$ homoclinic separatrices on the boundary $\partial Z$, then the center zone has $n$ ends at infinity on $\partial Z$.  The ends are all odd if the center zone is to the left of the oriented homoclinic separatrices and all even if the center zone is to the right. In these cases, the zones are called \emph{counter-clockwise} and \emph{clockwise} (or sometimes \emph{odd} and \emph{even}) center zones respectively (see Figure \ref{centerzone}).
\item
If a sepal zone $Z$ has $n$ homoclinic separatrices on the boundary $\partial Z$, then the sepal zone has $n+1$ ends at infinity on $\partial Z$.  They are all odd if the sepal zone is to the left of the oriented homoclinic separatrices and all even if the sepal zone is to the right.  In these cases, the zones are called \emph{odd} and \emph{even} sepal zones respectively (see Figure \ref{sepalzone}).
\item
If an $\alpha \omega$-zone $Z$ is both on the left of $n^+$ homoclinic separartrices on $\partial Z$ and on the right of $n^-$ homoclinic separatrices on $\partial Z$, then the $\alpha \omega$-zone has $n^++1$ odd ends and $n^-+1$ even ends on its boundary $\partial Z$ (see Figure \ref{alphomegzone}). 
\end{itemize}
\end{proposition}
\begin{remark}[Notation]
When we want to distinguish between the different types of zones, we denote an $\alpha \omega$-zone by $Z^{\alpha \omega}$, a center zone by $Z^c$,  and a sepal zone by $Z^s$.
\end{remark}
\subsection{Numbering conventions}
\label{numbconv}
We now want to be able to label the  homoclinic separatrices and zones with an index that contains some combinatorial information.  There are several ways to do this, but we use the following conventions.\par
Let $h_P \geq 0$ be the number of homoclinic separatrices for $\xi_P$, and let $k$ be the uniquely determined odd number so that $s=s_{k,j}$. We number the homoclinic separatrices by $s_{k_1,j_1},\dots,s_{k_{h_P},j_{h_P}}$ such that  
\begin{equation}
0< k_1< \dots < k_{h_P}\leq2d-3.
\end{equation}\par
Let $Z^{\alpha \omega}$ be an $\alpha \omega$-zone for $\xi_P$. From Theorem \ref{zones}, we know that $Z^{\alpha \omega}$ has at least one incoming (even) landing separatrix and at least one outgoing (odd) landing separatrix on $\partial Z^{\alpha \omega}$. 
 We label $Z^{\alpha \omega}_{k,j}$ by the landing separatrix $s_k$ outgoing from the odd end $e_k$ and the landing separatrix $s_j$ incoming to the even end $e_j$ on $\partial Z^{\alpha \omega}$.
Let $s_P\geq 0$ be the number of $\alpha \omega$-zones. We number the zones $Z^{\alpha \omega}_{k_1,j_1},\dots,Z^{\alpha \omega}_{k_{s_P},j_{s_P}}$ corresponding to
\begin{equation}
0\leq j_1< \dots < j_{s_P}<2d-3.
\end{equation}
\par
We label an odd sepal-zone for $\xi_P$ by $Z^s_k$, where $s_k$ is the unique odd landing separatrix on the boundary. Let $n^+_s$ be the number of odd sepal-zones. We number the zones $Z^s_{k_1},\dots,Z^s_{k_{n^+_s}}$ corresponding to
\begin{equation}
0< k_1< \dots < k_{n^+_s}\leq 2d-3.
\end{equation}\par
We label an even sepal-zone for $\xi_P$ by $Z^s_j$, where $s_j$ is the unique even landing separatrix on the boundary.  Let $n^-_s$ be the number of even sepal-zones. We number the zones $Z^s_{j_1},\dots,Z^s_{j_{n^-_s}}$ corresponding to
\begin{equation}
0\leq j_1< \dots < j_{n^-_s}<2d-3.
\end{equation}\par
We label an odd center-zone for $\xi_P$ by $Z^c_k$, where $k$ is the smallest odd index of the homoclinic separatrices on the boundary. Let $n^+_c$ be the number of odd center-zones. We number the zones $Z^c_{k_1},\dots,Z^c_{k_{n^+_c}}$ corresponding to
\begin{equation}
0< k_1< \dots < k_{n^+_c}\leq 2d-3.
\end{equation}\par
We label an even center-zone for $\xi_P$ by $Z^c_j$, where $j$ is the smallest even index of the homoclinic separatrices on the boundary. Let $n^-_c$ be the number of even center-zones. We number the zones $Z^c_{j_1},\dots,Z^c_{j_{n^-_c}}$ corresponding to
\begin{equation}
0\leq j_1< \dots < j_{n^-_c}<2d-3.
\end{equation}\par
\begin{theorem}[from \cite{Sent}]
\label{isos}
For a given vector field $\xi_P \in \Xi_d$, there exist holomorphic isomorphisms $\phi$
\begin{itemize}
\item
from each $\alpha \omega$-zone $Z^{\alpha \omega}_{k,j}$ to a horizontal strip $\Sigma_{k,j}$, where $j$ is the index of the separatrix $s_j$ described above  whose image is on the upper boundary of the strip $\Sigma_{k,j}$ and $k$ is the index of the separatrix $s_k$ described above whose image is on the lower boundary (see Figure \ref{analinv_esstrans}).
\item
from each odd sepal zone $Z^s_k$ to an upper half plane $\HH_k$ and from each even sepal zone $Z^s_j$ to a lower half plane $- \HH_j$,
\item
and from each  center zone minus the curve $\overline{\gamma}_e$,  which is orthogonal to the periodic trajectories in $Z^c$ and joins the center $\zeta$ and $\infty$ through the access defining $e$, to a vertical half-strip $C$. For an odd center zone, $\phi$ maps $Z^c_k \setminus \overline{\gamma}_{e_k}$ to an upper vertical half-strip $C_k=\{ z\in \HH \mid 0<\Re \left(z\right)<\tau\}$, and for an even center zone, $\phi$ maps $Z^c_j \setminus \overline{\gamma}_{e_j}$ to a lower vertical half-strip $C_j=\{ z\in -\HH \mid 0<\Re \left(z\right)<\tau\}$, where $\tau$ is the period of the periodic trajectories,
\end{itemize} 
such that these isomorphisms conjugate $\xi_P$ to $\vf$. 
\end{theorem}
These isomorphisms take the form
\begin{equation}
\label{phi}
\phi\left(z\right)=\int_e^z \frac{dw}{P\left(w\right)},
\end{equation}
where $e$ is one of the ends at infinity on the boundary of the zone. \par
\begin{remark}
\label{ccurves}
The $Z^{\alpha \omega}$, $Z^s$, and  $Z^c \setminus \overline{\gamma}_e$ are simply connected. Therefore, the integral in \eqref{phi} is path independent and the function $\phi$ is well-defined up to addtion by a constant when the end $e$ is changed. The function $\phi$ is actually well-defined in any simply connected domain which avoids the roots of $P$. 
\end{remark}
We then have 
\begin{equation}
\phi_{\ast}\left(\xi_P\right)=\phi'\left(z\right)P\left(z\right)\vf=\vf.
\end{equation}
\begin{definition}
\label{rectcoords}
The isomorphisms $\phi$ are called  \emph{rectifying coordinates} of $\left(\C,\xi_P\right)$.
\end{definition}
We elaborate on the structure of the rectified zones $\pm \HH_{\ell}$, $\Sigma_{k,j}$, and  $C_{\ell}$. 
\begin{itemize}
\item
For an $\alpha \omega$-zone $Z^{\alpha \omega}$, a branch of $\phi$ is chosen such that $e := e_{k_0}$ is one of the odd ends on its boundary.  Thus, $\phi_{k_0}:e_{k_0}\mapsto 0$, and further all of the odd ends on $\partial Z^{\alpha \omega}$ are mapped to $\R$. Any homoclinic separatrices on $\partial Z^{\alpha \omega}$ connecting two odd ends on $\partial Z^{\alpha \omega}$ are also mapped to $\R$.  The one separatrix on  $\partial Z^{\alpha \omega}$ incoming to an odd end and the one separatrix on  $\partial Z^{\alpha \omega}$ outgoing from an odd end are also mapped to $\R$. All even ends, homoclinic separatrices connecting two even ends, and the separatrices incoming to and outgoing from an even end on  $\partial Z^{\alpha \omega}$ have their images on $\R+ih$, where $h$ is the height of the strip.
\item
For an odd sepal zone $Z^s_k$, a branch of $\phi$ is chosen such that $e := e_{k_0}$ is one of the odd ends on its boundary. Thus, $\phi_{k_0}:e_{k_0}\mapsto 0$, $\phi_{k_0}:Z^s_k \rightarrow \HH_{k}$, and further all of the odd ends, all homoclinic separatrices, and the two landing separatrices on $\partial Z^s$ are mapped to $\R$.
\item
For an even sepal zone $Z^s_j$, a branch of $\phi$ is chosen such that $e := e_{j_0}$ is one of the even ends on its boundary. Thus, $\phi_{j_0}:e_{j_0}\mapsto 0$, $\phi_{j_0}:Z^s_j \rightarrow -\HH_{j}$, and further all of the even ends, all homoclinic separatrices, and the two landing separatrices on $\partial Z^s$ are mapped to $\R$. 
\item
For an odd  center zone $Z^c_k$ containing a center $\zeta$, a branch of $\phi$ is chosen such that $e := e_{k_0}$ is one of the odd  ends on its boundary. Thus, $\phi_{k_0}:e_{k_0}\mapsto 0$, $\phi_{k_0}:Z^c\setminus \overline{\gamma}_{e_{k_0}} \rightarrow C_{k_0}$, and further all of the odd ends and all homoclinic separatrices on $\partial Z^c_k$ are mapped to $\R$.
\item
For an even  center zone $Z^c_j$, a branch of $\phi$ is chosen such that $e := e_{j_0}$ is one of the even ends on its boundary. Thus, $\phi_{j_0}:e_{j_0}\mapsto 0$, $\phi_{j_0}:Z^c\setminus \overline{\gamma}_{e_{j_0}} \rightarrow C_{j_0}$, and further all of the even ends and all homoclinic separatrices on $\partial Z^c_j$ are mapped to $\R$.
\end{itemize}
The images of the ends $e_{\ell}$ and separatrices $s_{\ell}$ under $\phi$ are worth distinction and are denoted $E_{\ell}$ and $S_{\ell}$ respectively. \par
To each $\alpha \omega$-zone are associated specific curves joining the odd ends to the even ends.  It will become clear in the next section that these curves are complementary to homoclinic separatrices in more ways than one. The definition follows.
\begin{definition}[Transversals]
A \emph{transversal} $T_{k,j}$ of $\xi_P$ is a curve joining the ends $e_k$ and $e_j$, such that it does not cross any separatrices and crosses the trajectories of $\xi_P$ at a constant, non-zero angle.
\end{definition}
In rectifying coordinates, the images of the transversals are the straight line segments joining the $E_k$ to the $E_j$ in a strip.
  
\section{Combinatorial Data Set for $\xi_P$}
\label{eqcl}
The goal of this section is to define a \emph{combinatorial data set} which completely describes the topological structure of a vector field $\xi_P\in \Xi_d$.\par
The separatrices of $\xi_P$ are labeled by the integers $\{ 0,\dots,2d-3 \}$ or by elements of $\Z /\left(2d-2\right)$ where $0,\dots,2d-3$ are our preferred representatives.\par
The structure of the separatrix graph $\Gamma_P$ is reflected in an \eqre\ $\sim_P$ on $\Z /\left(2d-2\right)$, to be defined below.  Such \eqre s have certain properties that we synthesize into the definiton of a \emph{combinatorial data set of degree $d$} to be defined in Definition \ref{combintdef}. \par
Any $\xi_P \in \Xi_d$ gives rise to a combinatorial data set of degree $d$, and we will prove in Sections \ref{strthm}-\ref{finalpf} that any combinatorial data set is realized by polynomial vector fields from $\Xi_d$.  To a given \cds\ of degree $d$ is therefore associated a \emph{combinatorial class} of polynomial vector fields, i.e. all $\xi_P$ realizing the given \cds. \par
A \cds\ is also called a \emph{combinatorial invariant}.\par
\begin{figure}%
    \centering
    \resizebox{!}{8cm}{\input{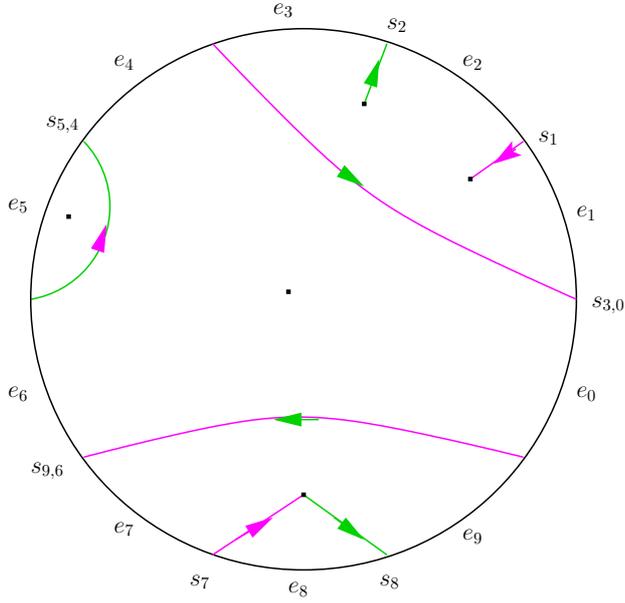}}
    \caption{Example of a disk representation for some $\xi_P \in \Xi_6$ containing one counter-clockwise closed $H_P$-chain $H_{\{[5]\}}$, one clockwise closed $H_P$-chain $H_{\{[0]\}}$, and two counter-clockwise open $H_P$-chains $H_{\{[3]\}}$ and $H_{\{[9]\}}$. }
   \label{closedHchain1}
 \end{figure}

\begin{definition} 
\label{eqredef}
\mbox{}
Let $\xi_P\in \Xi_d$ be given.  Then  the \eqre\ $\sim_P$ and the marked set $H_P$  on $\Z/\left(2d-2\right)$ are defined as follows:
\begin{itemize}
\item[1.] 
$H_P \subset \Z / \left(2d-2\right) $ is the subset defined by 
\begin{equation}
\ell \in H_P \Leftrightarrow s_{\ell} \text{ is a homoclinic separatrix of }\xi_P \nonumber
\end{equation}
\item[2.] The \eqre\ $\sim_P$ on $\Z / \left(2d-2\right) $ is defined by
\begin{align}
&\text{If } \ell',\ell'' \in H_P, \text{ then } \ell' \sim_P \ell'' \Leftrightarrow s_{\ell'}=s_{\ell''}. \nonumber \\
&\text{If } \ell',\ell'' \notin H_P, \text{ then } \ell' \sim_P \ell'' \Leftrightarrow \text{ the separatrices } s_{\ell'} \text{ and } s_{\ell''} \text{ of } \xi_P \text{ land}\nonumber \\
&\text{at the same point in }\C.\nonumber
\end{align}
\end{itemize}
\end{definition}
\begin{remark}[Remark and Notation]
It is convenient to use $L_P$ to name the complement of $H_P$: the set of $\ell$ such that $s_{\ell}$ lands in $\C$. 
\end{remark}
We sometimes represent $L_P$ and $H_P$ by a symbolic disk representation (see Figure \ref{closedHchain1}) where we draw the asymptotic directions of the separatrices on $\mathbb{S}^1$ and draw in curves with arrows representing the separatrices. We also add points corresponding to \eqpt s and sometimes transversals.\par
For $\xi_P$ given and $[\ell] \subseteq L_P$ an \eqcl\ of $\sim_P$, the \eqpt\ which is the common landing point of $s_{\ell'}$ for $\ell' \in [\ell]$ is denoted $\zeta_{[\ell]}\left(\xi_P\right)$ or $\zeta_{[\ell]}$ when the vector field is clear from the context.
\begin{definition}
For $\xi_P$ given, an \eqcl\ $[\ell]\subseteq L_P$ is called \emph{odd}, \emph{even}, or \emph{mixed} if and only if it consists respectively of only odd elements, only even elements, or both odd and even elements (denoted respectively by $[k]$, $[j]$, and $[m]$). An equivalence class $[\ell]\subseteq H_P$ consists of exactly two 
integers, one even and one odd and is called \emph{homoclinic}. 
\end{definition}
\noindent Note that $\sim_P$ has one \eqcl\ if and only if $P\left(z\right)=z^d$.\par
\begin{remark}[Notation]
We use $I=[\ell',\ell'']$ to denote the interval in $\Z/\left(2d-2\right)$ consisting of all the elements $\ell', \ell'+1,\dots,\ell''-1,\ell''$ (i.e. in the counter-clockwise direction). With the preferred representatives $0,\dots,2d-3$, it is possible to have $\ell''<\ell'$, in which case $0\in I$. 
\end{remark}
\begin{definition}
We define the \emph{shift map}
$\sigma: \Z / \left(2d-2\right) \rightarrow \Z / \left(2d-2\right)$ for $H_P \cup L_P$  so that it is a bijection on each \eqcl\ $[\ell]$ and defined by $\sigma\left(\ell\right)$ giving the 
next label in the \eqcl\ in the counter-clockwise direction and $\sigma^{-1}\left(\ell\right)$ giving
the next in the clockwise direction with respect to the asymptotic directions $\delta_{\ell}$.
\end{definition}
The mapping $\sigma$ is an involution on $H_P$, i.e. for $k \in H_P$, $k \sim_P j$, $\sigma\left(k\right)=j$ and $\sigma\left(j\right)=k$.\par
\begin{remark}[Notation]
For an \eqcl\ $[\ell]\subset L_P$, we denote by $p_{[\ell]}$ the number of parity changes
\begin{equation}
p_{[\ell]}=\begin{cases}
0&\text{if } [\ell]=[k] \text{ or } [j]\\
\sum \limits_{\ell'\in [m]}2\left( \left\lfloor \frac{\sigma\left(\ell'\right)-\ell'}{2}\right\rfloor - \frac{\sigma\left(\ell'\right)-\ell'}{2}\right)&\text{if } [\ell]=[m]\end{cases}.
\end{equation}
That is, we count each time $\sigma(\ell')-\ell'$ is odd.  
\end{remark} 
The number $p_{[\ell]}$ of parity changes (corresponding to the number of sepals for $\zeta_{[\ell]}$) is always even, and the number of times there is a change from odd to even is equal to the number of changes from even to odd. \par 
For $\sim_P$, $\sum \limits_{[\ell]\subseteq L_P}p_{[\ell]}$ corresponds to the total number of sepals for the vector field.  
For $\sim_P$ and $[m]$ a mixed \eqcl, the number of interpetals (recall Proposition \ref{zones}) is equal to the number of times the parity does not 
change, i.e. 
when $\sigma\left(\ell\right)-\ell$, $\ell \in [m]$ is even.  The interpetal is attracting when 
two adjacent 
elements in $[m]$ are odd, and the interpetal is repelling when two adjacent 
elements are 
even.
\par
\begin{definition}[Non-crossing \eqre]
\label{noncrossingeqre1}
\mbox{}
An \eqre\ $\sim$ is non-crossing if and only if, for any \eqcl 
es $[\ell]$ with arbitrary $\ell',\ell''\in [\ell]$, any other \eqcl \ $[\tilde{\ell}]$ is contained in either $[\ell'+1,\ell''-1]$ or $[\ell''+1,\ell'-1]$.
\end{definition}
Note that the \eqre\ $\sim_P$ associated to a vector field $\xi_P$ is non-crossing since the separatrices of $\xi_P$ are non-crossing. \par
Important structures induced by the \eqre\ on $H_P$ are \emph{$H_P$-chains}.
\begin{definition}
\label{closedHchain}
A set $\{[k_{i_q}]  \} \subseteq H_P$ of $n\geq 1$ distinct \eqcl es such that $\{i_q \mid q=1,\dots,n  \}\subseteq \{1,\dots , h_P  \}$ 
forms a \emph{counter-clockwise closed $H_P$-chain of length $n$} if
\begin{align}
k_{i_q}&=\sigma\left(k_{i_{q-1}}\right)+1, \quad q=2,\dots,n \nonumber \\
k_{i_1}&=\sigma\left(k_{i_{n}}\right)+1
\end{align}
or forms a \emph{clockwise closed $H_P$-chain of length $n$} if for $j_{i_q}=\sigma\left(k_{i_q}\right)$
\begin{align}
j_{i_q}&=\sigma\left(j_{i_{q-1}}\right)+1, \quad q=2,\dots,n \nonumber \\
j_{i_1}&=\sigma\left(j_{i_{n}}\right)+1
\end{align}
\end{definition}
\begin{definition}
\label{openHchain}
A set $\{[k_{i_q}]  \} \subseteq H_P$ of $n\geq 0$ distinct \eqcl es such that $\{i_q \mid q=1,\dots,n  \}\subseteq \{1,\dots , h_P  \}$ 
forms a \emph{counter-clockwise open $H_P$-chain of length $n$} if for $j_{i_q}=\sigma\left(k_{i_q}\right)$
\begin{align}
k_{i_q}&=\sigma\left(k_{i_{q-1}}\right)+1, \quad q=2,\dots,n \nonumber \\
k_{i_1}-1, & \ \sigma\left(k_{i_{n}}\right)+1 \notin  H_P 
\end{align}
or forms a \emph{clockwise open $H_P$-chain of length $n$} if
\begin{align}
j_{i_q}&=\sigma\left(j_{i_{q-1}}\right)+1, \quad q=2,\dots,n \nonumber \\
j_{i_1}-1, & \ \sigma\left(j_{i_{n}}\right)+1 \notin  H_P.
\end{align}
\end{definition}
Note that any \eqcl\ $\{ k,j\} \subseteq H_P$ is part of a unique counter-clockwise (open or closed) $H_P$-chain and of a unique clockwise (open or closed) $H_P$-chain.\par
The $H_P$-chains are natural structures if one looks at the boundary components of the rectified zones. \par
\begin{figure}%
    \centering
    \resizebox{!}{8cm}{\input{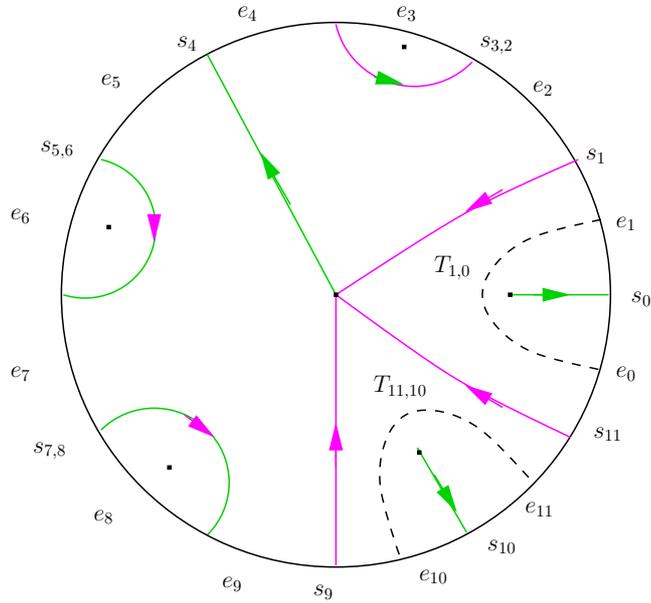}}
    \caption{Example of a disk representation for some $\xi_P \in \Xi_7$ containing one counter-clockwise open $H_p$-chain $H_{\{[5]\}}$, one clockwise open $H_P$-chain $H_{\{[2]\}}$, and two clockwise open $T_P$-chains $T_{\{[1]\}}$ and $T_{\{[11]\}}$. }
   \label{openHTchains1}
 \end{figure}

The homoclinic separatrices corresponding to a closed counter-clockwise (clockwise) $H_P$-chain form the boundary in $\C$ of a counter-clockwise (clockwise) center zone (see Figure \ref{closedHchain1}). In rectifying coordinates, the homoclinic separatrices lie on the lower (upper) boundary of the vertical half-strip. The separatrices corresponding to an open counter-clockwise (clockwise) $H_P$-chain are part of the boundary in $\C$ of an odd (even) sepal zone or an $\alpha \omega$-zone (see Figure \ref{openHTchains1}).  In rectifying coordinates, the separatrices lie on the lower (upper) boundary of the half-plane or the strip.\par
  Let $H_{\{[k_{i_1}]\}}$ denote a counter-clockwise $H_P$-chain and $H_{\{[j_{i_1}]\}}$ be a clockwise  $H_P$-chain, where it will be specified whether it is open or closed when needed. We label the length $n_{\{[\ell_{i_1}]\}}$  by the chain it is associated to, but we sometimes use only $n^+$ or $n^-$ for counter-clockwise and clockwise respectively, or sometimes $n$ for ease in notation when the context is clear.\par
\begin{figure}%
    \centering
    \resizebox{!}{8cm}{\input{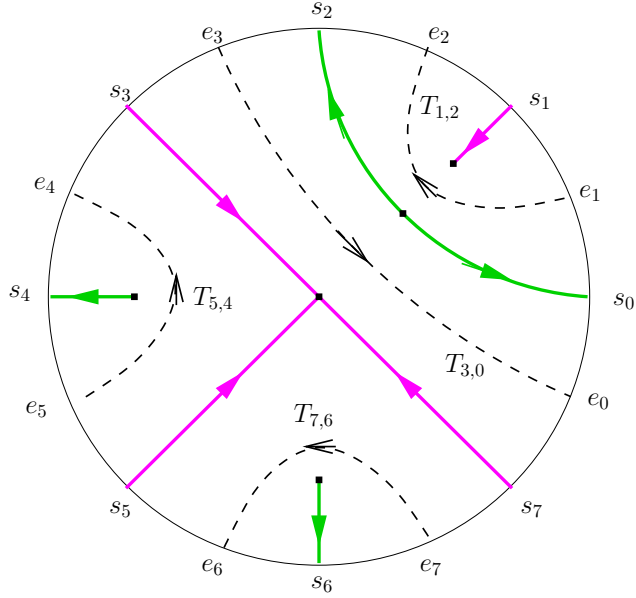}}
    \caption{Example of a disk representation of a vector field $\xi_P$ of degree $d=5$ having neither \multeq s nor homoclinic separatrices. The transversals are all essential and represented by the dashed lines, and their orientations are  denoted by the arrows. Figures \ref{transchain2} and \ref{transchain3} depict the corresponding transversal chains in rectifying coordinates.}
   \label{transchain1}
 \end{figure}

\begin{figure}%
    \centering
    \resizebox{!}{8cm}{\input{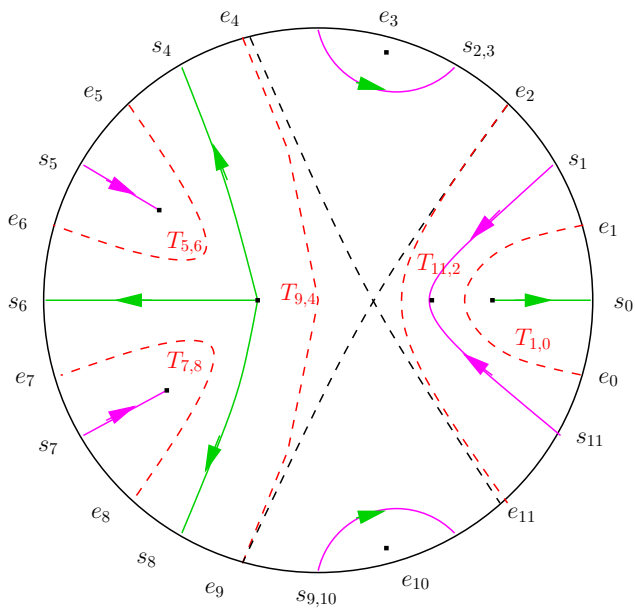}}
    \caption{Example of a disk representation of a $\xi_P \in \Xi_7$ that contains two counter-clockwise closed $T_P$-chains $T_{\{[0]\}}$ and $T_{\{[4]\}}$; three clockwise closed $T_P$-chains $T_{\{[1]\}}$, $T_{\{[5]\}}$, and $T_{\{[7]\}}$; one counter-clockwise open $H_P$-chain $H_{\{[9]  \}}$; and one clockwise open $H_P$-chain $H_{\{[2]  \}}$. The essential transversals are the grey dashed curves  and the non-essential transversals are the black dashed curves.}
   \label{closedTchain1}
 \end{figure}

Complementary to the $H_P$-chains are the structures induced by the transversals. Let 
\begin{equation}
T_P=\{\left(k,\sigma\left(k-1\right)\right)\mid \sigma\left(k-1\right)-(k-1)\in 2\Z\} \cup \{ \left(\sigma\left(j-1\right),j\right)\mid \sigma\left(j-1\right)-(j-1)\in 2\Z  \}.
\end{equation}
\begin{definition}
The \emph{essential transversals} are the transversals $T_{k,j}$ such that $\left(k,j\right) \in T_P$. 
Any other transversal will be called a \emph{non-essential transversal} (see Figure \ref{closedTchain1}). 
\end{definition}
\begin{remark}
Note that in an $\alpha \omega$-zone $Z^{\alpha \omega}_{k,j}$, the essential transversals are the transversals $T_{\sigma^{-1}\left(j\right)+1,j}$ and $T_{k,\sigma^{-1}\left(k\right)}$ (see Figure \ref{closedTchain1}).
Note that $T_{\sigma^{-1}\left(j\right)+1,j}=T_{k,\sigma^{-1}\left(k\right)}$  if and only if there are no homoclinic separatrices on the boundary of $Z^{\alpha \omega}_{k,j}$.
\end{remark}
The structures complementary to the $H_P$-chains are based on these essential transversals and are called \emph{transversal chains} or  \emph{$T_P$-chains}, defined below.
\begin{definition}
A set of $n\geq 1$ distinct transversals $T_{ k_i,j_i }$ such that $\left(k_i,j_i\right)\in T_P$, $i=1,...,n$, form a \emph{counter-clockwise closed $T_P$-chain of length $n$} if
\begin{align}
k_i&=j_{i-1}+1, \quad i=2,...,n \nonumber \\
k_1&=j_{n}+1
\end{align}
or forms a \emph{clockwise closed $T_P$-chain of length $n$} if
\begin{align}
j_i&=k_{i-1}+1, \quad i=2,...,n \nonumber \\
j_1&=k_{n}+1
\end{align}
\end{definition}
\begin{definition}
A set of $n\geq 1$ distinct transversals $T_{ k_i,j_i }$ such that $\left(k_i,j_i\right)\in T_P$, $i=1,...,n$, form a \emph{counter-clockwise open $T_P$-chain of length $n$} if
\begin{align}
k_i&=j_{i-1}+1, \quad i=2,...,n \nonumber \\
(\cdot, k_1-1) \cap T_P=\emptyset,\text{ and } &(j_{n}+1,\cdot)\cap T_P=\emptyset
\end{align}
or forms a \emph{clockwise open $T_P$-chain of length $n$} if
\begin{align}
j_i&=k_{i-1}+1, \quad i=2,...,n \nonumber \\
(j_1-1,\cdot) \cap T_P=\emptyset,\text{ and } &(\cdot,k_{n}+1)\cap T_P=\emptyset
\end{align}
\end{definition}
 If an equilibrium point is contained in a component of $\C \setminus \{T_{k,j}\mid \left(k,j\right)\in T_P  \}$ having a counter-clockwise closed $T_P$-chain on the boundary, it is a source. If an equilibrium point is contained in a component having a clockwise closed $T_P$-chain on the boundary, it is a sink (see Figure \ref{closedTchain1}). So to each closed $T_P$-chain is an associated equivalence class $[\ell]$ corresponding to the \eqpt\ $\zeta_{[\ell]}$, even for counter-clockwise, odd for clockwise.   When a \multeq\ has an interpetal (which is true in all but the case $P(z)=z^d$), then there will be open $T_P$-chains on the boundary. 
Let $T_{\{[k_{i_1}]\}}$ denote a counter-clockwise $T_P$-chain and $T_{\{[j_{i_1}]\}}$ be a clockwise  $T_P$-chain, where it will be specified whether it is open or closed when needed (see Figures \ref{transchain1}, \ref{transchain2}, and \ref{transchain3}). We label the length as for the $H_P$-chains. \par
  \begin{figure}%
    \centering
    \resizebox{!}{6cm}{\input{transchain2.pstex_t}}
    \caption{Depiction of the clockwise transversal chains (oriented line segments) $T_{\{[3]\}}$ and  $T_{\{[1]\}}$ in rectifying coordinates for the combinatorial invariant depicted in Figure \ref{transchain1}. }
   \label{transchain2}
 \end{figure} 
  \begin{figure}%
    \centering
     \resizebox{!}{6cm}{\input{transchain3.pstex_t}}
    \caption{Depiction of the counter-clockwise transversal chains (oriented line segments) $T_{\{[0]\}}$, $T_{\{[4]\}}$,  and  $T_{\{[6]\}}$ in rectifying coordinates for the combinatorial invariant depicted in Figure \ref{transchain1}. }
   \label{transchain3}
 \end{figure} 
For a \multeq , there are some open $H_P$-chains together with some open $T_P$-chains that form the boundary of a domain containing all (and only) the separatrices landing at the \multeq\ (see Figure \ref{openHTchains1}). \par
The set $\Z / \left(2d-2\right)$ has a natural geometric representation as the marked points (later called \emph{division points}) $\delta_{\ell}=\exp \left( \frac{2 \pi i\ell}{2d-2}\right)$ on the unit circle $\mathbb{S}^1$. For any \eqre\ $\sim$ on $\Z / \left(2d-2\right)$ and any \eqcl\ $[\ell]$, we let $[\ell]_{\overline{\D}}$ denote the convex hull of the \eqcl\ $[\ell]$ in the Poincar\'{e} metric in $\overline{\D}$, i.e. the smallest convex closed subset of $\overline{\D}$ containing the geodesics joining $\delta_{\ell'}$ and $\delta_{\ell''}$ for any $\ell',\ell'' \in [\ell]$ (see Figure \ref{diskmodelcells}).\par
In light of this representation in the disk, we give the following alternative definition of a non-crossing \eqre \ (compare with Definition \ref{noncrossingeqre1}).
\begin{definition}[Non-crossing \eqre]
\label{noncrossingeqre}
\mbox{}
An \eqre\ $\sim$ on  $\Z / \left(2d-2\right)$ is called \emph{non-crossing} if and only if for any pair of distinct \eqcl es $[\ell_1]$ and $[\ell_2]$, the corresponding convex hulls $[\ell_1]_{\overline{\D}}$ and $[\ell_2]_{\overline{\D}}$ are disjoint.
\end{definition}
\begin{definition}[The disk model associated to $\sim$]
\mbox{}
Let $\sim$ be any non-crossing \eqre\ on $\Z / \left(2d-2\right)$. The connected components of 
\begin{equation}
\overline{\D} \setminus \bigcup \limits_{\ell \in \Z /\left(2d-2\right)} [\ell]_{\overline{\D}}
\end{equation}
are called \emph{cells}.  The boundary in $\overline{\D}$ of a cell consists of some geodesics connecting $\delta_{\ell'}$ and $\delta_{\ell''}$ for $\ell' \sim \ell''$, some marked points $\delta_{\ell}$, and some arcs in $\mathbb{S}^1$ between marked points. The arc on $\mathbb{S}^1$ between $\delta_{\ell-1}$ and $\delta_{\ell}$ is denoted $\epsilon_{\ell}$ and referred to as an end in the disk model for $\sim$.  Every end is on the boundary of exactly one cell.
\end{definition}
We use a similar convention when drawing the transversals in the disk model.  For a transversal $T_{k,j}$, we join the midpoint of the arc $\epsilon_{k}$ to the midpoint of the arc $\epsilon_{j}$ by the geodesic in the Poincar\'{e} metric.
\begin{remark}
Note that the convex hulls corresponding to the \eqcl\ for $\sim_P$ drawn in the same disk as the transversals according to the above convention are also non-crossing.
\end{remark}
For a vector field $\xi_P \in \Xi_d$, the cells in the disk model are in one-to-one correspondence with the zones of $\xi_P$. We formulate this in the following proposition by giving natural names to the cells and characterizing the different types.
\begin{proposition}
\label{cellchar}
Let $\xi_P \in \Xi_d$ be given and let $\left(\sim_P,H_P\right)$ be defined as in definition \ref{eqredef} The disk model of $\sim_P$ can have up to five different types of cells called an \emph{$\alpha \omega$-cell}, an odd or even \emph{sepal cell}, and an odd or even \emph{center cell} (see Figure \ref{diskmodelcells}). 
 \begin{figure}
  \centering
     \resizebox{!}{8cm}{\input{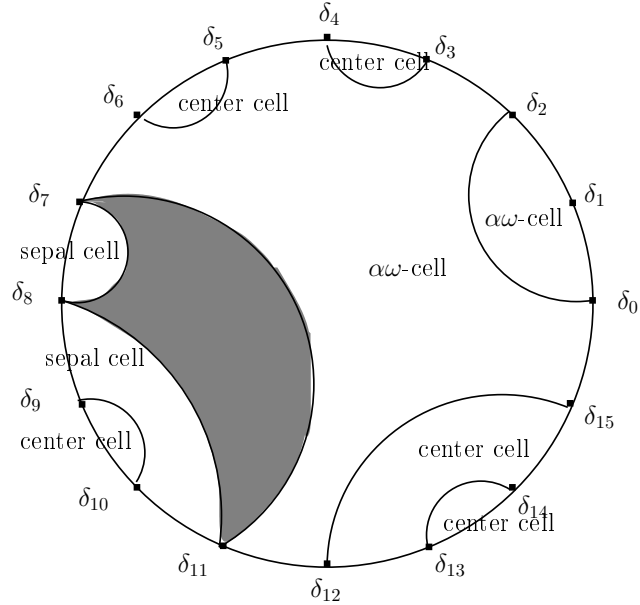}}
    \caption{Example of a disk model for $d=9$ and the different types of cells that can occur for the \eqcl es $H=\{\{3,4\},\{5,6\},\{9,10\},\{12,15\},\{13,14\} \}$, $[0]=\{0,2\}$, $[1]=\{ 1\}$, and $[7]=\{7,8,11\}$. }
   \label{diskmodelcells}
 \end{figure}
The cells are characterized as follows:
\begin{itemize}
\item
A cell is an \emph{$\alpha \omega$-cell} $\Leftrightarrow$ its boundary in $\overline{\D}$ is associated with one counter-clockwise open $H_P$-chain $H_{\{[k_{i_1}]\}}$ of length $n^+ \geq 0$ and the corresponding $k_{i_1}-1,\sigma\left(k_{i_{n^+}}\right)+1 \notin H_P$; one clockwise open $H_P$-chain $H_{\{[j_{i_1}]\}}$ of length $n^- \geq 0$ and the corresponding $j_{i_1}-1,\sigma\left(j_{i_{n^-}}\right)+1 \notin H_P$; the corresponding division points $\delta_{\ell}$, $\ell \in H_{\{[k_{i_1}]\}},\ H_{\{[j_{i_1}]\}}$; $n^++1$ odd ends and $n^-+1$ even ends; and so that $k_{i_1}-1 \sim_P \sigma\left(j_{i_{n^-}}\right)+1$ and $\sigma\left(k_{i_{n^+}}\right)+1 \sim_P j_{i_1}-1$.
\item
A cell is an odd or an even \emph{sepal cell} $\Leftrightarrow$ its boundary in $\overline{\D}$ is associated with one counter-clockwise (resp. one clockwise) open $H_P$-chain $H_{\{[k_{i_1}]\}}$ (resp. $H_{\{[j_{i_1}]\}}$) of length $n \geq 0$ and the corresponding $k_{i_1}-1,\sigma\left(k_{i_{n}}\right) +1 \notin H_P$ (resp.  $j_{i_1}-1,\sigma\left(j_{i_{n}}\right) +1 \notin H_P$ ), the corresponding marked points $\delta_{\ell}$, $\ell\in H_{\{[k_{i_1}]\}}$ (resp. $\ell \in H_{\{[j_{i_1}]\}}$),  $n+1$ odd (resp. even) ends, and so that $k_{i_1}-1\sim_P \sigma\left(k_{i_{n}}\right)+1$  (resp. $j_{i_1}-1 \sim_P \sigma\left(j_{i_{n}}\right)+1)$.
\item
A cell is an odd or an even \emph{center cell} $\Leftrightarrow$ its boundary in $\overline{\D}$ is associated with one counter-clockwise (resp. one clockwise) closed $H_P$-chain $H_{\{[k_{i_1}]\}}$ (resp.  $H_{\{[j_{i_1}]\}}$) of length $n \geq 1$, the corresponding marked points $\delta_{\ell}$, $\ell \in H_{\{[k_{i_1}]\}}$ (resp.  $\ell \in H_{\{[j_{i_1}]\}}$), and the $n$ odd (resp. even) ends. 
\end{itemize}
\end{proposition}
We now have enough notation in order to give the definition of a combinatorial data set of degree $d$.
\begin{definition}
\label{combintdef}
A \emph{combinatorial data set} $\left(\sim,H\right)$ of degree $d \geq 2$ consists of an \eqre\ $\sim$ on $\Z/\left(2d-2\right)$ and a marked subset $H \subset \Z/\left(2d-2\right)$ satisfying:
\begin{itemize}
\item[1)]
$\sim$ is non-crossing.
\item[2)]
If $\ell' \neq \ell''$, then $\ell' \sim \ell''$ and $\ell' \in H \Leftrightarrow \ell'' \in H$ and $\ell'$ and $\ell''$ have different parity. 
\item[3)]
Every cell in the disk-model realization of $\left(\sim,H\right)$ is one of the five types: an $\alpha \omega$-cell, an odd or even sepal-cell, or an odd or even center-cell characterized as above.
\end{itemize}
\end{definition}
\noindent Let $\DD_d$ be the set of combinatiorial data sets of degree $d$. \par
\begin{remark}
If $H = \emptyset$ and there are no mixed \eqcl es so that all cells are $\alpha \omega$-cells, then the definition of a combinatorial data set can be formulated equivalently as a non-crossing pairing of even and odd ends. This is the definition used in \cite{Sent}.
\end{remark}
\begin{remark}[Abstractifications]
\mbox{}
A definition of the different types of cells for $\left(\sim,H\right)\in \mathcal{D}_d$ is taken from Proposition \ref{cellchar} with $\left(\sim_P,H_P\right)$ replaced by $\left(\sim,H\right)$. The numbering conventions of the cells for $\left(\sim,H\right)\in \mathcal{D}_d$ are likewise analagous to those for zones in \ref{numbconv}, with the appropriate replacement of terms.  The properties of the \eqcl es, $H$-chains, and $T$-chains for $\left(\sim,H\right)\in \mathcal{D}_d$ are completely analagous to the properties stated for $\left(\sim_P,H_P\right)$.
\end{remark}
For a given $\left(\sim,H\right)\in \mathcal{D}_d$, we let $s=s\left(\sim,H\right)$ denote the number of $\alpha \omega$-cells in the disk model of $\left(\sim,H\right)$, $p=p\left(\sim,H\right)$ denote the number of sepal-cells,  $h=h\left(\sim,H\right)=\frac{1}{2}|H|$ the number of \eqcl es in $H$, and $c\left(\sim,H\right)=n_c^++n_c^-$ the number of closed $H$-chains.\par
We analyze the \eqre\ $\sim$ for a given $\left(\sim,H\right)\in \mathcal{D}_d$.  \par 
\begin{remark}
\label{homsplit}
For $\left(\sim, H\right)\in \mathcal{D}_d$, if $H \neq \emptyset$, then  $| H|/2=h$ and $\overline{\D}\setminus \bigcup \limits_{[\ell]\subseteq H} [\ell]_{\overline{\D}}$ consists of $h+1$ connected components.  By identifying division points $\delta_j$ and $\delta_k$ where $\{j,k\}=[k] \subseteq H$ and removing them (and hence identifying the ends $\epsilon_k$ with $\epsilon_{j+1}$ or $\epsilon_{k+1}$ with $\epsilon_j$), we obtain $h+1$ circles $\mathbb{S}_i$, $i=0,\dots,h$, each with an even number $2d_i-2$ (possibly 0) of division points adding up to $2d-2-2h$. Let $L_i=\{\ell \in L \mid \delta_{\ell}\in \mathbb{S}_i  \}$ and $[\ell',\ell'']_{L_i}=[\ell',\ell'']\cap L_i$. Each $\overline{\D}_i$ carries an induced \eqre\ $\sim_i$. A connected component in $\overline{\D}\setminus \bigcup \limits_{[\ell]\subseteq H} [\ell]_{\overline{\D}}$ is a \emph{center cell}  if and only if there are no further division points on the boundary of the corresponding  $\overline{\D}_i$, i.e. $d_i=1$, where the $[\ell]\subseteq H$ on the boundary of the component form a closed $H$-chain.
The remaining circles satisfy the following properties outlined in Propositions \ref{numeqcl} and \ref{alphomegeqcl} for a combinatorial data set without homoclinic separatrices (see Figure \ref{decomp1}).
 \begin{figure}
  \centering
     \resizebox{!}{7cm}{\input{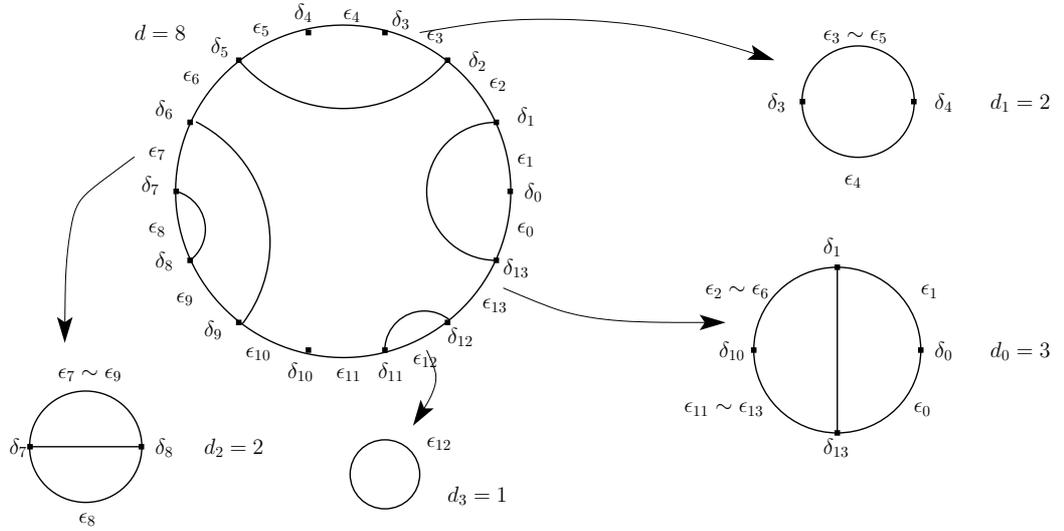}}
    \caption{Example of how a disk model can be decomposed by the \eqcl es in $H=\{\{2,5\},\{6,9\},\{11,12\} \}$ to $h+1$ induced disk models where $H=\emptyset$. }
   \label{decomp1}
 \end{figure}
\end{remark}
\begin{proposition}
\label{numeqcl}
Suppose $\left(\sim,H\right)\in \mathcal{D}_d$ with $H=\emptyset$.
 The equivalence relation $\sim$ on $L=\Z/\left(2d-2\right)$ has 
\begin{equation}
\label{numeqclvf}
q=d- \frac{1}{2}\sum \limits_{[\ell]\subseteq L} p_{[\ell]}
\end{equation}
\eqcl es. 
\end{proposition}
\begin{proof}
There are two possible types of zones: $\alpha \omega$-zones with two ends, and sepal zones with one end each. By comparing the number of ends, we get the relation
\begin{equation}
2\left(d-1\right)=2 s\left(\sim,H\right)+p\left(\sim,H\right).
\end{equation}
 Note that if $s\left(\sim,H\right)=0$, then $\sim$ necessarily has only one \eqcl . Assume that $\sim$ has more than one \eqcl , so that $s\left(\sim,H\right)>0$. The transversals in the $\alpha \omega$-cells divide $\overline{\D}$ into $s\left(\sim,H\right)+1$ subsets, corresponding to the $s\left(\sim,H\right)+1$ \eqcl es. One can see this by noting each subset must correspond to at least one \eqcl \ since there will be at least one division point on the boundary of the subset, and there cannot be more than one \eqcl , since then there would exist a cell with more than two ends on the boundary. Therefore, $\sim$ has 
\begin{equation}
q=s\left(\sim,H\right)+1=d- \frac{1}{2}\sum \limits_{[\ell]\subseteq L} p_{[\ell]},
\end{equation}
\eqcl es since $p\left(\sim,H\right)=\sum \limits_{[\ell]\subseteq L} p_{[\ell]}$.
\end{proof}
\begin{proposition}
\label{alphomegeqcl}
Suppose $\left(\sim,H\right)\in \mathcal{D}_d$ with $H=\emptyset$, and $\sim$ has more than one \eqcl. For each $\alpha \omega$-cell, the associated transversal $T_{k,j}$ divides the disk model into two connected components $U_k$ and $U_j$, numbered so that $\delta_i \in \partial U_i$, $i=k,j$.  Correspondingly, $L=\Z/\left(2d-2\right)$ is divided into two disjoint subsets $I_k$ and $I_j$ where $\ell \in I_i \Leftrightarrow \delta_{\ell} \in \partial U_i$. Define $d_i$ by $|I_i|=\left(2d_i-2\right)+1$, since $|I_i|$ is odd. Then the number of \eqcl es in $I_i$ is
\begin{equation}
q_i=d_i-\frac{1}{2}\sum \limits_{[\ell]\subseteq I_i}p_{[\ell]},\quad i=k,j.
\end{equation}
\end{proposition}
\begin{proof}
Since we assumed $\sim$ has more than one \eqcl, there is always at least one $\alpha \omega$-cell.  Choose an $\alpha \omega$-cell, say the one with transversal $T_{k,j}$. It follows that $k \sim j-1$ and $j \sim k-1$.  
Since $\sim$ is non-crossing, $\sim$ induces an \eqre\  $\sim_k$ on $I_k=[k,j-1]$ and $\sim_j$ on $I_j=[j,k-1]$ by 
\begin{equation}
\label{indeqrel}
\begin{cases}\ell', \ell'' \in I_k & : \ell' \sim_k \ell'' \\
\ell', \ell'' \in I_j & : \ell' \sim_j \ell'' \end{cases} \Leftrightarrow \ell' \sim \ell''.
\end{equation} 
We think of $\delta_k$ as being identified with $\delta_{j-1}$ in $I_k$, and we think of $\delta_j$ as being identified with $\delta_{k-1}$ in $I_j$. In the induced circles $\mathbb{S}_k$ and $\mathbb{S}_j$, the number of ends is even: $j-1-k \pmod{ 2d-2}$ in $\mathbb{S}_k$ and $k-1-j \pmod{ 2d-2}$ in $\mathbb{S}_j$.   \par
By comparing the number of ends, we must have
\begin{align}
j-1-k\ \pmod{2d-2}&=2 \left({}_ks\left(\sim,H\right)\right)+{}_kp\left(\sim,H\right) \nonumber \\
k-1-j\ \pmod{2d-2}&=2 \left({}_js\left(\sim,H\right)\right)+{}_jp\left(\sim,H\right),
\end{align}
where ${}_is\left(\sim,H\right)$ and ${}_ip\left(\sim,H\right)$, are the number of $\alpha \omega$-cells and sepal cells for $I_i$ respectively, $i=k,j$.  The number of \eqcl es for $I_i$ is $q_i={}_is\left(\sim,H\right)+1$, $i=k,j$.  From these relations, one can deduce that the \eqre\  on $I_i$ induces 
\begin{equation}
q_i=d_i-\frac{1}{2}\sum \limits_{[\ell]\subseteq I_i} p_{[\ell]}
\end{equation}
\eqcl es.
\end{proof} 
Below are presented some properties of \eqre s and marked subsets $H$ satisfying 1) and 2) in Definition \ref{combintdef}. These \emph{decomposition properties}  will prove to be equivalent with 3) in Definition \ref{combintdef}.
\begin{definition}[Decomposition Properties]
\label{decompprops}
Let $\sim$ be an \eqre\ on $\Z/\left(2d-2\right)$ and $H \subseteq \Z/\left(2d-2\right)$ a marked subset satisfying 1) and 2) in  Definition \ref{combintdef}. The \eqre \ $\sim$ together with $H$ is said to satisfy the  \emph{decomposition properties} if the following conditions hold.
\begin{itemize}
\item[i.]
Let $h=|H|/2$ (possibly 0). Then in each of the $h+1$ connected components $\overline{\D}_i$ of $\overline{\D}\setminus \bigcup \limits_{[\ell]\subseteq H}[\ell]_{\overline{\D}}$ in the disk model, where we define $d_i$ by setting $2d_i-2$ equal to the number of division points, there are exactly
\begin{equation}
d_i-\frac{1}{2}\sum \limits_{[\ell]\subseteq L_{i}}p_{[\ell]}, \quad i=0,\dots,h
\end{equation}
\eqcl es in that component if $d_i>1$, and no \eqcl es in a component where $d_i=1$.
\item[ii.]
For every partition $I_0=[\ell_0,\ell_1]$ and $I_1=[\ell_1+1,\ell_0-1]$ of the \eqcl es of $\sim$ with $|I_0|$ and $|I_1|$ odd (i.e. $\ell_0$ and $\ell_1$ of the same parity), then $I_0$ has
\begin{equation}
\frac{\ell_1-\ell_0+2}{2}-\frac{1}{2}\sum \limits_{[\ell]\subseteq I_0}p_{[\ell]}
\end{equation}
\eqcl es and $I_1$ has
\begin{equation}
\frac{2d-2+\ell_0-\ell_1}{2}-\frac{1}{2}\sum \limits_{[\ell]\subseteq I_1}p_{[\ell]}
\end{equation}
\eqcl es.
\item[iii.]
For every $\ell \in L_i$, $\sigma\left(\ell\right)-\ell$ is either even or $[\ell+1,\sigma\left(\ell\right)-1]_{L_i}=\emptyset$.
\end{itemize}
\end{definition}
We now show the above properties are equivalent  with $3)$ in Definition \ref{combintdef}, given $1)$ and $2)$.
\begin{theorem}[Characterization of a combinatorial data set $\left(\sim, H\right)\in \mathcal{D}_d$]
\label{selsim}
Let $\sim$ be an \eqre\ and a marked set $H \subseteq \Z/\left(2d-2\right)$ satisfying $1)$ and $2)$ in Definition \ref{combintdef}.  Then  $\left(\sim, H\right)$ satisfies $3)$ in Definition \ref{combintdef} if and only if  $\left(\sim, H\right)$ has the decomposition properties in Definition \ref{decompprops}.
\end{theorem}
We first need a definition and several lemmas to prove this.
\begin{lemma}
\label{sampareq}
\mbox{}
\begin{itemize}
\item[(1)]
Given an \eqre\ $\sim$ on $\Z/\left(2d-2\right)$, $d>2$ with $H=\emptyset$ satisfying property 1) in Definition \ref{combintdef} and the decomposition properties, there exists some distinct $\ell' \sim \ell''$ of the same parity.
\item[(2)]
For $I=[\ell+1,\sigma\left(\ell\right)-1]$, if $\sigma\left(\ell\right)-\ell>2$, then there exists some $\ell' \sim \ell''\in I$ of the same parity.
\end{itemize}
\end{lemma}
\begin{proof}[Proof by contradiction]
Assume not. Then an \eqcl\ can contain at most 2 elements: one even and one odd. Let $n$ be the number of elements in \eqcl es by themselves.
\begin{itemize} 
\item[(1)]
Note that $n$ must necessarily be even. Then the number of \eqcl es is $n+ \frac{2\left(d-1\right)-n}{2}$. Now by the first counting property, the number of \eqcl es should be equal to $d-\frac{1}{2} \sum \limits_{[\ell]\subseteq L}p_{[\ell]}=d-\left( \frac{2\left(d-1\right)-n}{2} \right)=\frac{2+n}{2}$, which implies that $d=2$.
\item[(2)]
Note that $n$ must necessarily be odd. Let $2s+1=\sigma\left(\ell\right)-\ell-1=|I|$. The number of \eqcl es is $n+\frac{2s+1-n}{2}$. By the second decomposition property, the number of \eqcl es should be equal to $s+1-\frac{1}{2}\sum \limits_{[\ell]\subseteq I}p_{[\ell]}=\frac{1+n}{2}$.  This implies that $s=0$ and $\sigma\left(\ell\right)-\ell=2$. 
\end{itemize}
\end{proof}
\begin{lemma}
\label{lsimlp2}
Given an \eqre\ $\sim$ on $\Z/\left(2d-2\right)$, $d>2$, with $H=\emptyset$ satisfying property 1) in Definition \ref{combintdef} and the decomposition properties, there exists some $\ell$ such that $\ell \sim \ell+2$.
\end{lemma}
\begin{proof}
By Lemma \ref{sampareq}, there exists some $\ell,\sigma\left( \ell\right)$ of the same parity with $\sigma\left(\ell\right)\geq \ell+2$. 
If $\sigma\left(\ell\right)= \ell+2$, then we are finished. If $\sigma\left( \ell\right) \geq \ell+4$, then we consider the integers $I=[\ell+1,\sigma\left( \ell\right)-1]$.  Since $I$ is disjoint from $[\sigma\left(\ell\right),\ell]$, we can use the second decomposition property to conclude that they cannot all be in \eqcl es by themselves.  If  all $\ell \in I$ are in the same \eqcl\, then we are again finished. If they are not all in the same \eqcl, then we can  use Lemma \ref{sampareq} (2) to conclude that there must again exist $\ell' \sim \sigma\left(\ell'\right)\in I$ of the same parity with $ \sigma\left(\ell'\right)-\ell' \leq \sigma\left( \ell\right)-\ell-2$.  Repeating this process forms nested subsets of integers, whose maximum difference is strictly decreasing and always even.  Since these are never empty, we must have the claim fulfilled.
\end{proof}
\begin{definition}
Given a non-crossing \eqre\ $\sim$ on $\Z/\left(2d-2\right)$, $d>2$, the \emph{basic partition with basepoint $\ell_0$} is the set of consecutive closed  intervals
\begin{equation}
\left. \begin{array}{cccc}
\underbrace{[\ell_0,\ell_1]}&\underbrace{[\ell_1+1,\ell_2]}&\dots&\underbrace{[\ell_{n-1}+1,\ell_n=\ell_0+2d-3]} \\
I_0 &I_1& &I_{n-1}
\end{array}
\right.
\end{equation}
such that
\begin{align*}
\ell_0&\sim \ell_1, \quad \ell_1 \nsim \ell \in \sum \limits_{i=1}^{n-1}I_i\\
\ell_1+1 &\sim \ell_2, \quad \ell_2 \nsim \ell \in \sum \limits_{i=2}^{n-1}I_i \\
&\vdots \\
\ell_{n-1}+1 &\sim \ell_n=\ell_0+2d-3. 
\end{align*}  
\end{definition}
\begin{lemma}
Given an \eqre\ $\sim$ on $\Z/\left(2d-2\right)$, $d>2$, with $H=\emptyset$ satisfying property 1) in Definition \ref{combintdef} and the decomposition properties,  then for any $\ell_0$, the basic partition with basepoint $\ell_0$ consists of exactly two intervals if $\ell_0 \nsim \ell_0+2d-3$ and one interval if $\ell_0 \sim \ell_0+2d-3$.
\end{lemma}
\begin{proof}
If $\ell_0 \sim \left(\ell_0+2d-3\right)$, then $I_0=[\ell_0,  \ell_0+2d-3]$. Suppose $\ell_0 \nsim \left(\ell_0+2d-3\right)$ and let $I_0,\dots,I_{n-1}$ be the basic partition with basepoint $\ell_0$. Note that $\ell_i+1 \sim \ell_{i+1}$, $i=0,\dots,n-1$, have the same parity by the third decomposition property, and $n\geq 2$. If each $I_i$, $i=0,\dots,n-1$ gives rise to $q_i$ equivalence classes, each having a total number of parity changes
\begin{equation}
\sum \limits_{[\ell]\subseteq I_i}p_{[\ell]},
\end{equation}
then we must have the relation 
\begin{equation}
\label{multmoth}
d-\sum \limits_{i=0}^{k-1} q_i = \sum \limits_{i=0}^{k-1} \frac{1}{2} \sum \limits_{[\ell]\subseteq I_i}p_{[\ell]}.
\end{equation}
We set
\begin{equation}
\label{dast}
 d_i:=\frac{\ell_i-\ell_{i-1}+1}{2}, \quad i=0,\dots,n-1. 
\end{equation}
 Now each $I_i$ must satisfy the decomposition properties, so we further have the relation
\begin{equation}
\label{subprop}
d_i-q_i=\frac{1}{2}\sum \limits_{[\ell]\subseteq I_i}p_{[\ell]},
\end{equation}
and Equation \eqref{dast} gives
\begin{equation}
2\sum \limits_{i=0}^{n-1} d_i=n+\ell_n+1=2d-2+n.
\end{equation}
Now \eqref{subprop} and \eqref{multmoth} gives
\begin{equation}
\sum \limits_{i=0}^{k-1} d_i=d,  
\end{equation}
so 
\begin{equation}
2d=2d-2+n,
\end{equation}
which implies that $n=2$. 
\end{proof}
\begin{lemma}
\label{multhug}
Given an \eqre\ $\sim$ on $\Z/\left(2d-2\right)$, $d>2$, with $H=\emptyset$ satisfying property 1) in Definition \ref{combintdef} and the decomposition properties, if $[\ell]\subseteq [\ell_0,\ell_1]$ and $\ell_0, \ \ell_1 \in [\ell]$, then $\ell_0-1 \sim \ell_1+1$.
\end{lemma}
\begin{proof}
If $\ell_1=\ell_0+2d-3$, this is trivially true. If $\ell_1\neq \ell_0+2d-3$, then $\ell_0$ and $\ell_1$ have the same parity and $I_0=[\ell_0,\ell_1]$ together with $I_1=[\ell_1+1,\ell_0+2d-3]$ must be the basic partition with basepoint $\ell_0$. Hence, $\ell_1 \sim \ell_0+2d-3=\ell_0-1 \left(\mod 2d-2\right)$.
\end{proof}
We are now ready to prove the theorem.
\begin{proof}[Proof of Theorem \ref{selsim}]
\mbox{}
We first show that for $\left(\sim,H\right)\in \mathcal{D}_d$, the  decomposition properties hold.\par
Remark \ref{homsplit} and Proposition \ref{numeqcl} imply property $i.$, and the fact that a sepal cell has only one end on its boundary when $H=\emptyset$ implies property $iii.$.
Proposition \ref{alphomegeqcl} proves property $ii.$ and finishes the first half of the proof.
\par
We now prove that if $\sim$ and $H$ satisfy  1) and 2) from Definition \ref{combintdef} and the decomposition properties, then $\left(\sim,H\right)\in \mathcal{D}_d$ (i.e. 3) from \ref{combintdef} holds). That is, we need to prove that we can only have cells of the types specified in 3). If $H \neq \emptyset$, then the disk model is decomposed into $h+1$ components $\overline{\D}_i$ with associated $d_i$, $i=0,\dots, h$. If $d_i=1$, then the $[\ell]\subseteq H$ on $\partial\overline{\D}_i $ must necessarily form a closed $H$-chain, corresponding to a center cell (see Remark \ref{homsplit}).  It is otherwise enough to consider the case $H=\emptyset$.  An $\alpha \omega$-cell exists when there is one even end $\epsilon_j$ and one odd end $\epsilon_k$ on its boundary.  A sepal-cell exists when there is exactly one end on its boundary.  We prove 3) by induction on $d$. It is true for $d=2$. Indeed, for $d=2$, we have for $H=\emptyset$ either $[0]=\{ 0 \}$ and $[1]=\{ 1\}$ or $[0]= \{0 , 1\} $.   In the first case, $\epsilon_0$ and $\epsilon_1$ are the only ends on the boundary of the same cell (i.e. an $\alpha \omega$-cell). In the second case, $\epsilon_0$ and $ \epsilon_1$ are each the only end on the boundary of different cells (i.e. sepal cells).  Assume now that there are only $\alpha \omega$-cells with one odd and one even end on the boundary and sepal cells with one end on the boundary for all $d$.   Then for 
$d+1$, we have 
that some \eqcl\ either contains only one element or there is an equivalence class having at least three consecutive elements by Lemma \ref{lsimlp2}. If the first case, we can assume without loss of generality that $[0]=\{0\}$.  Then we know by Corollary \ref{multhug} that $1 \sim 2d-3$ and hence $\epsilon_0$ and $\epsilon_1$ are on the boundary of an $\alpha \omega$-cell.  If we identify $\delta_1$ and $\delta_{2d-3}$ and remove $\delta_0$, then we have a degree $d$ \eqre, whose induced equivalence relation on the ends has only $\alpha \omega$-cells and sepal cells  by assumption. Hence,  we again only have cells of the desired types.\par
   If the second case, then we can assume without loss of generality that $\left(2d-3\right) \sim 0 \sim 1$. Then the ends $\epsilon_0$ and $\epsilon_1$ are each on the boundary of a sepal cell. Again, we can identify $\delta_1$ and $\delta_{2d-3}$ and remove $\delta_0$, then we have a degree $d$ \eqre, whose induced equivalence relation on the ends  has only $\alpha \omega$-cells and sepal cells  by assumption. Hence, there must be only cells of the desired type for all $d\geq 2$. 
\end{proof}
  
\section{Analytic Invariants}
\label{analinvsec}
We have defined a combinatorial invariant describing the topological properties of any $\xi_P \in \Xi_d$. Our aim is to define a set of analytic invariants describing the geometric properties.\par
\begin{figure}%
    \centering
    \resizebox{!}{8cm}{\input{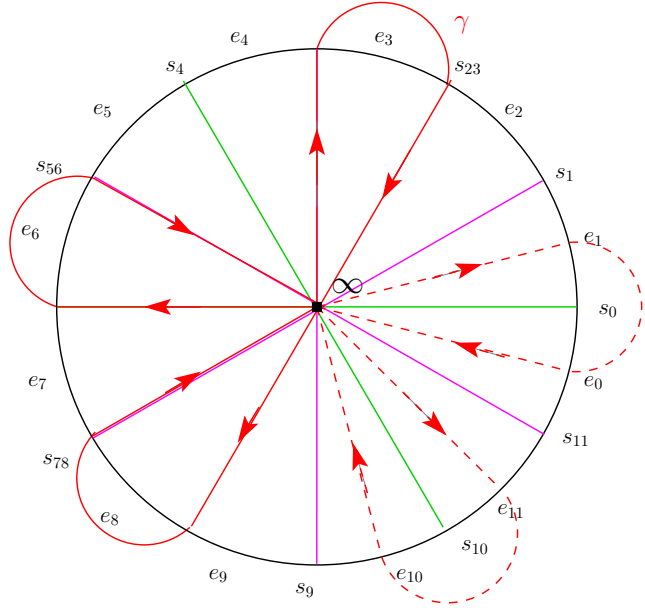}}
    \caption{Depiction of how $\gamma$ might be made from a union of homoclinic separatrices (solid oriented curves) and transversals (dashed oriented curves). Compare with Figure \ref{openHTchains1}. }
   \label{gamma_dynres}
 \end{figure}

It seems natural to consider the \emph{dynamical residues} (defined below) of the \eqpt s as analytic invariants, as these determine the local normal form of the vector field in a neighborhood of each \eqpt\ (see for instance ~\cite{BT76}), but we will see that there is in fact a better choice (see Definition \ref{analinvsdef}).
\begin{definition}[Dynamical residues of $\xi_P$]
Let $\zeta$ be an  arbitrary \eqpt\ of $\xi_P$, i.e. a root of $P$. Let $\gamma$ be a simple, closed, oriented curve not intersecting any \eqpt s and winding counter-clockwise around exactly one \eqpt \ $\zeta$. Then the \emph{dynamical residue} of $\xi_P$ at $\zeta$ is defined as
\begin{equation}
\rho\left(\zeta\right)=\int_{\gamma}\frac{dz}{P\left(z\right)}=2 \pi i \Res \left(\frac{1}{P},\zeta \right).
\end{equation}
\end{definition}
\begin{remark}
Applying the usual formula for residues, we have 
\begin{equation}
\rho\left(\zeta\right)=2 \pi i \lim \limits_{z \rightarrow \zeta}\frac{1}{\left(m-1\right)!}\frac{d^{m-1}}{dz^{m-1}}\left[\left(z-\zeta\right)^m\frac{1}{P\left(z\right)}\right]
\end{equation}
if $\zeta$ is a multiple root of multiplicity $m$. In particular, 
if $\zeta$ is a simple root, then 
\begin{equation}
\rho\left(\zeta\right)=\frac{2 \pi i}{P'\left(\zeta\right)}.
\end{equation}
\end{remark}
Note that we could allow the above curve $\gamma$ to be the piecewise smooth (closed, but not simple) curve consisting of $T_P$-chains and $H_P$-chains bounding the \eqpt\ and the separatrices it receives, together with the point at infinity (see Figure \ref{gamma_dynres}). 
More specifically,
\begin{itemize}
\item
For a center, we let $\gamma$ be the counter-clockwise (or clockwise) $H_P$-chain on the boundary of its zone, together with $\{\infty  \}$.
\item
For a \multeq\ $\zeta_{[m]}$, we let 
\begin{equation}
\gamma=\{ \infty \}\ \cup \bigcup \limits_{\sigma\left(j\right)-j \text{ odd}}H_{\{[j+1]\}}\cup \bigcup \limits_{\sigma\left(k\right)-k \text{ odd}}H_{\{[k+1]\}}\cup \bigcup \limits_{\sigma\left(\ell\right)-\ell \text{ even}}T_{\{[\ell+1]\}},
\end{equation}
for $j,\ k , \ \ell \in [m]$.
\item
For a source or sink $\zeta_{[\ell]}$, we let 
\begin{equation}
\gamma=\{ \infty \} \cup T_{\{[\ell+1]\}}.
\end{equation}
\end{itemize}
The point at infinity does not cause any problems. The singularity for $\frac{dz}{P\left(z\right)}$ at $\infty$ is for $d\geq 2$ a removable singularity, in particular a zero for $d>2$. Indeed, 
\begin{equation}
\left(\frac{-1}{z^2}\right)\left( \frac{1}{P\left(\frac{1}{z} \right)}\right)=\frac{-z^{d-2}}{1+a_{d-2}z^2+\cdots+a_0z^d}.
\end{equation}
\begin{remark}
\label{lincomb}
Since the dynamical residues are simply sums of the integrals of $\frac{dz}{P\left(z\right)}$ over transversals and homoclinic separatrices, it makes sense for these to be the analytic invariants instead. However, if we considered integrals over all transversals, these would contain superfluous information if at least one of the $\alpha \omega$-zones has a homoclinic separatrix on the boundary, since any integral over a transversal is a sum of the integral over some other transversal and integrals over appropriate homoclinic separatrices. Therefore, for each  $\alpha \omega$-zone, the integral over only one transversal is needed to retain the same information. 
\end{remark}
We therefore choose a convention after the numbering of the zones (see Subsection \ref{numbconv}) and state the definition.
\begin{definition}[Definition of analytic invariants of $\xi_P$]
\label{analinvsdef}
\mbox{}
Let $\xi_P \in \Xi_d$ be given.
\begin{itemize}
\item[1.] The \emph{analytic invariant}  of an $\alpha \omega$-zone $Z^{\alpha \omega}_{k,j}$ of $\xi_P$ is defined as 
\begin{equation}
\alpha\left(Z^{\alpha \omega}_{k,j}\right)=\int_{e_k}^{e_j}\frac{dz}{P\left(z\right)},
\end{equation}
where  the integral is along any curve $\gamma_Z$ in $Z^{\alpha \omega}_{k,j}$ connecting $e_k$ to $e_j$ (see Figure \ref{analinv_esstrans}).
  \begin{figure}%
    \centering
    \resizebox{!}{3.5cm}{\input{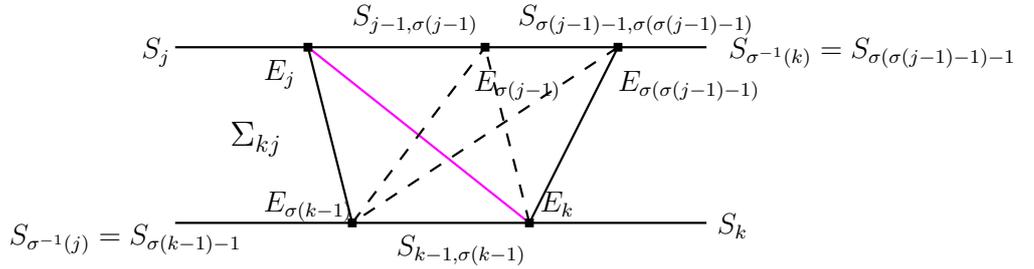}}
    \caption{The strip $\Sigma_{k,j}$ associated to an $\alpha \omega$-zone $Z^{\alpha \omega}_{k,j}$. The transversals are drawn in the strip where the solid black lines are the essential transversals, the solid grey line is the transversal defining the analytic invariant associated to $Z^{\alpha \omega}_{k,j}$, and the dashed lines are some other non-essential transversals. }
   \label{analinv_esstrans}
 \end{figure} 
\item[2.]
The \emph{analytic invariant} of a homoclinic separatrix $s_{k,\sigma\left(k\right)}$ of $\xi_P$ is defined as
\begin{equation}
\tau\left(s_{k,\sigma\left(k\right)}\right)=\int_{s_{k,\sigma\left(k\right)}}\frac{dz}{P\left(z\right)}.
\end{equation}
\item[3.]
The \emph{total analytic invariant} of $\xi_P$ is 
\begin{equation}
A\left(\xi_P\right)=\left(\alpha_1\left(\xi_P\right),\dots,\alpha_{s_P}\left(\xi_P\right),\tau_1\left(\xi_P\right),\dots,\tau_{h_P}\left(\xi_P\right)\right),
\end{equation}
where $\alpha_i\left(\xi_P\right)=\alpha\left(Z^{\alpha \omega}_{k_i,j_i}\right)$, $i=1,\dots,s_P$, and $\tau_i\left(\xi_P\right)=\tau\left(s_{k_i ,\sigma\left(k_i\right)}\right)$, $i=1,\dots,h_P$.
\end{itemize}
\end{definition}
\begin{proposition}[from ~\cite{Sent}]
\label{stripheight}
Let $\xi_P \in \Xi_d$ be given.
\begin{itemize}
\item[1.]
For any $\alpha \omega$-zone $Z$ of $\xi_P$, the invariant $\alpha\left(Z\right)$ belongs to $\HH$ and $\Im \left(\alpha\left(Z\right)\right)$ is equal to the height of the horizontal strip described in Theorem \ref{isos}.
\item[2.]
For any homoclinic separatrix $s$ of $\xi_P$, the invariant $\tau\left(s\right)$ belongs to $\R_+$.
\end{itemize}
\end{proposition}
\begin{remark}
\mbox{}
\begin{itemize}
\item[1.]
It follows from Proposition \ref{stripheight} that $A\left(\xi_P\right)\in \HH^{s_P}\times \R_+^{h_P}$.
\item[2.]
We may refer to the invariant $\alpha\left(Z^{\alpha \omega}_{k,j}\right)$ as the complex ``time'' it takes to go from the end $e_{k}$ to the end $e_j$ in the $\alpha \omega$-zone $Z^{\alpha \omega}_{k,j}$  (see  \cite{Sent} or \cite{BT2007} for elaboration).
\item[3.]
The invariant $\tau\left(s_{k, \sigma\left(k\right)}\right)$ is the (real) time it takes to go along the homoclinic separatrix $s_{k, \sigma\left(k\right)}$ from $\infty$ to $\infty$ in the direction given by $k$ to $\sigma\left(k\right)$.
\end{itemize}
\end{remark}
Later, we will be interested in the sum of all analytic invariants associated to an $H_P$-chain.  Therefore, we now define
\begin{equation}
\tau_{\{[k_{i_1}]\}}\left(\xi_P\right)=\sum \limits_{\nu=1}^{n}\tau_{i_{\nu}}\left(\xi_P\right),
\end{equation}
corresponding to the counter-clockwise $H_P$-chain $H_{\{[k_{i_1}]\}}$ of length $n$ and
\begin{equation}
\tau_{\{[j_{i_1}]\}}\left(\xi_P\right)=\sum \limits_{\nu=1}^{n}\tau_{i_{\nu}}\left(\xi_P\right), 
\end{equation}
corresponding to the clockwise $H_P$-chain $H_{\{[j_{i_1}]\}}$ of length $n$.\par 
We will also often be interested in the integral of $\frac{dz}{P\left(z\right)}$ over an essential transversal which is a linear combination of the analytic invariants associated to a strip. For a strip $\Sigma_{k,j}$, the complex time associated to the essential transversals $T_{\sigma^{-1}\left(j\right)+1,j}$ and $T_{k,\sigma^{-1}\left(k\right)+1}$ is defined as
\begin{align}
\alpha^-_i(j)&=\alpha_i+\tau_{\{[\sigma^{-1}\left(j\right)+1]\}} \quad \text{and} \\
\alpha^+_i(k)&=\alpha_i+\tau_{\{[\sigma^{-1}\left(k\right)+1]\}}
\end{align}
respectively.\par
Applying the Residue Theorem,  we obtain the following relations between the analytic invariants and the dynamical residues.
\begin{itemize}
\item[1.]
For an $\alpha \omega$-zone $Z^{\alpha \omega}_{k,j}$ and for any curve $\gamma_Z$ in $Z^{\alpha \omega}_{k,j}$ connecting $e_k$ to $e_j$,
\begin{equation}
\alpha\left(Z^{\alpha \omega}_{k,j}\right)= \sum \limits_{\zeta \text{ left of } \gamma_Z}\rho\left(\zeta\right).
\end{equation}
\item[2.]
For a homoclinic separatrix $s$:
\begin{equation}
\tau\left(s\right)=\sum \limits_{\zeta \text{ left of } s} \rho\left(\zeta\right).
\end{equation}
\item[3.] For a source $\zeta_{[j]}$:
\begin{equation}
\rho\left(\zeta_{[j]}\right)=\sum \limits_{j'\in [j]}\alpha^-_i(j').
\end{equation}
\item[4.]
For a sink $\zeta_{[k]}$:
\begin{equation}
\rho\left(\zeta_{[k]}\right)=\sum \limits_{k'\in [k]} \alpha^+_i(k')
\end{equation}
(see Remark \ref{lincomb} and Figure \ref{analinv_esstrans}).
\item[5.]
For a center $\zeta$:
\begin{equation}
\rho\left(\zeta\right)=\pm \sum \limits_{s\subset \partial \bs \left(\zeta\right)}\tau\left(s\right)
\end{equation}
with ``$+$'' if $\Im \left(P'\left(\zeta\right)\right)>0$ and  ``$-$'' if $\Im \left(P'\left(\zeta\right)\right)<0$ and where we sum over all homoclinic separatrices on the boundary of $\bs \left(\zeta\right)$, the basin of the center.
\end{itemize}
If there are $q$ equilibrium points which are not centers, having muliplicities 
$m_i$, 
$i=1,\dots,q$, then $s_P$ and $h_P$ satisfy
 \begin{equation}
\label{qshrelation}
 s_P+h_P=d-1-\frac{p_P}{2}.
 \end{equation}
 The number of \eqpt s including centers but not counting multiplicity $|P^{-1}\left(0\right)|=d-\frac{p_P}{2}$ is equal to the number of dynamical residues, so we can conclude from Equation \eqref{qshrelation} that the number of analytic invariants $s+h$ is always one less than the number of dynamical residues. The superfluous information in one of these dynamical residues is due to the fact that by centering, the position of one of the roots is completely determined by the positions of the others.\par
\begin{remark}[Abstractifications]
\mbox{}
The abstractifications of dynamical residue and analytic invariant and the respective numbering of these for $\left(\sim,H\right)\in \mathcal{D}_d$ are analagous to the definitions given for $\xi_P\in \Xi_d$, by replacing $\alpha_i\left(\xi_P\right)$ with $\alpha_i\left(\sim,H\right)$, $\tau_i\left(\xi_P\right)$ with $\tau_i\left(\sim,H\right)$, $H_P$ and $T_P$-chains with $H$ and $T$-chains, and $\rho\left(\zeta_{[\ell]}\right)$ with $\rho_{[\ell]}$.
\end{remark}
The analytic data sets associated to $\left(\sim,H\right)$ are all $(s+h)$-tuples in $\HH^s \times \R_+^h$, i.e.
\begin{equation}
\mathcal{A}\left(\sim,H\right) = \HH^s \times \R_+^h.
\end{equation}
\begin{definition}
For any $\left(\sim,H\right)\in \mathcal{D}_d$ and any $A \in \mathcal{A}\left(\sim,H\right)$, we say that $\xi_P \in \Xi_d$ \emph{realizes} $\left(\sim,H\right)$ and $A$ if and only if $\left(\sim_P,H_P\right)=\left(\sim,H\right)$ (hence $s_P=s \text{ and } h_P=h$) and $A\left(\xi_P\right)=A$.
\end{definition}
\section{Equivalence of Flows within a Combinatorial Class}
\label{floweq}
To every polynomial vector field $\xi_P \in \Xi_d $, we associate a combinatorial invariant and
  analytic invariants.  This section will show that vector fields within a combinatorial class have equivalent flows, and that a fixed analytic invariant within this class can only correspond to one vector field. \par
Let $P_1,\ P_2 \in \mathcal{P}_d$. For simplicity, we use for $i=1,2$ the notation $\xi_{P_i}=\xi_i$, $\sim_{P_i}=\sim_i$, $H_{P_i}=H_i$, $s_{P_i}=s_i$, $h_{P_i}=h_i$, etc.
\begin{theorem}
The vector fields associated with two polynomials $P_1$ and $P_2$ having the same combinatorial data set, i.e. $(\sim_1,H_1)=(\sim_2,H_2)$, have quasi-conformally equivalent flows.
\end{theorem}
\begin{proof}
We will in several steps construct a quasi-conformal homeomorphism $\psi:\C \rightarrow \C$ that makes the flows of $\xi_1$ and $\xi_2$ equivalent. In particular, \eqpt s of $\xi_1$ are mapped to \eqpt s of $\xi_2$ and separatrices to separatrices, but the parameterization by time is not necessarily preserved but is always piecewise linear (to be shown). The mapping $\psi$ is constructed such that the numbering of the separatrices is preserved, in particular $\psi\left(s_0\left(\xi_1\right)\right)=s_0\left(\xi_2\right)$. \par
Since $\left(\sim_1,H_1\right)=\left(\sim_2,H_2\right)$, we have the same \eqcl es in $H:=H_1=H_2$ and $\sim:=\sim_1=\sim_2$, the same number of zones of the different types, and the same open and closed $H$-chains. We will construct
\begin{equation}
\psi:\C \setminus \left( \Gamma_1 \cup \overline{\gamma}_e\left(\xi_1\right) \right) \rightarrow \C \setminus \left( \Gamma_2 \cup \overline{\gamma}_e\left(\xi_2\right) \right),
\end{equation} 
where the $ \overline{\gamma}_e\left(\xi_i\right)$ are the curves removed from the center zones for $\xi_i$ as described in Remark \ref{ccurves} together with the corresponding center.
The isomorphism $\psi$ is constructed by mapping zones of $\xi_1$ to zones of $\xi_2$, respecting the different types and numbering as follows.  \par
For any $\nu=1,\dots,n_c^+$, let 
\begin{equation}
\psi: \left(Z^c_{k_{\nu}}\left(\xi_1\right)\setminus \overline{\gamma}_{e_{k_{\nu}}}\left(\xi_1\right)\right)\rightarrow \left(Z^c_{k_{\nu}}\left(\xi_2\right) \setminus \overline{\gamma}_{e_{k_{\nu}}}\left(\xi_2\right)\right)
\end{equation}
be equal to $\phi_2^{-1} \circ A \circ \phi_1$ where $\phi_1$ and $\phi_2$ are the rectifying coordinates from Proposition \ref{isos} mapping the sliced odd center-zones $\left(Z^c_{k_{\nu}}\left(\xi_i\right) \setminus \overline{\gamma}_{e_{k_{\nu}}}\left(\xi_i\right)\right)$ to the upper vertical half-strips $C_{k_{\nu}}\left( \xi_i\right)$, and let $A:C_{k_{\nu}}\left( \xi_1\right)\rightarrow C_{k_{\nu}}\left( \xi_2\right)$ be a piecewise-affine mapping constructed as follows. If the number of homoclinic separatrices on the boundaries of $C_{k_{\nu}}\left( \xi_i\right)$ is $n$, then for simplicity we renumber the ends on the boundaries of the $C_{k_{\nu}}\left( \xi_i\right)$ from left to right
\begin{equation}
E_{k_{\nu}}:=E_{k_1}=0<E_{k_2}<\cdots<E_{k_n}<E_{k_{n+1}}=\tau_{\{[k_{\nu}]\}},
\end{equation} 
where $E_{k_{n+1}}$ also represents $E_{k_1}$.
We then decompose both $C_{k_{\nu}}\left( \xi_i\right)$ into $n$ vertical half-strips $C^m_{k_{\nu}}\left( \xi_i\right)$, $m=1,\dots ,n$, where 
\begin{equation}
C^m_{k_{\nu}}\left( \xi_i\right)= \left\{ z \in  \HH \mid E_{k_m}\left( \xi_i\right)\leq\Re (z)\leq E_{k_{m+1}}\left( \xi_i\right)\right\}.
\end{equation}
We define $A_m$ on each $C^m_{k_{\nu}}\left( \xi_1\right)$ so that
\begin{equation}
\left[ \begin{array}{c}
\tau\left(S_{k_m,k_{m+1}-1}(\xi_1)\right)  \\
0  \end{array} \right]\mapsto \left[ \begin{array}{c}
\tau\left(S_{k_m,k_{m+1}-1}(\xi_2)\right) \\
0  \end{array} \right] \quad\text{and}\quad \left[ \begin{array}{c}
0  \\
1  \end{array} \right] \mapsto \left[ \begin{array}{c}
0  \\
1  \end{array} \right],
\end{equation}
so that the associated matrix $\underline{A}_m$ becomes
\begin{equation}
\underline{A}_m=\left[ \begin{array}{cc}
\frac{\tau\left(S_{k_m,k_{m+1}-1}(\xi_2)\right)}{\tau\left(S_{k_m,k_{m+1}-1}(\xi_1)\right) } & 0  \\
0 & 1  \end{array} \right],
\end{equation}
and
\begin{equation}
\label{Am}
A_m(x+iy)=\frac{\tau\left(S_{k_m,k_{m+1}-1}(\xi_2)\right)}{\tau\left(S_{k_m,k_{m+1}-1}(\xi_1)\right)}x+iy +\mathcal{T}_m,
\end{equation}
where $\mathcal{T}_m$ is the appropriate positive real constant.
It is similar for the even center zones.\par
For any $\nu=1,\dots,\frac{1}{2}p(\sim,H)$, let 
\begin{equation}
\psi: Z^s_{k_{\nu}}\left(\xi_1\right)\rightarrow Z^s_{k_{\nu}}\left(\xi_2\right)
\end{equation}
be equal to $\phi_2^{-1} \circ A \circ \phi_1$ where $\phi_1$ and $\phi_2$ are the rectifying coordinates from Proposition \ref{isos} mapping the odd sepal zones $ Z^s_{k_{\nu}}\left(\xi_i\right)$ to the upper half-planes $\HH_{k_{\nu}}\left( \xi_i\right)$, and $A:\HH_{k_{\nu}}\left( \xi_1\right)\rightarrow \HH_{k_{\nu}}\left( \xi_2\right)$ is a piecewise-affine mapping constructed as follows. If the number of homoclinic separatrices on the boundaries of $\HH_{k_{\nu}}\left( \xi_i\right)$ is $n$, then for simplicity we renumber the $n+1$ ends on the boundaries of the $\HH_{k_{\nu}}\left( \xi_i\right)$ from left to right
\begin{equation}
E_{k_{\nu}}:=E_{k_1}=0<E_{k_1}<\cdots<E_{k_n}<E_{k_{n+1}}=\tau_{\{[k_{\nu}]\}}.
\end{equation} 
We then decompose both $\HH_{k_{\nu}}\left( \xi_i\right)$ into $n+2$ subsets  $H^m_{k_{\nu}}\left( \xi_i\right)$, $m=1,\dots ,n$, where 
\begin{align}
H^0_{k_{\nu}}\left( \xi_i\right)&=\left\{ z \in \HH \mid -\infty<\Re (z)\leq E_{k_1}\left( \xi_i\right)\right\}\\
H^m_{k_{\nu}}\left( \xi_i\right)&= \left\{ z \in  \HH \mid E_{k_m}\left( \xi_i\right)\leq\Re (z)\leq E_{k_{m+1}}\left( \xi_i\right)  \right\}\\
H^{n+1}_{k_{\nu}}\left( \xi_i\right)&=\left\{ z \in  \HH \mid E_{k_{n+1}}\left( \xi_i\right)\leq\Re (z)<\infty \right\},
\end{align}
We define $A_m$ on each $H^m_{k_{\nu}}\left( \xi_1\right)$ to be
\begin{equation}
A(x+iy)=\begin{cases} 
\id&\text{for } m= 0\\ 
A_m(x+iy) \text{ from Equation }\eqref{Am} &\text{for } m=1,\dots,n\\
\id + \tau_{\{[k_{\nu}]\}}(\xi_{2})- \tau_{\{[k_{\nu}]\}}(\xi_1)& \text{for }m=n+1\end{cases}.
\end{equation}
It is similar for the even sepal zones.\par
For any $\nu=1,\dots,s(\sim,H)$, let
\begin{equation}
\psi: Z^{\alpha \omega}_{k_{\nu},j_{\nu}}\left(\xi_1\right)\rightarrow Z^{\alpha \omega}_{k_{\nu},j_{\nu}}\left(\xi_2\right) 
\end{equation}
be equal to $\phi_2^{-1} \circ A \circ \phi_1$ where $\phi_1$ and $\phi_2$ are the rectifying coordinates from Proposition \ref{isos} mapping the $\alpha \omega$-zones $Z^{\alpha \omega}_{k_{\nu},j_{\nu}}\left(\xi_i\right)$ to the horizontal strips $\Sigma_{k_{\nu},j_{\nu}}\left(\xi_i\right)$, and $A:\Sigma_{k_{\nu},j_{\nu}}\left( \xi_1\right)\rightarrow \Sigma_{k_{\nu},j_{\nu}}\left( \xi_2\right)$ is a piecewise-affine mapping constructed as follows. 
If the number of homoclinic separatrices on the upper boundary of the $\Sigma_{k_{\nu},j_{\nu}}\left(\xi_i\right)$ is $n^-$ and the number of homoclinic separatrices on the lower boundary of the $\Sigma_{k_{\nu},j_{\nu}}\left(\xi_i\right)$  is $n^+$, then we set $n=n^++n^-$ and decompose the $\Sigma_{k_{\nu},j_{\nu}}\left(\xi_i\right)$ into $n+2$ subsets $\Sigma^m_{k_{\nu},j_{\nu}}\left(\xi_i\right)$, $m=0,\dots ,n+1$, where $\Sigma^0_{k_{\nu},j_{\nu}}$ is the left most part of the strip bounded by the separatrices $S_{j_{\nu}}$ and $S_{\sigma^{-1}(j_{\nu})}$ and the transversal $T_{\sigma^{-1}(j_{\nu})+1,j_{\nu}}$, and $\Sigma^{n+1}_{k_{\nu},j_{\nu}}\left(\xi_i\right)$ is the right most part of the strip bounded by the separatrices $S_{k_{\nu}}$ and $S_{\sigma^{-1}(k_{\nu})}$ and the transversal $T_{\sigma^{-1}(k_{\nu})+1,k_{\nu}}$. The rest of the $\Sigma_{k_{\nu},j_{\nu}}\left(\xi_i\right)$ is divided  into triangles. The triangles have 3 $E_{\ell}$s as vertices (at least one even and odd) and two transversals and one homoclinic separatrix as edges. There are several ways to do this, but the important thing to demand is that the triangles must be disjoint (except for possibly shared edges) and cover the area. There can be two types of triangles: one type with an edge contained in a clockwise $H$-chain and one type with an edge contained in a counter-clockwise $H$-chain (see Figure \ref{Striang}). 
 \begin{figure}%
   \centering
   \resizebox{!}{3.5cm}{\input{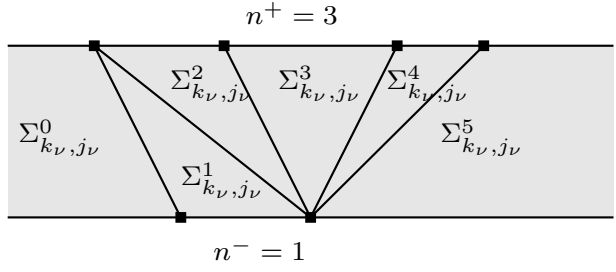}}%
   \caption{Depiction of a possible decomposition of $\Sigma_{k_{\nu},j_{\nu}}$ into $n+2=6$ subsets $\Sigma_{k_{\nu},j_{\nu}}^m$, $m=0,\dots,5$. There is one triangle $\Sigma_{k_{\nu},j_{\nu}}^1$ that has an edge contained in the counter-clockwise $H$-chain on the boundary, and there are three triangles $\Sigma_{k_{\nu},j_{\nu}}^m$, $m=2,\dots,4$ having an edge in the clockwise $H$-chain on the boundary.}
  \label{Striang}
\end{figure} 
We number these triangles from left to right as $\Sigma^m_{k_{\nu},j_{\nu}}\left(\xi_i\right)$, $m=1,\dots ,n$. For simplicity of notation, let $\alpha_1$ be the analytic invariant associated to $Z^{\alpha \omega}_{k_{\nu},j_{\nu}}\left(\xi_1\right)$ and $\alpha_2$ the analytic invariant associated to $Z^{\alpha \omega}_{k_{\nu},j_{\nu}}\left(\xi_2\right)$.
We define $A_0$ on $\Sigma^0_{k_{\nu},j_{\nu}}\left(\xi_1\right)$ so that 
\begin{equation}
\left[ \begin{array}{c}
1 \\
0  \end{array} \right]\mapsto \left[ \begin{array}{c}
1  \\
0  \end{array} \right] \quad\text{and}\quad \left[ \begin{array}{c}
\Re (\alpha^-_1)  \\
\Im(\alpha^-_1)  \end{array} \right] \mapsto \left[ \begin{array}{c}
\Re (\alpha^-_2)  \\
 \Im(\alpha^-_2)  \end{array} \right],
\end{equation}
and we define $A_{n+1}$ on $\Sigma^{n+1}_{k_{\nu},j_{\nu}}\left(\xi_1\right)$ so that 
\begin{equation}
\left[ \begin{array}{c}
1 \\
0  \end{array} \right]\mapsto \left[ \begin{array}{c}
1  \\
0  \end{array} \right] \quad\text{and}\quad \left[ \begin{array}{c}
\Re (\alpha^+_1)  \\
\Im(\alpha^+_1)  \end{array} \right] \mapsto \left[ \begin{array}{c}
\Re (\alpha^+_2)  \\
 \Im(\alpha^+_2)  \end{array} \right],
\end{equation}
so that
\begin{equation}
\underline{A}=\left[ \begin{array}{cc}
1&\frac{\Re(\alpha^{\pm}_2-\alpha^{\pm}_1)}{\Im(\alpha^{\pm}_1)} \\
0 & \frac{\Im(\alpha^{\pm}_2)}{\Im(\alpha^{\pm}_1)}\end{array} \right]
\end{equation}
and
\begin{align}
A_0(x+iy)&=x+\frac{\Re(\alpha^-_2-\alpha^-_1)}{\Im(\alpha^-_1)}y+i \frac{\Im(\alpha^-_2)}{\Im(\alpha^-_1)}y\\
A_{n+1}(x+iy)&=x+\frac{\Re(\alpha^+_2-\alpha^+_1)}{\Im(\alpha^+_1)}y+i \frac{\Im(\alpha^+_2)}{\Im(\alpha^+_1)}y+\mathcal{T},
\end{align}
where $\mathcal{T}=\tau_{\{[\sigma^{-1}(j_{\nu})+1]\}}\left(\xi_2\right)-\tau_{\{[\sigma^{-1}(j_{\nu})+1]\}}\left(\xi_1\right)$ is the appropriate positive real constant.
Now let 
\begin{equation}
v_{\tilde{k}\tilde{j}}\left(\xi_i\right)=\alpha^-_i+\sum\limits_{j'=\tilde{j}}^{\sigma(j-1)}\tau\left(S_{j',\sigma(j')}\left(\xi_i\right)\right)-\sum\limits_{k'=\sigma^{-1}(j)+1}^{\sigma(\tilde{k}-1)}\tau\left(S_{k',\sigma(k')}\left(\xi_i\right)\right).
\end{equation}
Then on each triangle $\Sigma^m_{k_{\nu},j_{\nu}}\left(\xi_1\right)$, $m=1,\dots,n$, with $v_{\tilde{k}\tilde{j}}(\xi_1)$ as the left edge and $S_{\tilde{j}-1,\sigma(\tilde{j}-1)}(\xi_1)$ as the upper (or $S_{\tilde{k},\sigma(\tilde{k})}(\xi_1)$ as the lower) edge, we define $A_m$ by
\begin{equation}
\left[ \begin{array}{c}
\tau\left(S_{\ell,\sigma(\ell)}\left(\xi_1\right)\right) \\
0  \end{array} \right]\mapsto \left[ \begin{array}{c}
 \tau\left(S_{\ell,\sigma(\ell)}\left(\xi_2\right)\right)  \\
0  \end{array} \right] \quad\text{and}\quad \left[ \begin{array}{c}
\Re (v_{\tilde{k}\tilde{j}}\left(\xi_1\right))  \\
\Im(v_{\tilde{k}\tilde{j}}\left(\xi_1\right))  \end{array} \right] \mapsto \left[ \begin{array}{c}
\Re (v_{\tilde{k}\tilde{j}}\left(\xi_2\right))  \\
 \Im(v_{\tilde{k}\tilde{j}}\left(\xi_2\right))  \end{array} \right], 
\end{equation}
$ \ell=\tilde{k}, \tilde{j}-1$, so that 
\begin{equation}
\underline{A}=\left[ \begin{array}{cc}
\frac{ \tau\left(S_{\ell,\sigma(\ell)}\left(\xi_2\right)\right) }{\tau\left(S_{\ell,\sigma(\ell)}\left(\xi_1\right)\right)}&\frac{\Re(v_{\tilde{k}\tilde{j}}\left(\xi_1\right)-v_{\tilde{k}\tilde{j}}\left(\xi_1\right))}{\Im(v_{\tilde{k}\tilde{j}}\left(\xi_1\right))} \\
0 & \frac{\Im(v_{\tilde{k}\tilde{j}}\left(\xi_2\right))}{\Im(v_{\tilde{k}\tilde{j}}\left(\xi_1\right))}\end{array} \right]
\end{equation}
and
\begin{equation}
A_{m}(x+iy)=x+\frac{\Re(v_{\tilde{k}\tilde{j}}\left(\xi_2\right)-v_{\tilde{k}\tilde{j}}\left(\xi_1\right))}{\Im(v_{\tilde{k}\tilde{j}}\left(\xi_1\right))}y+i \frac{\Im(v_{\tilde{k}\tilde{j}}\left(\xi_2\right))}{\Im(v_{\tilde{k}\tilde{j}}\left(\xi_1\right))}y+\mathcal{T}_m, \quad m=1,\dots,n,
\end{equation}
where $\mathcal{T}_m$ is the proper translation as above.\par
Every above map is an orientation-preserving linear map and is hence a quasi-conformal homeomorphism. Since each $\phi_1$, $\phi_2$ is a conformal isomorphism, then the compositions 
\begin{equation}
\phi_2^{-1}\circ A \circ \phi_1
\end{equation}
are quasi-conformal on the $Z^c\setminus \{ \zeta_c\}$, $Z^s$, and $Z^{\alpha \omega}$. These maps  extend continuously to the \eqpt s $\zeta\left( \xi_i \right)$. 
\end{proof}
\begin{corollary}
Suppose $\xi_1,\ \xi_2 \in \Xi_d$ have the same combinatorial and analytic invariants, i.e. $\left(\sim_1,H_1\right)=\left(\sim_2,H_2\right)$ and $A\left(\xi_1\right)=A\left(\xi_2\right)$. Then $\xi_1=\xi_2$ (i.e. $P_1=P_2$).
\end{corollary}
\begin{proof}
When the analytic invariants are identical, each $A=\id$ above, and in this case we can conclude that $\psi$ is a holomorphic automorphism of $\C$.
The only isomorhpisms $\psi$ that conjugate vector fields in $\Xi_d$ are of the form $z \mapsto \sqrt[d-1]{1}z$. Since both $s_0\left(\xi_1\right)$ and $s_0\left(\xi_2\right)$ are asymptotic to $\R_+$ at infinity, it follows that $\psi$ is the identity.
\end{proof}
Therefore, the combinatorial and analytic invariants assigned to a vector field 
determine the vector field uniquely in $\Xi_d$. 

\section{The Structure Theorem}
\label{strthm}
We have defined analytic and combinatorial invariants, given a vector field.  We now wish to show that given a combinatorial invariant and a set of 
  analytic invariants with the properties described above,  there exists a
 unique polynomial vector field $\xi_P \in \Xi_d$ realizing the given invariants.
\begin{theorem}[Structure Theorem]
Given a combinatorial data set $\left(\sim, H\right) \in \DD_d$ and an analytic data set $A \in \mathcal{A}\left(\sim, H\right)$, there exists a unique $\xi_P \in \Xi_d$ realizing $\left(\sim, H\right)$ and $A$, i.e. $\left(\sim_P,H_P\right)=\left(\sim, H\right)$ and $A\left(\xi_P\right)=A$.
 \end{theorem}
The proof is constructed as follows:      
\begin{itemize}
\item[1.] Construct the \emph{rectified surface} $\M$ with vector field $\xi_{\M}$ from the data $\left(\sim,H\right)$ and $A$.
\item[2.] Define an atlas on $\M$ to show it is a
Riemann surface
\item[3.] Show that $\M$ is isomorphic to $\rs$.
\item[4.] Show existence of the unique isomorphism $\Phi:\M \rightarrow \rs$ that induces 
$\Phi_{\ast}\left(\xi_{\M}\right)=\xi_P \in \Xi_d$ having the given invariants.
\end{itemize} 
\section{Definition of and Atlas for $\M$}
\label{atlassec}
Let $(\sim,H)$ (a combinatorial data set) and $\mathcal{A}(\sim,H)$ (a corresponding total analytic invariant) be given.  Out of these data, we will construct a Riemann surface $\M$ with a vector field $\xi_{\M}$. \par
\subsection{Construction of the Rectified Space}
We first construct the rectified zones:  vertical half-strips,  half-planes, and strips  as follows.\par
\begin{figure}%
    \centering
    \resizebox{!}{3cm}{\input{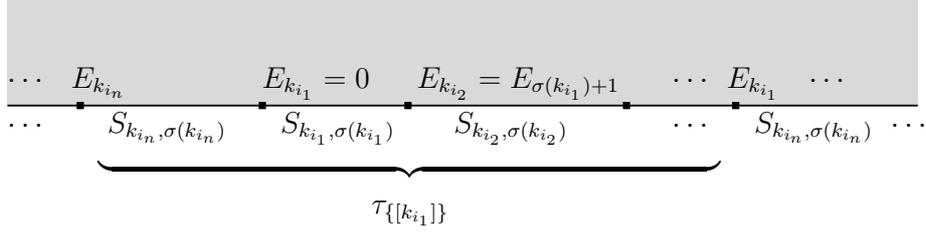}}
    \caption{Depiction of how to build an upper vertical half-strip from a counter-clockwise closed $H$-chain. }
   \label{buildccwcyl}
 \end{figure}
\begin{figure}%
    \centering
    \resizebox{!}{3cm}{\input{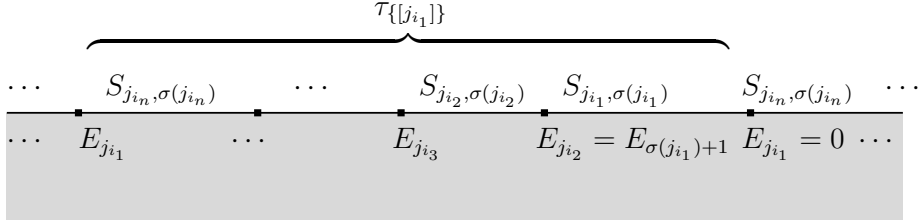}}
    \caption{Depiction of how to build a lower vertical half-strip from a clockwise closed $H$-chain.  }
   \label{buildcwcyl}
 \end{figure}

For every counter-clockwise (respectively clockwise) closed $H$-chain $H_{\{[k_{i_1}]  \}}$ (resp. $H_{\{[j_{i_1}]  \}}$), we construct an upper (lower) vertical half-strip 
\begin{equation} 
C_{\{[k_{i_1}]\}}=\left\{z\in \HH \mid 0<\Re (z)< \tau_{\{[k_{i_1}]\}} \right\}\nonumber
\end{equation}
\begin{equation}
\left(\text{resp. } C_{\{[j_{i_1}]\}}=\left\{ z\in -\HH \mid  0<\Re (z)< \tau_{\{[j_{i_1}]\}} \right\} \right) \nonumber
\end{equation}
Its lower (resp. upper) boundary is $]0,\tau_{\{[k_{i_1}]\}}[\subset \R$ (resp. $]0,\tau_{\{[j_{i_1}]\}}[\subset \R$), labeled according to the corresponding $H$-chain as in Definition \ref{closedHchain}.  That is, 
 $]0,\tau(S_{k_{i_1},\sigma(k_{i_1})})[$ (resp. $]0,\tau(S_{j_{i_n},\sigma(j_{i_n})})[$) is labeled by $S_{k_{i_1},\sigma(k_{i_1})}$ (resp. $S_{j_{i_n},\sigma(j_{i_n})}$), and 
\begin{equation}
 ]\tau(S_{k_{i_1},\sigma(k_{i_1})}),\tau(S_{k_{i_1},\sigma(k_{i_1})})+\tau(S_{k_{i_2},\sigma(k_{i_2})})[  \nonumber 
\end{equation}
\begin{equation}
\left(\text{resp. }  ]\tau(S_{j_{i_n},\sigma(j_{i_n})}),\tau(S_{j_{i_n},\sigma(j_{i_n})})+\tau(S_{j_{i_{n-1}},\sigma(j_{i_{n-1}})})[\right)\nonumber
\end{equation}
is labeled as in Definition \ref{openHchain} by $S_{k_{i_2},\sigma(k_{i_2})}$ (resp. $S_{j_{i_{n-1}},\sigma(j_{i_{n-1}})}$), etc. (see Figure \ref{buildccwcyl} resp. Figure \ref{buildcwcyl}). 
Each homoclinic separatrix $S_{k,j}$ in the $H$-chain connects the  marked points $E_k$ and $E_{j+1}$ (resp. $E_{k-1}$ and $E_j$) in $\R$.   \par
For every counter-clockwise (resp. clockwise) open $H$-chain $H_{\{[k_{i_1}]\}}$ (resp. $H_{\{[j_{i_1}]\}}$) of length $n\geq 0$, the odd separatrix $S_k$, $k=\sigma(k_{i_n})+1$ (resp. $k=j_{i_n}-1$),  and the even separatrix $S_j$, $j=k_{i_1}-1$ (resp. $j=\sigma(j_{i_1})+1$), are attached to the $H$-chain by the ends $E_k$ and $E_{j+1}$ (resp. $E_{k+1}$ and $E_j$) respectively.  If $k=\sigma(j)$ (resp. $j=\sigma(k)$), then there are two possibilities: either the $H$-chain together with  $S_k$, $S_j$, and all the ends they connect make up the lower (resp. upper) boundary of an upper (resp. lower) half-plane $\HH_k$ (resp. $-\HH_j$), or they are on the lower (resp. upper) boundary of the strip $\Sigma_{k,\sigma(j)}$ (resp. $\Sigma_{\sigma(k),j}$), where $S_{\sigma(j)}$ and $S_{\sigma^{-1}(k)}$ (resp. $S_{\sigma^{-1}(j)}$ and $S_{\sigma(k)}$) are attached to the same open $H$-chain, and these together with the ends they connect make up the upper (resp. lower) boundary of the strip. \par
\begin{figure}%
    \centering
    \resizebox{!}{4cm}{\input{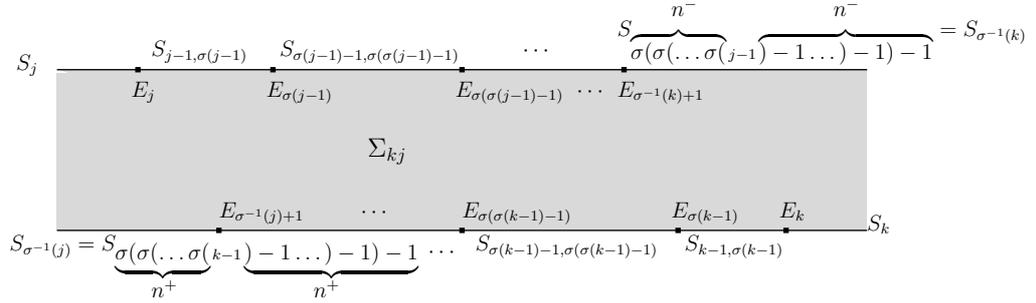}}
    \caption{Depiction of how to build a strip with one clockwise and one counter-clockwise $H$-chain on the boundary. The construction is the same for  open $H$-chains on the boundaries of  half-planes. }
   \label{buildstrip}
 \end{figure}

In the case of a strip $\Sigma_{k,j}$, we label $\R_-$ by $S_{\sigma^{-1}(j)}$, and we label $\R_-+\alpha^-_i(\sim,H)$ by $S_j$ ($\tau:=0$ when the $H$-chain has length 0). Now $S_j$ and $S_{\sigma^{-1}(j)}$ are attached to the clockwise, respectively counter-clockwise, open $H$-chains $H_{\{[\sigma(j-1)]\}}$ and $H_{\{[\sigma^{-1}(j)+1]\}}$. The interval $[0,\tau_{\{[\sigma^{-1}(j)+1]\}} ]$ is labelled by the homoclinic separatrices in $H_{\{[\sigma^{-1}(j)+1]\}}$ and the ends they connect, and $[0,\tau_{\{[\sigma(j-1)]\}}]+\alpha^-_i(\sim,H)$ is labelled by the homoclinic separatrices in $H_{\{[\sigma(j-1)]\}}$ and the ends they connect. Finally, $]\tau_{\{[\sigma^{-1}(j)+1]\}},\infty[$ is labelled by $S_k$, and $]\tau_{\{[\sigma(j-1)]\}},\infty[+\alpha^-_i(\sim,H)$ is labelled by $S_{\sigma^{-1}(k)}$ (see Figure \ref{buildstrip}).\par
\begin{figure}%
    \centering
    \input{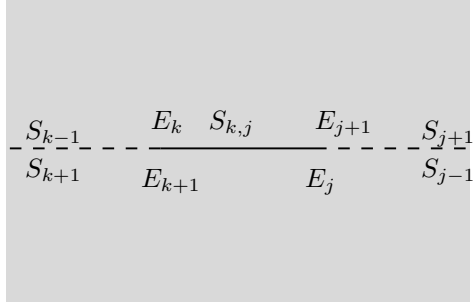}%
    \caption{The neighborhood of an \eqcl\ in  $H$ on the rectified surface $\M$ in the case where the neighborhood consists of two half-planes, along with the
 line segment corresponding to the \eqcl.}
   \label{Uhom1}
 \end{figure}

\begin{figure}%
    \centering
    \input{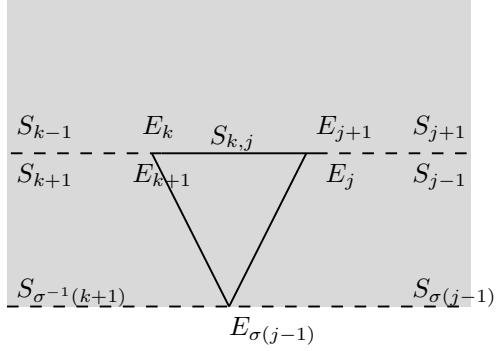}%
    \caption{The neighborhood of an \eqcl\ in $H$ on the rectified surface $\M$ in the case where the
 neighborhood of the \eqcl\ 
 consists of one half-plane, a strip, and the
 line segment corresponding to the \eqcl. }
   \label{Uhom2}
 \end{figure}

\begin{figure}%
    \centering
    \input{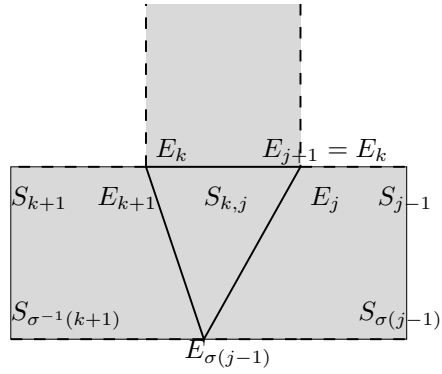}%
    \caption{The neighborhood of an \eqcl\ in $H$ on the rectified surface $\M$ in the case where the neighborhood of the \eqcl\  
 consists of one vertical half-strip sewn to the strip along the
 line segment corresponding to the \eqcl. }
   \label{Uhom3}
 \end{figure}

\begin{figure}%
    \centering
    \input{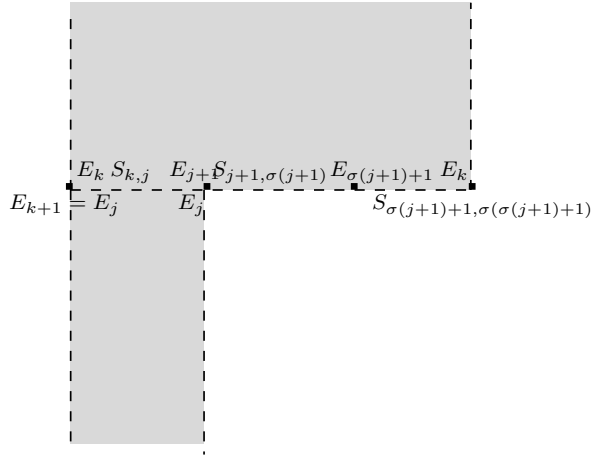}%
   \caption{The neighborhood of an \eqcl\ in $H$ on the rectified surface $\M$ in the case where the \eqcl\ belongs to two different closed $H$-chains. The neighborhood of 
 the \eqcl\
 consists of an upper and lower vertical half-strip, possibly having different widths, along with the
 line segment corresponding to the \eqcl.}
   \label{Uhom4}
 \end{figure}

\begin{figure}%
    \centering
    \input{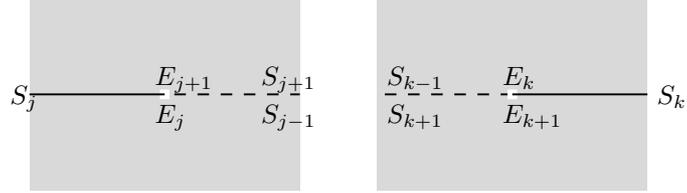}%
    \caption{The figure on the left is a neighborhood of a separatrix $S_j$, $j \in [m]$. It consists of an upper and lower half-plane sewn along the separatrix $S_j$, so it can be represented by $\C \setminus \left(\R_+
    \cup \{ 0 \}\right)$. The figure on the right is a neighborhood of a separatrix $S_k$, $k \in [m]$. It consists of an upper and lower half-plane sewn along the separatrix $S_k$, so it can be represented by $\C \setminus \left(\R_-
    \cup \{ 0 \}\right)$. }
   \label{Upetals1}
 \end{figure}

Every separatrix is on the upper boundary of  exactly one rectified zone and on the lower boundary of exactly one rectified zone. 
In the non-homoclinic case, it is  on the boundary of:
\begin{itemize}
\item
two strips,
\item
one strip and one half-plane, or
\item
two half-planes (see Figure \ref{Upetals1}).
\end{itemize}
In the homoclinic separatrix case, the separatrix is on the boundary of: 
\begin{itemize}
\item
two strips,
\item
one strip and one half-plane  (see Figure \ref{Uhom2}), 
\item
one strip and one vertical half-strip  (see Figure \ref{Uhom3}),
\item
two half-planes (see Figure \ref{Uhom1}),
\item
one half-plane and one vertical half-strip, or
\item
two vertical half-strips (see Figure \ref{Uhom4}).
\end{itemize}
Every separatrix admits a unit-speed parameterization, one starting from the marked point labeled by $E_{\ell}$ and the other from $E_{\ell + 1}$ always such that if $E_j$ is seen from below, then $E_{j+1}$ is seen from above corresponding to the ends on either side of $S_j$, and if $E_{k+1}$ is seen from below, then $E_k$ is seen from above corresponding to the ends on either side of $S_k$ (compare Figures \ref{buildccwcyl} and  \ref{buildcwcyl}).  \par
Now that we have constructed the vertical half strips $C_{\{[\ell]\}}$, half-planes $\pm \HH_{\ell}$, and strips $\Sigma_{k,j}$ with proper labeling, 
\begin{equation}
\M^{\ast}=\left(\coprod \overline{C}_{\{[\ell]\}}\coprod \pm\overline{ \HH}_{\ell}\coprod \overline{\Sigma}_{k,j}\right)/\sim,
\end{equation}
where $\sim$ is the \eqre \ such that all ends are identified, $E=\{E_{\ell} \}/\sim$, each pair of points $z\in S_{\ell}$ or $z \in S_{k,j}$ on the two occurences of any separatrix respecting the unit-speed parameterization are identified, and $z\sim z+\tau_{\{[\ell]\}}$  for $z\in i\R_{\geq 0}\subset \partial C_{\{[\ell]\}}$.\par
\begin{proposition}
The set $\M^{\ast} $ is a Hausdorff topological space.
\end{proposition}
\begin{proof}
Each $\overline{C}_{\{[\ell]\}}$, $\pm\overline{ \HH}_{\ell}$, and $\overline{\Sigma}_{k,j}$ is a Hausdorff topological space with subspace topology from $\C$.  Every disjoint union of Hausdorff spaces is again Hausdorff, so the disjoint union $\coprod \overline{C}_{\{[\ell]\}}\coprod \pm\overline{ \HH}_{\ell}\coprod \overline{\Sigma}_{k,j}$ is a Hausdorff topological space. Now $\M^{\ast}$ is a topological space with the quotient topology, and one can convince oneself that the identified points under $\sim$ do not cause any separability problems.
\end{proof}
As a topological space, $\M^{\ast}$ has $q$ ends we call $\infty_{[\ell]}$ corresponding to the  \eqcl es  $[\ell]\in L$ and $c=c(\sim,H)$ ends $\infty_{\{[\ell]\}}$ corresponding to the closed $H$-chains $H_{\{ [\ell] \}}$.
We will make charts of the neighborhoods of the ends to show that each end corresponds to one point.  The natural $(q+c)$-point closure of $\M^{\ast}$ is denoted $\M$ and is called the \emph{rectified surface}. The construction of the Riemann surface structure on $\M$ is contained in the construction of the atlas that follows.\par
\subsection{An atlas for $\M$}
We  show that $\M$ is a Riemann surface by constructing an atlas for $\M$.  
We construct overlapping open sets $U \subset \M$ covering 
$\M$, and define charts $\eta_U:U \rightarrow \C$.
There are four types of points on $\M$ we need to consider:
\begin{itemize}
\item
points in rectified zones,
\item
points on separatrices, 
\item
the point $E$, and
\item
the ends of $\M^{\ast}$.
\end{itemize}
For the first type of point, we only need to note that each $C_{\{[\ell]\}}$,  $\pm \HH_{\ell}$, and $\Sigma_{k,j}$  has as a natural chart in $\C$ by projection. \par
For points on separatrices, we note the following. Since each separatrix is on the upper boundary of exactly one rectified zone and on the lower  boundary of exactly one rectified zone, we let  $U$ in the non-homoclinic case be the open set in $\M^{\ast}$ which corresponds to the two rectified zones in $\C$ identified along $S_{\ell}$ and let $\eta_U:U\rightarrow \C$ be the chart for which $\eta_U(S_{\ell})=\R_+$ (resp. $\R_-$) for $\ell$ odd (resp. even). In the homoclinic case, we let $U$ be the open set in $\M^{\ast}$ which corresponds to the two rectified zones in $\C$ identified along  $S_{k,j}\subseteq ]0, \tau[$, where $\tau$ is the length of $S_{k,j}$.\par
We define a chart $\eta_E:U_E\rightarrow \D_{\frac{\sqrt[d-1]{r}}{d-1}}$ around $E$ as follows.\par
\begin{figure}%
   \centering
   \resizebox{!}{3cm}{\input{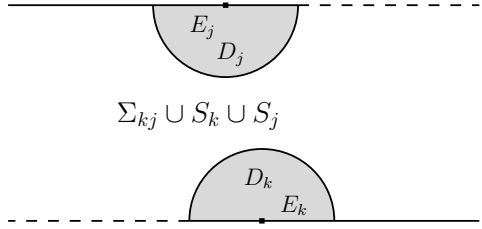}}%
   \caption{Half disks $D^+_k$ and $D^-_j$.}
  \label{DpmB}
\end{figure}
\begin{figure}%
   \centering
   \resizebox{!}{6cm}{\input{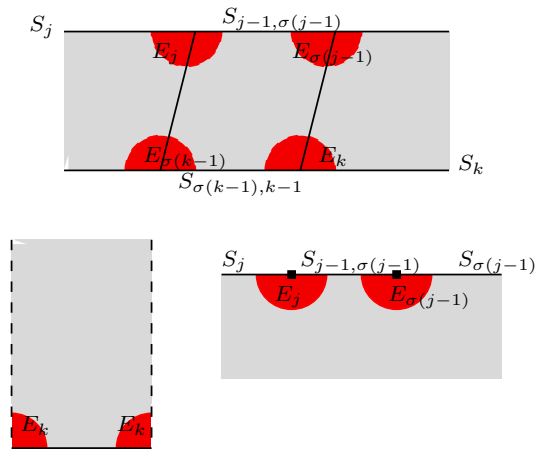}}%
   \caption{Examples of how one can construct half-disks for the ends in each
   rectified zone. }
  \label{nongenend}
\end{figure}

For each $\ell$, 
let
$D_{\ell}$ be the semi-disk in either a strip, half-plane, or 
vertical half-strip with center $E_{\ell}$, ``positive'' (see Figure \ref{DpmB}) if $\ell$ is odd, ``negative'' if
$\ell$ is even, and with radius $r$ sufficiently small, i.e. 
\begin{equation}
r<\frac{1}{2}\left(\min \limits_{i}\{\Im(\alpha_i(\sim,H)) \},\{\tau_i(\sim,H) \}\right) 
\end{equation}
for $i=1,\dots,s$, and $i=1,\dots,h$, respectively
(see Figures~\ref{DpmB} and \ref{nongenend}).  
The neighborhood $U_E \subset \M$ of $E$ is then the half-disks taken with  proper identification, $\{ D_{\ell} \}/\sim$. \par
First consider the map $E \mapsto 0$, giving a $d-1$ covering of $\D_r$ corresponding to  $\arg(Z-E_{\ell}) \in [ -(\ell-1)\pi, -(\ell-2)\pi]$ in each $D_{\ell}$.
To get a univalent covering, we apply the appropriate $(d-1)^{st}$ root by setting $\phi_E:\D_r\rightarrow \D_{\frac{\sqrt[d-1]{r}}{d-1}}$, where
\begin{equation}
\label{phiE}
\phi_E(z)=\exp \left(\frac{1}{d-1}\mathcal{L}_{(-\ell+\frac{1}{2})\pi}  \right) 
\end{equation}
and $\mathcal{L}_{\left(-\ell+\frac{1}{2}\right)\pi}$ is an appropriate branch of the logarithm.
This takes our $d-1$ covering about 0 to a univalent
covering of $\D_{\frac{\sqrt[d-1]{r}}{d-1}}$ where $S_0$ is asypmtotic to $\R_+$ at infinity 
\begin{figure}%
   \centering
   \resizebox{!}{6.5cm}{\input{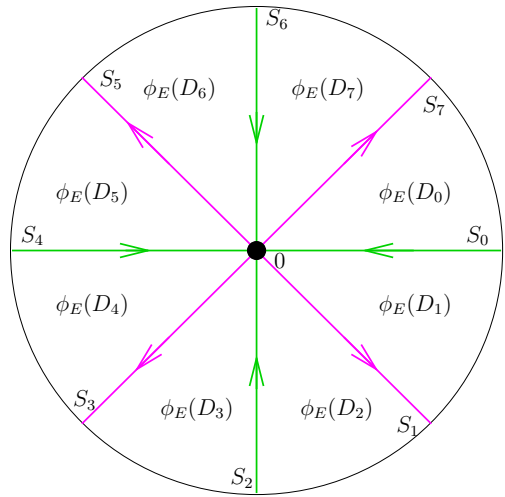}}%
   \caption{A neighborhood of $\eta_E(E)=0$ in the chart $\eta_E$ exemplified for $d=5$.}
  \label{Dr1}
\end{figure} 
(see Figure \ref{Dr1}).  \par
We now construct overlapping open sets $U \subset \M$ covering the ends of 
$\M^{\ast}$, and
define charts $\eta_U:U \rightarrow \C$ determined by $\rho_{[\ell]}$ for each odd or even \eqcl\ and closed $H$-chain.  The chart about mixed \eqcl es is constructed via Fatou coordinates. 

\subsection{Neighborhoods of odd or even \eqcl es and closed $H$-chains}
 Let $\Sigma_{[\ell]}$ be the open set in $\C$ corresponding to an odd or even
equivalence class $[\ell]$ constructed in the following manner. \par
Let $I_{[\ell]}=\{ i \in 1,\dots,s\left(\sim,H\right)\mid j_i \in [\ell] \}$. \par
\begin{figure}%
   \centering
   \input{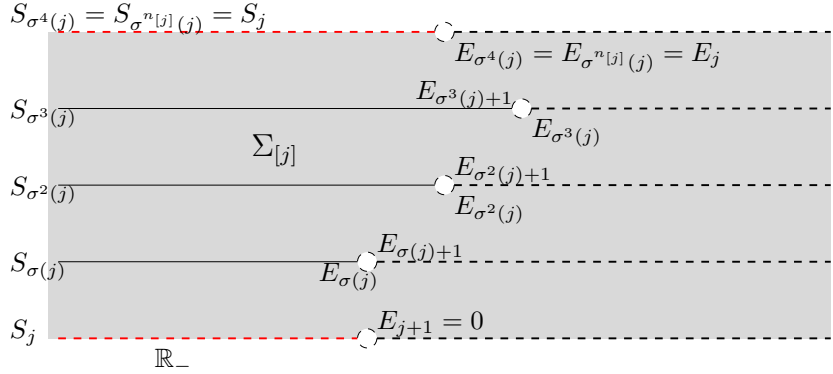}%
   \caption{Example of the set $\Sigma_{[j]}$ for an even \eqcl\ with four elements. }
  \label{fingers1}
\end{figure}

For an even \eqcl\  $[j]$ (see Figure~\ref{fingers1}), we make an open strip in the 
plane $\left\{ z \in \C \mid 0<\Im \left( z \right)<\Im\left( \rho_{[j]}\right)\right\}$.  We take away the half 
line $\{ z \in \C \mid
z=\alpha^-_i\left(\sim,H\right)+t, \quad t \in \mathbb{R}_+ \}$ such that $i \in I_{[j]}$.  We continue this 
process for the 
remainder of
$j\in [j]$.\par
\begin{figure}%
   \centering
   \input{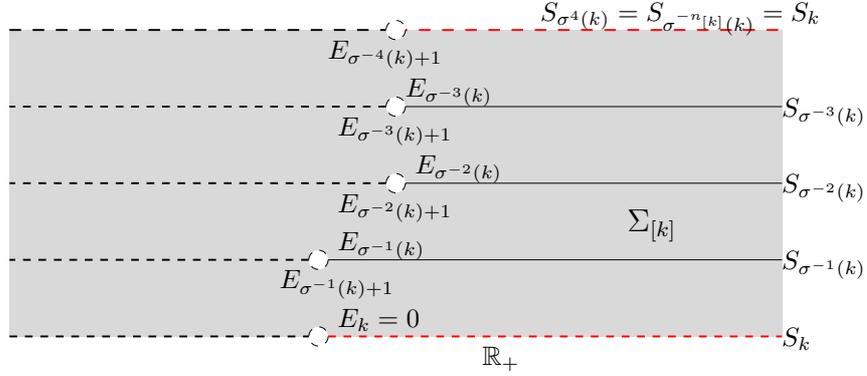}%
   \caption{Example of the set $\Sigma_{[k]}$ for an odd \eqcl\ with four elements.}
  \label{fingers2}
\end{figure}

For an odd \eqcl \ $[k]$ (see Figure~\ref{fingers2}), we make an open strip in the plane $\left\{ z \in \C \mid 0<\Im \left( z 
\right)<\Im \left(\rho_{[k]}\right)\right\}$.  We take away the half line $\{ z \in \C \mid
z=\alpha^+_i\left(\sim,H\right)+t, \quad t \in \mathbb{R}_- \}$.
 We continue this 
process for the 
remainder of
$k \in [k]$ (see Figure~\ref{fingers2}). 
\par
\begin{figure}%
   \centering
  \input{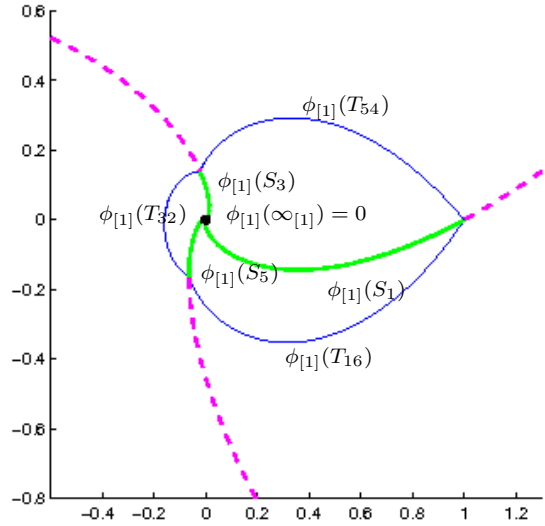}%
   \caption{$V_{[k]}$ for an odd \eqcl\ $[k]$.  Specifically, $[1]=\{1,3,5\}$. The thin dark curves are the images of the rectified
   transversals $T_{16}$, $T_{32}$, and $T_{54}$ under $\phi_{[1]}$. The solid light curves are the images of the
   rectified separatrices $S_1$, $S_3$, and $S_5$. The dashed curves are tails of
   logarithmic spirals which are the images of $\R_-$, $\R_- + \alpha^+_{i_3}$,  $\R_- + \alpha^+_{i_3}+\alpha^+_{i_2}$, and $\R_- + \alpha^+_{i_3}+ \alpha^+_{i_2}+ \alpha^+_{i_1}$, where $i_2=i_1+1$ and $i_3=i_2+1$.}
  \label{Omega_l}
\end{figure}

We first define $U_{[\ell]}^{\ast}=\left(\Sigma_{[\ell]}/ \sim\right)$, where $\sim$ means identifying the appropriate points of $S_{\ell}$. 
We define for $[\ell]$ odd (respectively, $[\ell]$ even) $\tilde{U}_{[\ell]}=U^{\ast}_{[\ell]} \setminus S_{\ell}$. Now $\tilde{U}_{[\ell]}$ is an open subset in $\M^{\ast}$, 
which is 
identified by
$q_{[\ell]}:\tilde{U}_{[\ell]} \rightarrow \Sigma_{[\ell]}$ with 
the open set $\Sigma_{[\ell]} \subset \mathbb{C}$. 
This domain in $\C$ is then mapped by $\phi_{[\ell]}: \Sigma_{[\ell]} \rightarrow \C$ 
defined by
$\phi_{[\ell]}\left(Z\right)=\exp \left(\frac{\mp 2 \pi i}{\rho_{[\ell]}}Z \right)$.  The image $\phi_{[\ell]}\left(\Sigma_{[\ell]}\right)$ is $\C$ 
except for
one complete logarithmic spiral going through the point $z=0$, the origin, and
tails of logarithmic spirals corresponding to $\phi_{[\ell]}\left(\partial 
\Sigma_{[\ell]}\right)$.  Note 
that
$\partial \Sigma_{[\ell]}$ consists of all the dotted lines in Figures \ref{fingers1}
and \ref{fingers2} along with the $E_{\ell_i}$ and  $E_{\ell_i+1}$ where $\ell_i \in [\ell]$. 
The map $\phi_{[\ell]}$ 
extends continuously to
$\varphi_{[\ell]}:\Sigma_{[\ell]}
\cup \mathbb{R}_{\pm} \rightarrow V_{[\ell]}^{\ast}$ where 
$V_{[\ell]}^{\ast}$ 
is $\C$ minus the origin and tails of
logarithmic spirals corresponding to $\partial \Sigma_{[\ell]}\setminus \{ 
\mathbb{R}_{\pm} 
\cup
 \left( \rho_{[\ell]}+\mathbb{R}_{\pm}\right) \}$ (see Figure \ref{Omega_l}). 
 It is defined as follows:
\begin{equation}
\varphi_{[\ell]}\left(Z\right)= \begin{cases} \phi_{[\ell]} & \text{for } Z \in \Sigma_{[\ell]} \\
\exp \left(\frac{\mp 2 \pi
i}{\rho_{[\ell]}} Z \right) & \text{for } Z\in
\mathbb{R}_{\pm}  \end{cases}.
\end{equation} 
For an even or odd \eqcl\ $[\ell]$ , we define $\eta_{[\ell]}:U^{\ast}_{[\ell]} \rightarrow V_{[\ell]}^{\ast}$ by
$\eta_{[\ell]}=\varphi_{[\ell]} \circ q_{[\ell]}$ ($Z \neq 0$). This is a chart on $U^{\ast}_{[\ell]}$ which can be
extended
to a chart
$\eta_{[\ell]}:U_{[\ell]}\rightarrow V_{[\ell]}$ on $U_{[\ell]}$ by mapping the puncture $\infty_{[\ell]}\mapsto 0$.\par
The next case corresponds to a closed $H$-chain $H_{\{[\ell]\}}$, where we remind that clockwise or counter-clockwise orientation are represented by $\ell$ even or odd respectively. We define $U^{\ast}_{\{[\ell]\}}$ to be the open set on $\M$, constructed by
$C_{\{[\ell]\}}/\sim$. 
We define for $[\ell]$ odd (respectively, $[\ell]$ even) 
$\tilde{U}_{\{[\ell]\}}=U^{\ast}_{\{[\ell]\}} \setminus 
 S$ where $S$ is the curve on $\M$ defined by the equivalence 
relation $z \sim z+\tau_{\{[\ell]\}}$, $z \in i \R_{\pm}\subset \partial C_{\{[\ell]\}}$. Now $\tilde{U}_{\{[\ell]\}}$ is an open subset in $\M$, which is 
identified by
$q_{\{[\ell]\}}:\tilde{U}_{\{[\ell]\}} \rightarrow C_{\{[\ell]\}}$ with 
the
open set $C_{\{[\ell]\}}=\{ z \in \pm \HH \mid 0<\Re\left(z\right)<\tau_{\{[\ell]\}} \}$. 
This domain $C_{\{[\ell]\}}\subset \C$ is then mapped by $\phi_{\{[\ell]\}}: C_{\{[\ell]\}} \rightarrow \D$ defined by
$\phi_{\{[\ell]\}}\left(Z\right)=\exp \left( \frac{\mp 2 \pi i}{\tau_{\{[\ell]\}}}Z \right)$.  The image is  $\phi_{\{[\ell]\}}\left(C_{\{[\ell]\}}\right)=\D \setminus [0,1[$. The map $\phi_{\{[\ell]\}}$ 
extends continuously to
$\phi_{\{[\ell]\}}:C_{\{[\ell]\}} \cup  i\R_{\pm} \rightarrow \D^{\ast}$.
It is defined as follows:
\begin{equation}
\varphi_{\{[\ell]\}}\left(Z\right)= \begin{cases} \phi_{\{[\ell]\}} & \text{for } Z \in 
C_{\{[\ell]\}} \\
\exp \left(\frac{\mp 2 \pi
i}{\tau_{\{[\ell]\}}} Z \right) & \text{for } Z\in
 i \R_{\pm}  \end{cases}.
\end{equation} 
For $[\ell]$ odd (resp. even), we define $\eta_{\{[\ell]\}}:U^{\ast}_{\{[\ell]\}} \rightarrow \D^{\ast}$ by
$\eta_{\{[\ell]\}}=\varphi_{\{[\ell]\}} \circ q_{\{[\ell]\}}$ ($Z \neq 0$). 
This is a chart which extends 
to
$\eta_{\{[\ell]\}}:U_{\{[\ell]\}}\rightarrow \D$ by $\infty_{\{[\ell]\}}\mapsto 0$.\par

\subsection{Neighborhoods around Mixed Equivalence Classes $[m]$}
\label{vfmult}
For a mixed \eqcl\ $[m]$ with $p_{[m]}$ parity changes, we define $U^{\ast}_{[m]}$ to be a suitably restricted  domain  on $\M^{\ast}$, identified with appropriate restrictions (to be determined below) of  upper and lower  half-planes $ \HH_k$ and $-\HH_j$ with $k \sim j \in [m]$ and half-strips of the $\Sigma_{k,j}$ such that either $j \sim \sigma^{-1}\left(j\right) \in [m]$ or $k \sim \sigma^{-1}\left(k\right) \in [m]$, with proper identification along the separatrices.  The goal is to show the existence of a chart $\eta_{[m]}:U^{\ast}_{[m]}\rightarrow V_{[m]}^{\ast}$, where $V_{[m]}$ is some neighborhood of 0, and the puncture $\infty_{[m]}$ in $U^{\ast}_{[m]}$ corresponds to 0.\par 
 \begin{figure}%
   \centering
   \resizebox{!}{5cm}{\input{gap.pstex_t}}%
   \caption{Possible depiction of an almost $p_{[m]}/2$  times ``covering'' of $\infty$ in rectifying coordinates with $\Im \left(\rho_{[m]}\right) >0$ and $p_{[m]}=2$.}
  \label{gap}
\end{figure} 
 \begin{figure}%
   \centering
   \resizebox{!}{5cm}{\input{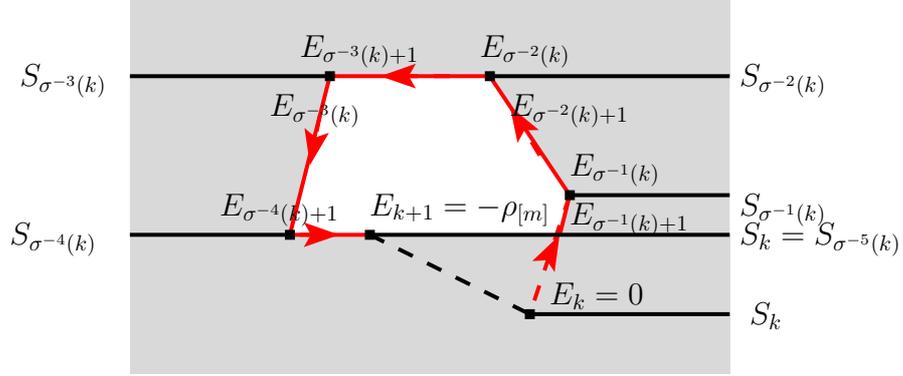}}%
   \caption{Possible depiction of a more than $p_{[m]}/2$  times ``covering'' of $\infty$ in rectifying coordinates with $\Im \left(\rho_{[m]}\right) <0$ and $p_{[m]}=2$. }
  \label{overlap}
\end{figure} 
Pick some $k \in [m]$ and fix in $\C$ the half-plane or half-strip having $S_k$ on its lower boundary such that $S_k$ is equal  to $\R_+$. Continue building counter-clockwise in this way with identification corresponding to $[m]$ until one comes back to the other represenation of the separatrix $S_k$ on the upper boundary of some half-plane or half-strip, which will be identified with the half-line $\R_+- \rho_{[m]}$.The half-planes and half-strips will wind around approximately $p_{[m]}/2$ times in the following sense: a neighborhood of $\infty$ is ``covered'' $p_{[m]}/2$ times if $\Im \left(\rho_{[m]}\right) =0$, almost $p_{[m]}/2$ times if $\Im \left(\rho_{[m]}\right) >0$ (see Figure \ref{gap}), and a bit more than $p_{[m]}/2$ times if $\Im \left(\rho_{[m]}\right) <0$ (see Figure \ref{overlap}).
Let $\D_R$ be the disk with radius
\begin{equation}
\label{DR}
R=  \sum\limits_{i=1}^s \left|\alpha_i \right|+ \sum\limits_{i=1}^h \left|\tau_i \right|.
\end{equation}
Then through several steps we can construct a chart $\eta_{[m]}:U^{\ast}_{[m]}\rightarrow V_{[m]}^{\ast}$. The following steps in the construction of the chart $\eta_{[m]}$ are inspired by Shishikura's treatment of Fatou coordinates in \cite{Shi}. \par
Consider the map $\pi:\C \rightarrow \C$ defined by
\begin{equation}
\pi\left(z\right)=\exp\left(\pi i p_{[m]} z\right).
\end{equation}
Note that $\pi \left(z+ \frac{n}{p_{[m]}/2}  \right)=\pi\left(z\right)$ for all $n \in \Z$.\par
\begin{figure}%
   \centering
   \resizebox{!}{6cm}{\input{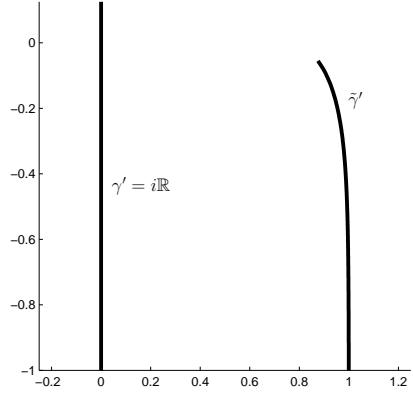}}%
   \caption{Graph of  $\gamma'$, the unique preimage of $\R_+$ under $\pi$ such that $\Re \left(z\right)=0$, and $\tilde{\gamma}'$, the unique preimage of the half-line $\R_+- \rho_{[m]}$  under $\pi$ for which  $\Re \left(z\right) \rightarrow 1$  for  $\Im \left(z\right) \rightarrow - \infty$ on $\tilde{\gamma}'$. In this example, $\rho_{[m]}=-1+2i$ and $p_{[m]}=2$ (compare with Figure \ref{gap}),  so  the preimage contained in the strip $0<\Re \left(z\right)<1$.  }
  \label{l_lprime}
\end{figure}

  Set $\gamma '=i \R$, i.e. the unique preimage of $\R_+$ under $\pi$ such that $\Re \left(z\right)=0$ for $z \in \gamma'$, and let $\tilde{\gamma}'$
 be the unique preimage of the half-line $\R_+- \rho_{[m]}$  under $\pi$ contained in the strip $\frac{p_{[m]}/2-1}{p_{[m]}/2}<\Re \left(z\right)<1$ when $\Im \left(\rho_{[m]}\right) > 0$ ($1<\Re \left(z\right)<\frac{p_{[m]}/2+1}{p_{[m]}/2}$ when $\Im \left(\rho_{[m]}\right)< 0$, and $\Re \left(z\right)=1$ when $\Im\left(\rho_{[m]}\right) = 0$).  Note that in each case $\Re \left(z\right) \rightarrow 1$ for  $\Im \left(z\right) \rightarrow - \infty$ on $\tilde{\gamma}'$ (see Figure \ref{l_lprime}).  \par
We define the map $\tilde{f}:\gamma' \rightarrow \tilde{\gamma}'$ such that it satisfies
\begin{equation}
\left(\pi \circ \tilde{f} \right)\left(z\right)= \pi\left(z\right) +\tau, \quad z \in \gamma'
\end{equation}
with $\tau=-\rho_{[m]}$, i.e. so that the proper identification is preserved.
From this relation, one can deduce the expression 
\begin{equation}
\label{ftilde}
\tilde{f}\left(z\right)=z + 1+ \frac{\Log \left( 1+ \tau \exp\left(- \pi i p_{[m]}z\right) \right)}{ \pi i p_{[m]}}, \quad z \in \gamma'.
\end{equation}
The function $\tilde{f}$ can be extended holomorphically in some lower half plane $H_1=\{ z\in \C \mid \Im\left(z\right)<c_1<0 \}$.  The constant $c_1$ can be determined by forcing the condition $|\tau \exp\left(- \pi i p_{[m]}z\right)|<1$ so that $\Log$ is analytic in $H_1$. Then $c_1$ becomes
\begin{equation}
c_1=\frac{-\ln|\tau|}{ \pi p_{[m]}}.
\end{equation}
\begin{proposition}
There exists a half-plane $H_2\subset H_1$ such that  $\tilde{f}$ satisfies the properties
\begin{align}
\label{properties1}
|\tilde{f}\left(z\right)-\left(z+1\right)|&<1/4 \\
\label{properties2}
|\tilde{f}'\left(z\right)-1|&<1/4. 
\end{align}
\end{proposition}
\begin{proof}
Suppose $|\tau| \exp\left( \pi p_{[m]} \Im\left(z\right)\right)<R<1$ for some $R$.  
This can be rewritten as
\begin{equation}
\Im \left(z\right)< c_1+\frac{\ln \left(R\right)}{\pi p_{[m]}},
\end{equation}
for which we have
 \begin{equation}
 |\tilde{f}'\left(z\right)-1|=\frac{|-\tau \exp\left(- \pi i p_{[m]} z\right)|}{|1+\tau \exp\left(- \pi i p_{[m]}z\right)|}\leq \frac{R}{1-R}. 
 \end{equation}
For the other condition, we can calculate
\begin{equation}
|\tilde{f}\left(z\right)-\left(z+1\right)|=\frac{|\Log\left(1+\tau \exp\left(- \pi i p_{[m]}z\right)\right)|}{| \pi i p_{[m]}|}
\end{equation}
and  $|\Arg\left(1+\tau \exp\left(- \pi i p_{[m]}z\right)\right)|<\Arcsin\left(R\right)$, giving
\begin{equation}
|\tilde{f}\left(z\right)-\left(z+1\right)|^2<\frac{1}{\left( \pi p_{[m]} \right)^2}\left[ \ln^2\left(1+R\right)+\Arcsin^2\left(R\right) \right].
\end{equation}
Since $\frac{R}{1-R}$, $\ln\left(1+R\right)$, and $\Arcsin\left(R\right)$ tends to 0 as $R \rightarrow 0$, we can conclude that there exists $R<1$ such that $|\tilde{f}\left(z\right)-\left(z+1\right)|<\frac{1}{4}$ and $|\tilde{f}'\left(z\right)-1|<1/4$. 
  Choose  $c_2$ such that  the above conditions are satisfied in the half-plane $H_2=\{z \mid \Im \left(z\right)<c_2 \}$. 
\end{proof}
\begin{lemma}
\label{funi}
$\tilde{f}$ is univalent in $H_2$.
\end{lemma}
\begin{proof}
Take two points $z_1 \neq z_2$ in $H_2$.  We show that $\tilde{f}\left(z_1\right)\neq \tilde{f}\left(z_2\right)$.  One can always connect two points in the half-plane $H_2$ with a straight line segment.  We parameterize this line segment by
\begin{equation}
l\left(t\right)=\left(1-t\right)z_1+t z_2, \quad t \in [0,1].
\end{equation}
Let $L\left(t\right)=\tilde{f}\left(l\left(t\right)\right)$. Then
\begin{equation}
\tilde{f}\left(z_2\right)-\tilde{f}\left(z_1\right)=\int_0^1 L'\left(t\right)dt,
\end{equation}
and since $l' \left(t\right)=z_2- z_1\neq 0$ is a constant,
\begin{equation}
  \tilde{f}\left(z_2\right)-\tilde{f}\left(z_1\right)=0 \Leftrightarrow \int_0^1 \tilde{f}'\left(l\left(t\right)\right)dt=0.
\end{equation}
But 
\begin{equation}
\int_0^1 \Re \left( \tilde{f}'\left(l\left(t\right)\right)\right) dt >0
\end{equation}
by Equation \eqref{properties2},
so $\tilde{f}\left(z_2\right) \neq \tilde{f}\left(z_1\right)$. 
\end{proof}
We now construct a univalent function $\Phi$ in a closed region $\Omega$ (defined below) satisfying $\Phi\left(\tilde{f}\left(z\right)\right)=\Phi\left(z\right)+1$ on $\gamma=\Omega\cap \gamma'$.
Using the same method as in \cite{Shi}, let $\Sigma=\{ z \mid 0 \leq \Re \left(z\right) \leq 1 , \Im \left(z\right) \leq c_3< c_2\}$ and $\Omega=h_1\left(\Sigma\right)$ where $h_1:\Sigma \rightarrow \Omega$ is defined by
\begin{equation}
h_1\left(x,y\right)=\left(1-x\right)\left(iy\right)+x \tilde{f}\left(iy\right), \quad 0 \leq x \leq 1, y\leq c_3.
\end{equation}
We restrict to $y\leq c_3<c_2$ since we want $\tilde{f}$ to be univalent and satisfy properties \eqref{properties1} and \eqref{properties2}  for $y=c_3$. \par
For $x=0$, we are on $\gamma$, and for any $y \leq c_3$ and $x=1$, we are on $\tilde{f}\left(\gamma\right)=\tilde{\gamma}=\Omega \cap \tilde{\gamma}'$ and the image of the horizontal line segments for $x \in [0,1]$ in $\Sigma$ are the line segments joining $iy$ and $\tilde{f}\left(iy\right)$ in $\Omega$. \par
Now $h_1$ is $\mathcal{C}^1$ in $\Sigma^{\circ}$ as a function of the two real variables $\left(x,y\right)$ into $\R^2$ since $\tilde{f}$ is holomorphic.  We show that $h_1$ is a homeomorphism.
\begin{lemma}
 $h_1:\Sigma^{\circ} \rightarrow \Omega^{\circ}$ is a homeomorphism. 
\end{lemma}
\begin{proof}
It is enough to prove that $H_1$ is bijective since continuity follows from $h_1\left(x,y\right)$ being $\mathcal{C}^1$.\par
The function $h_1:\Sigma^{\circ}\rightarrow \Omega^{\circ}=h_1\left(\Sigma^{\circ}\right)$ is trivially surjective, so injectivity is all that remains to be proven.\par
There is a natural parameterization of $\gamma$  by  $y$, and  $\Im \left(\tilde{f}\left(iy\right)\right)$ must decrease when $y$ decreases by \eqref{properties2} for all $iy\in \gamma\subset H_2$.
Choose any two points $z_1$ and $z_2$ in $\partial \Sigma$ such that $\Re \left(z_1\right)=0$ and $\Re \left(z_2\right)=1$ and $\Im \left(z_1\right) =\Im \left(z_2\right)$.  These points are  mapped by $h_1$ respectively to $z \in \gamma$ and $\tilde{f}\left(z\right) \in \tilde{f}\left(\gamma\right)$, and the horizontal line segment joining $z_1$ and $z_2$ is mapped under $h_1$ onto the line segment which joins $z$ and $\tilde{f}\left(z\right)$. This line segment is  completely contained in $\Omega$ by the above parameterization property.  
Now the only way $h_1$ could fail to be injective is if there were two distinct line segments in $\Sigma$ mapped to two distinct line segments in $\Omega$ that intersect at some point in the interior of $\Omega$.  This cannot happen due to the above argument of the respective orientations of $\gamma$ and $\tilde{f}\left(\gamma\right)$. 
\end{proof}
By calculations identical to those in \cite{Shi}, it can be shown that $h_1$ has dilatation
\begin{equation}
\left| \frac{\frac{\partial h_1}{\partial \overline{z}}}{\frac{\partial h_1}{\partial z}} \right| <1/3.
\end{equation}
Hence, $h_1$ is a quasi-conformal mapping of $\Sigma^{\circ}$ onto the half-strip $\Omega^{\circ}$, where $h_1^{-1}\left(\tilde{f}\left(z\right)\right)=h_1^{-1}\left(z\right)+1$ for $z \in \gamma$.\par
Set $T\left(z\right)=z+1$. Note that for $iy \in \gamma$ we have $h_1^{-1}\left(\tilde{f}\left(iy\right)\right)=h_1^{-1}\left(iy\right)+1$, i.e. the diagram
\begin{equation}
\begin{CD}
\gamma @>h_1^{-1}=id >>\gamma\\
@VV\tilde{f}V @VVTV\\
 \tilde{f}\left(\gamma\right) @>h_1^{-1}>> T\left(\gamma\right)
\end{CD}
\end{equation}
commutes.\par
Use $h_1: \Sigma \rightarrow \Omega$ to pullback the standard complex structure $\sigma_0$ on $\Omega^{\circ}$ to $\sigma$ on $\Sigma^{\circ}$ by $\sigma=h_1^{\ast}\sigma_0$. Extend $\sigma$ to $H_3 \setminus \{ z \in \C \mid \Re\left(z\right)\in \Z \}$, where $H_3=\{z \mid \Im \left(z\right)<c_3  \}$, through further pullback of $\sigma$ using $T^n\left(z\right)$. For any $n \in \Z$,
\begin{equation}
\left(T^n\right)^{-1}\left(\Sigma^{\circ}\right)=\{ z \in \C \mid \Re \left(z\right)\in \left(]-n,-\left(n-1\right)[\right) \wedge \left(y<0\right) \}.
\end{equation}
Define $\sigma$ on $\left(T^n\right)^{-1}\left(\Sigma^{\circ}\right)$ by $\left(T^n\right)^{\ast}\sigma$, i.e. 
\begin{equation}
 \sigma \restrictto{\left(T^n\right)^{-1}\left(\Sigma^{\circ} \right)} =\left(T^n\right)^{\ast}\left(\sigma \restrictto{\Sigma^{\circ}}\right).
\end{equation}
\par 
Then extend $\sigma$ to $\tilde{\sigma}$ on $\C$ by
\begin{equation}
\tilde{\sigma}=
\begin{cases} \sigma, & z \in H_3 \setminus \{ z \in \C \mid \Re\left(z\right)\in \Z \}.
\\
\sigma_0, & \text{elsewhere in $\C$}
\end{cases}
\end{equation}
By construction, $\tilde{\sigma}$ is a $T$-invariant almost complex structure on $\C$, i.e. $T^{\ast}\tilde{\sigma}=\tilde{\sigma}$.\par
By the  measurable Riemann Mapping Theorem, there exists a unique quasi-conformal mapping $h_2:\C \rightarrow \C $ such that $h_2^{\ast}\sigma_0=\tilde{\sigma}$ and normalized by fixing $i\tilde{c}$ and $1+i\tilde{c}$ for some $\tilde{c}<c_3$. \par
Note that $h_2 \circ T \circ h_2^{-1}$ is a holomorphic bijection 
\begin{equation}
\left(\C,\sigma_0\right) \xrightarrow{h_2^{-1}} \left(\C,\tilde{\sigma}\right) \xrightarrow{T}\left(\C,\tilde{\sigma}\right)\xrightarrow{h_2}\left(\C,\sigma_0\right)
\end{equation} 
without fixed points, since $h_2$ is a homeomorphism and $T$ has no fixed points in $\C$.  Hence, $h_2 \circ T \circ h_2^{-1}$ is an affine map of the form $z \mapsto z+b$. By the chosen normalization of $h_2$, we have $h_2 \circ T \circ h_2^{-1}\left(i\tilde{c}\right)=h_2 \circ T \left(i\tilde{c}\right)=h_2\left(1+i\tilde{c}\right)=1+i\tilde{c}$. 
 Therefore, $b=1$ so $h_2 \circ T \circ h_2^{-1}=T$.  \par
Define $\Phi=h_2 \circ h_1^{-1}$ on $\Omega$ satisfying $\Phi\left(\tilde{f}\left(z\right)\right)=\Phi\left(z\right)+1$ on $\gamma=\Omega\cap \gamma'$ respecting the diagram
\begin{equation}
\begin{CD}
\gamma @>h_1^{-1}=id >>\gamma @> h_2 >> h_2\left(\gamma\right)\\
@VV\tilde{f}V @VVTV @VVTV\\
 \tilde{f}\left(\gamma\right) @>h_1^{-1}>> T\left(\gamma\right)@> h_2 >> h_2\left(\gamma+1\right).
\end{CD}
\end{equation}
 Then $\Phi$ is well-defined  and homeomorphic on $\Omega$ by the properties of $h_1$ and $h_2$. Furthermore, $\Phi$ is holomorphic on $\Omega^{\circ}$.\par
\begin{figure}%
   \centering
   \resizebox{!}{8cm}{\input{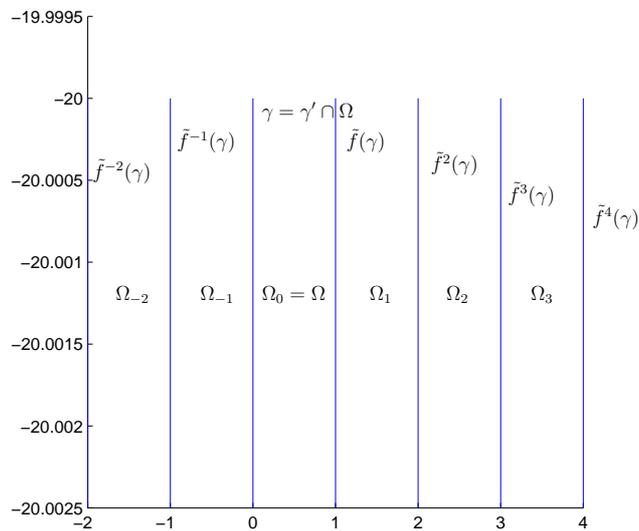}}%
   \caption{Depiction of the domain $U$ with labelling of the $\Omega_j$s and the images of $\gamma$. Note that these curves resemble straight lines rather quickly for $\Re (z) \rightarrow -2,\dots,4$.}
  \label{Sjs}
\end{figure}

We next extend $\Phi$  to $U=interior \left(\bigcup \limits_{n=-2}^3 \Omega_n \right)$ where $\Omega_n=\tilde{f}^n\left(\Omega\right)$ by 
\begin{equation}
\label{phiextend}
\Phi\left(z\right)=\Phi\left(\tilde{f}^{-n}\left(z\right)\right)+n, \quad z \in \Omega_n.
\end{equation}
Even though the domain $U$ is much larger than we need in order to prove that $\Phi$ is holomorphic on $\gamma$ and $\tilde{f}\left(\gamma\right)$, it is chosen in this way so that we have that $\Phi$ is univalent in $\{ \zeta \mid |\zeta-z|<5/4+\epsilon, z \in \Omega \}$ for later use in some deformation estimates (see Section \ref{assvfs}).   We have that $\Phi\left(z\right)$ is  well-defined by Equation \eqref{phiextend} in $U_1=U \cap H_4$ where $H_4=\{ z \mid \Im \left(z\right)<c_4=c_3-3/4 \}$ if there exists a unique integer $-2 \leq n \leq 3 $ such that $\tilde{f}^n\left(z\right) \in \Omega$ for $z \in U$ (see Figure \ref{Sjs}).
\begin{proposition} 
\label{uniqueSn}
There exists a unique integer $-2 \leq n \leq 3 $ such that $\tilde{f}^n\left(z\right) \in \Omega$ for $z \in U$.  
\end{proposition}
\begin{proof}
Consider points $z_1$ and $z_2$ in $\Omega_{n-1}$ and $\Omega_n$ respectively such that $\Im\left(z_1\right)=\Im\left(z_2\right)$ and connect them by a horizontal line segment intersecting the point $z_0$ on $\gamma_{n}=\tilde{f}^n\left(\gamma\right)$. 
The image of this horizontal line segment under $\tilde{f}$ is naturally parameterized by 
\begin{equation}
\gamma_n\left(x\right)=\tilde{f}\left(x+iy\right),
\end{equation}
where $y$ is constant. To see how much $\gamma_n\left(x\right)$ varies from the horizontal, we calculate
\begin{equation}
|\gamma_n'\left(x\right)-1|=\left|\tilde{f}'\left(x+iy\right)-1\right|<1/4,
\end{equation}
so it follows since $\Re \left(\tilde{f}'\right)>0$ and $\left|\Im \left(\tilde{f}'\right)\right|<1/4$ that $\gamma_n\left(x\right)$ is contained in the region bounded by the two lines with slope $\pm 1/4$ that intersect at $\tilde{f}\left(z_0\right)$, and furthermore, that $\gamma_n\left(x\right)|_{\Re \left(z_1\right)<x<\Re \left(z_0\right)}$ is contained in the left sector and  $\gamma_n\left(x\right)|_{\Re \left(z_0\right)<x<\Re \left(z_2\right)}$ is contained in the right sector (see Figure \ref{gammaj}).  
 \begin{figure}%
   \centering
   \resizebox{!}{6cm}{\input{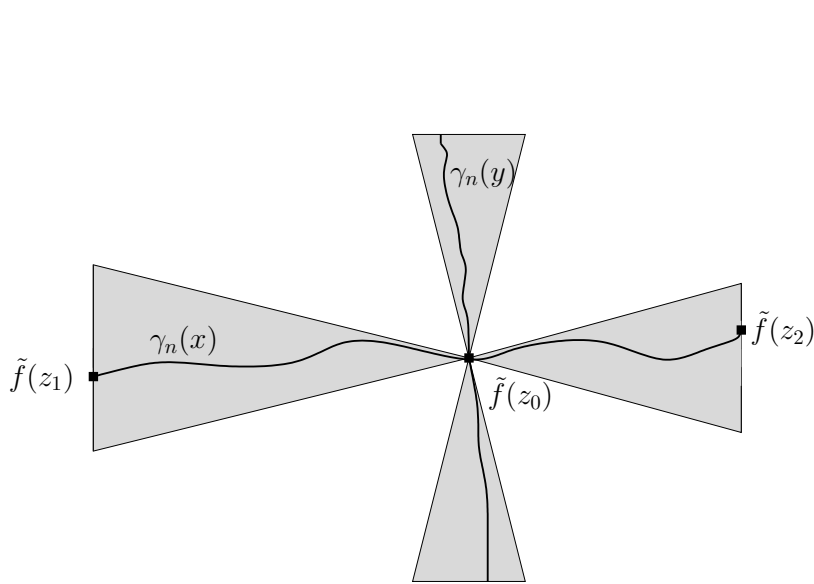}}%
   \caption{Illustration of the proof of Proposition \ref{uniqueSn} }
  \label{gammaj}
\end{figure} 
Now to see that  $\gamma_n\left(x\right)|_{\Re \left(z_1\right)<x<\Re \left(z_0\right)}\subset \Omega_{n}$ and  $\gamma_n\left(x\right)|_{\Re \left(z_0\right)<x<\Re \left(z_2\right)}\subset \Omega_{n+1}$, we need to show that $\gamma_n$, naturally parameterized by $\gamma_n\left(y\right)=\tilde{f}^n\left(iy\right)$, is contained in the region bounded by the two lines with slope $\pm 4$ that intersect at $\tilde{f}\left(z_0\right)$. That is
\begin{equation}
\label{lslope}
|\gamma_n'\left(y\right)-i|<1/4, \quad n=-2,\dots,3.
\end{equation}
For $n=1$, we have 
\begin{equation}
|\gamma_1'\left(y\right)-i|=|\tilde{f}'\left(iy\right)-1|<1/4,
\end{equation}
and noticing that
\begin{equation}
\gamma_n'\left(y\right)=\tilde{f}'\left(\gamma_{n-1}\left(y\right)\right)\cdot\gamma_{n-1}'\left(y\right),
\end{equation}
one can show that inequality \eqref{lslope} holds.  Since the vertical and horizontal sectors are disjoint, we can conclude that each point $z \in \bigcup \limits_{n=-2}^3\Omega_n$ is mapped into $\Omega=\Omega_0$ exactly once by $\tilde{f}^{-n}$, for $-2 \leq n \leq 3$.
\end{proof}
We have that $\Phi$ is holomorphic and univalent in $U \setminus \{ \bigcup \limits_j\gamma_j\}$ and is continuous and injective on the $\gamma_j$. Therefore, $\Phi$ is holomorphic and univalent on all of $U$ as desired by Morera's Theorem.\par
Let $y_0<c_4$ be such that the image of the line segment from $i y_0$ to $\tilde{f}(iy_0)$ under $\pi$ is completely contained in $\C \setminus \overline{\D}_R$ (see \eqref{DR}). Define 
\begin{equation}
U_{[m]}^{\ast}=\pi \left( S \cap \left\{ \Im (z)-y_0<\frac{\Im \left(\tilde{f}(iy_0) \right)-y_0}{\Re \left(\tilde{f}(iy_0) \right)}\Re(z) \right\} \right).
\end{equation}
The complex structure on $\Phi\left(\gamma\right)$ is translation invariant, which allows post-compositon with the map $\exp\left(-2 \pi i (\cdot)\right):\Phi\left(\Omega\right) \rightarrow V_{[m]}^{\ast}$, $\eta_{[m]}=\exp\left(-2 \pi i (\cdot)\right) \circ \Phi \circ \pi^{-1}$,  defining the neighborhood $V_{[m]}$ of zero. This defines a chart $\eta_{[m]}:U_{[m]}^{\ast}\rightarrow V_{[m]}^{\ast}$, which can be continuously extended so that $\infty_{[m]}\mapsto 0$. 

\subsection{Holomorphic Overlap}
The images of each connected component of the non-empty intersections $U_1\cap U_2$ of the open sets in $\M$ defining its charts are subsets of each $V_i=\eta_i\left(U_i\right)$, $i=1,2$ of one of the following types:
\begin{itemize}
\item
a strip, half-plane, or vertical half-strip, correspoinding to regular points on $\M$,
\item
combinations of any two of the above glued along the common bounding separatrix, corresponding to neighborhoods of separatrices on $\M$,
\item
 $\D_{\frac{\sqrt[d-1]{r}}{d-1}}$, corresponding to $U_E$,
\item
$V_{[\ell]}$, corresponding to a source or sink,
\item
$\D$ for a center, or
\item
$V_{[m]}$ for a \multeq .
\end{itemize}
If $V_1$, $V_2$ are one of the first two types, then both $\eta_1\left(U_1 \cap U_2\right)$ and $\eta_2\left(U_1 \cap U_2\right)$ are whole strips, planes, or  vertical half-strips and 
\begin{equation}
\eta_2 \circ \eta_1^{-1}=T, 
\end{equation}
for some translation $T$.\par
If $V_1$ is one of the first two types and $V_2=\eta_E\left(U_2=U_E\right)$, then $\eta_1\left(U_1 \cap U_E\right)$ is a union of at most two half-disks, possibly joined by a common bounding separatrix and
\begin{equation}
\eta_1 \circ \eta_E^{-1}=T \circ z^{d-1}, 
\end{equation}
for some translation $T$. This is clearly biholomorphic in the domain of definition when the appropriate $d-1^{st}$ branch is chosen.\par
If $V_1=V_{[\ell]}$ and $V_2=\eta_E\left(U_2=U_E\right)$, then each connected component of $U_1\cap U_E$ corresponds to exactly two half-disks $D_{\ell}$ and $D_{\ell+1}$ joined by part of the common bounding separatrix $S_{\ell}$ and
\begin{equation}
\eta_1 \circ \eta_E^{-1}=\varphi_{[\ell]}\circ T \circ z^{d-1}, 
\end{equation}
for some translation $T$. This is clearly biholomorphic in each connected component of the domain of definition when the appropriate $d-1^{st}$ branches are chosen. \par
If $V_1=V_{[k]}$ and $V_2=V_{[j]}$, then $\eta_1\left(U_1 \cap U_2\right)$ and $\eta_2\left(U_1 \cap U_2\right)$ are whole strips (maybe 2 whole strips if $d=2$), and
\begin{equation}
\eta_2 \circ \eta_1^{-1}=\varphi_{[j]}\circ T \circ\phi_{[k]}^{-1}, 
\end{equation}
for some translation $T$. The map 
\begin{equation}
\label{phil1}
 \phi_{[k]}^{-1}= \frac{\rho_{[k]}}{-2 \pi i} \circ\mathcal{L}_{[k]}
\end{equation}
is well-defined and biholomorphic in the domain of definition, since the logarithm $\mathcal{L}_{[k]}:\tilde{V}_{[k]} \rightarrow \Sigma_{[k]}$  taking the head 
of
the 
logarithmic spiral passing through $z=1$ to $\frac{-2 \pi i}{\rho_{[k]}} 
\mathbb{R}_{+}$ is well defined and holomorphic. \par
If $V_1=\D$ corresponding to an $H$-chain $H_{\{[\ell]\}}$, then $U_2$ must be of one of the first 2 types where $U_1\cap U_2$ can only correspond to a vertical half-strip and
\begin{equation}
\eta_2 \circ \eta_{\{[\ell]\}}^{-1}=T \circ \frac{\tau_{\{[\ell]\}}}{\mp2 \pi i} \circ \mathcal{L}_0,
\end{equation}
which is biholomorphic on $\D \setminus 
[0,1[$ since there is a well-defined and holomorphic branch of the logarithm $\LL_0:\D 
\setminus [0,1[ \rightarrow ]-\infty,0[ \times ]0, 2 \pi[$.\par
If $V_1=V_{[m]}$ corresponding to a \multeq , and $V_2$ is one of the first 2 types, then at worst $U_1 \cap U_2$ corresponds to part of one of the second types and 
\begin{equation}
\eta_{[m]} \circ \eta_{2}^{-1}=\exp\left(-2\pi i\left(\cdot\right)\right)\circ \Phi \circ \pi^{-1}\circ T, 
\end{equation}
for some translation $T$, and $\Phi$ and $\pi$ as in Section \ref{vfmult}. This map is biholomorphic since one can choose a well-defined branch of $\pi^{-1}$ in the domain of definition, so biholomorphicity of $\exp\left(-2\pi i\left(\cdot\right)\right)\circ \Phi \circ \pi^{-1}$ follows from injectivity and the holomorphic Inverse Function Theorem. \par
If $V_1=V_{[m]}$ corresponding to a \multeq , and $V_2=V_{[\ell]}$, then $U_1 \cap U_2$ corresponds to part of exactly one strip and
\begin{equation}
\eta_{[m]} \circ \eta_{[\ell]}^{-1}=\exp\left(-2\pi i\left(\cdot\right)\right)\circ \Phi \circ \pi^{-1}\circ T \circ\phi_{[\ell]}^{-1}, 
\end{equation}
for some translation $T$. This map is biholomorphic on its domain of definition by the previous arguments for $\exp\left(-2\pi i\left(\cdot\right)\right)\circ \Phi \circ \pi^{-1}$ and $\phi^{-1}_{[\ell]}$.

Therefore, we have made an atlas  with \hol\ overlap covering $\M$, and we can conclude that $\M$ is a Riemann surface. 
\section{$\M$ isomorphic to $\rs$}
\label{MisoC}
\begin{proposition}
$\M$ is homeomorphic to $\rs$.
\end{proposition} 
\begin{proof}
We utilize the Euler characteristic. Note that 
$\M$ has a topology induced by the topology of $\C$.  The end $E$ gives one 
vertex, and 
all other equilibrium points apart from centers contribute a vertex.  That is, 
there are 
a total of $1+d-c-\frac{p_P}{2}$ vertices. The number of edges is 
$2\left(d-1\right)-h$, the number of separatrices. Each strip gives one face, each basin of 
a center 
gives one face, and each upper and lower half-plane gives a face. That is, 
altogether we 
have $s+c+2\frac{p_P}{2}$ faces. In total,
\begin{align}
\chi\left(\M\right)=&\left(1+d-c-\frac{p_P}{2}\right)-\left(2\left(d-1\right)-h\right)+\left(s+c+2\frac{p_P}{2}\right) \nonumber \\
=&1-d+2+h+s+\frac{p_P}{2}\nonumber \\
=&2
\end{align}
since $s=d-1-\frac{p_P}{2}-h$ (see Equation \ref{qshrelation}).  Therefore, $\M$ is homeomorphic to $\rs$. 
\end{proof}
\begin{corollary}
$\M$ is isomorphic to $\rs$.
\end{corollary}
\begin{proof}
Since $\M$ is a simply connected Riemann surface that is 
compact, we have that it must be conformally equivalent to the Riemann sphere $\rs$ by The 
Uniformization Theorem.
\end{proof}
We now analyze the vector field $\xi_{\M}$ on $\M$.
\section{The Associated Vector Field $\xi_{\M}$}
\label{assvfs}
For $U_{[\ell]}$ such that $[\ell]$ is an odd or even \eqcl , the vector field $\xi_{\M}$ is defined on $U_{[\ell]}^{\ast}$ by
$\left(\eta_{[\ell]}\right)_{\ast}\left(\xi_{\M}\right)=\left(\varphi_{[\ell]}\right)_{\ast}\left(\vf\right)$.
On $\Sigma_{[\ell]} \cup \{ \mathbb{R}_{\pm} \cup
 \left( \rho_{[\ell]}+\mathbb{R}_{\pm}\right) \}$, there is a constant vector field $\vf$ 
that conjugates by $\varphi_{[\ell]}$ to the
linear vector field $\frac{\mp 2 \pi i}{\rho_{[\ell]}}z \vf$ in 
$V_{[\ell]}
\setminus \{ 0\}$.  For $[\ell]$ odd, the puncture $-\infty_{[\ell]}\in \M$ corresponds to a sink with multiplier $\frac{-2 \pi i }{\rho_{[\ell]}}$, and for $[\ell]$ even, the puncture $+\infty_{[\ell]}\in \M$ corresponds to a source with multiplier $\frac{2 \pi i }{\rho_{[\ell]}}$.  We see this by observing that
\begin{equation}
\left(\eta_E\right)_{\ast}\left(\xi_{\M}\right)=\left(\varphi_{[\ell]}\right)_{\ast}\left(\vf\right)=\frac{\mp 2 \pi i}{\rho_{[\ell]}}z \cdot \vf
\end{equation}
 since  $\frac{\mp 2
\pi i}{\rho_{[\ell]}}\varphi_{[\ell]}\left(z\right)=\frac{d}{dz}\varphi_{[\ell]}\left(z\right)$.  The map $\varphi_{[\ell]}$ 
extends holomorphically to 
$V_{[\ell]}$ by taking the value 0 at $z=0$. Hence, $\xi_{\M}$ extends holomorphically to the puncture $\infty_{[\ell]}$ by $\xi_{\M}=0$.\par
For $U_{\{[\ell]\}}$ corresponding to a closed $H$-chain $H_{\{[\ell]\}}$, 
we define the vector field $\xi_{\M}$ on $U_{\{[\ell]\}}^{\ast}$ by
$\left(\eta_{\{[\ell]\}}\right)_{\ast}\left(\xi_{\M}\right)=\left(\varphi_{\{[\ell]\}}\right)_{\ast}\left(\vf\right)$.
On the vertical half-strips $C_{\{[\ell]\}}$, there is a constant vector field 
$\vf$ that conjugates by $\varphi_{[\ell]}$ to the
linear vector field $\frac{\mp 2 \pi i}{\tau_{\{[\ell]\}}}z \vf$ 
in $\D^{\ast}$. The puncture corresponds to a center with multiplier $\frac{\mp 2 \pi i}{\tau_{\{[\ell]\}}}$. The map $\varphi_{\{[\ell]\}}$ extends holomorphically to 
$\mathbb{D}$ by taking the value 0 at $z=0$. Hence, $\xi_{\M}$ extends holomorphically to the puncture $\infty_{\{[\ell]\}}$ by $\xi_{\M}=0$.\par
The induced vector field in a neighborhood of $E$ can be calculated: 
\begin{equation}
\left(\eta_{E}\right)_{\ast}\left(\xi_{\M}\right)=\left(\phi_E\right)_{\ast}\left( \vf \right)=\frac{1}{d-1}w^{-\left(d-2\right)}\vfw 
\end{equation}
in $\D^{\ast}_{\frac{\sqrt[d-1]{r}}{d-1}}$ since the $d-1$ covering
of $\D^{\ast}_r$ has constant vector field $\vf$.  Therefore, the vector field  $\xi_{\M}$ has a pole of order $d-2$ at $E$.\par
For $U_{[m]}$ corresponding to a mixed \eqcl\ $[m]$, we wish to show that the vector field in the chart $\eta_{[m]}:U^{\ast}_{[m]}\rightarrow V_{[m]}^{\ast}$ is of order $\mathcal{O}\left(z^{p[m]/2+1}\right)$ in $V_{[m]}^{\ast}$, a neighborhood of zero. That is, we will show $\infty_{[m]} \in \M$, corresponding to the puncture in $U^{\ast}_{[m]}$ is a \multeq\ of multiplicity $p_{[m]}/2+1$ for $\xi_{\M}$.\par
We will use Koebe's Distortion Theorem and the extension of $\Phi$ from Section \ref{vfmult} to  prove the following lemma that $\Phi'\left(z\right)$ is uniformly bounded in $V$. 
\begin{lemma}
\label{phibdd}
Suppose we have $\Phi$ as above in $U$.  
Then there exist constants $k_1$ and $k_2$ in $\R_+$ such that 
\begin{equation}
k_1<|\Phi'\left(z\right)|<k_2, \quad z \in V.
\end{equation}
\end{lemma}
We first need to use a modified version of Koebe's Distortion Theorem.  It is stated in \cite{Pom}, for functions univalent in $\D$ and of the form $g\left(z\right)=z+a_2z^2+...$,
\begin{equation}
\frac{|z|}{\left(1+|z|\right)^2}\leq |g\left(z\right)|\leq \frac{|z|}{\left(1-|z|\right)^2}. 
\end{equation}
This is restricted to functions with $g\left(0\right)=0$ and $g'\left(0\right)=1$ defined in $\D$ which we want to extend to functions univalent in a disk $\D_{R}\left(z_0\right)$ centered at $z_0$ with radius $R$. We show that if $f$ is univalent in $\{z \mid |z-z_0|<R \}$ then
\begin{equation}
\frac{|z-z_0|/R}{\left(1+|z-z_0|/R\right)^2} \leq \frac{|f\left(z\right)-f\left(z_0\right)|}{|Rf'\left(z_0\right)|} \leq \frac{|z-z_0|/R}{\left(1-|z-z_0|/R\right)^2}.
\end{equation}
This is true by the following. Let $h:\{|\tilde{z}-f\left(z_0\right)|<Rf'\left(z_0\right)  \} \rightarrow \D$, $h\left(\tilde{z}\right)=\frac{\tilde{z}-f\left(z_0\right)}{Rf'\left(z_0\right)}$ and note that the mapping $z \mapsto \frac{z-z_0}{R}$ maps the disk centered at $z_0$ and with radius $R$ to the disk $\D$. Let $g:\D\rightarrow \D$ be defined by $g\left(w\right)=h \left( f \left( Rw+z_0\right) \right)$. This is the composition of univalent functions, so $g$ is univalent. Furthermore, $g(0)=0$ and $g'(0)=1$, so $g$ must take the form $g(w)=w+h.o.t$,  and hence by Koebe's Distortion Theorem, 
\begin{equation}
\frac{|w|}{\left(1+|w|\right)^2}\leq \frac{|f\left(Rw+z_0\right)-f\left(z_0\right)|}{|R f'\left(z_0\right)|} \leq \frac{|w|}{\left(1-|w|\right)^2}.
\end{equation}
Remembering that $w=\frac{z-z_0}{R}$ gives the result.\par
\begin{proof}[Proof of Lemma \ref{phibdd}]
It follows from the above that for some $z \in V$, if $\Phi$ is univalent in $\{\zeta \mid |\zeta-z|<R  \}$, where $\zeta=z+v(z)=\tilde{f}(z)$, then
\begin{equation}
\frac{|v\left(z\right)|}{\left(1+|v\left(z\right)|/R\right)^2} \leq \frac{|\Phi\left(z+v\left(z\right)\right)-\Phi\left(z\right)|}{|\Phi'\left(z\right)|} \leq \frac{|v\left(z\right)|}{\left(1-|v\left(z\right)|/R\right)^2},
\end{equation}
and hence
\begin{equation}
\frac{\left(1-|v\left(z\right)|/R\right)^2}{|v\left(z\right)|}\leq \left| \Phi'\left(z\right) \right|\leq \frac{\left(1+|v\left(z\right)|/R\right)^2}{|v\left(z\right)|}.
\end{equation}
From the condition $|v\left(z\right)-1|<1/4$, we have $3/4<|v(z)|<5/4$, and to have $|\Phi'|$ bounded away from 0 and $\infty$, we need therefore $R\geq 5/4 +\epsilon$. Since $\Phi$ is univalent in $U$, we have
\begin{equation}
\frac{\left(2/5 + \epsilon_1\right)^2}{5/4}<\left| \Phi'\left(z\right) \right|< \frac{\left(2-\epsilon_2\right)^2}{3/4}
\end{equation}
satisfied for all $z \in V \cap \{z \mid \Im(z)<c_4-(5/4+\epsilon)  \}$. 
\end{proof}
 We next examine the vector field in $U^{\ast}_{[m]}$ to determine what form $\xi_{\M}$ takes at $\infty_{[m]}$.  We will show that $\left(\eta_{[m]}\right)_{\ast}\left(\xi_{\M}\right)=o\left(z\right)$ for $z \in V_{[m]}^{\ast}$.\par 
From $\pi _{\ast }\left(j \vf  \right)=\vfw$,  the vector field $j \vf$ in $S$ is expressed by
\begin{equation}
 j\left(z\right)=\frac{1}{ \pi i p_{[m]} \exp\left( \pi i p_{[m]} z\right)}.
\end{equation}
The expression of the vector field $g \frac{d}{dz}$ in  $\Phi$ coordinates is determined by 
\begin{equation}
g\left(\Phi\left(z\right)\right)=\Phi'\left(z\right)j\left(z\right).
\end{equation}
Now by Lemma \ref{phibdd} $|\Phi'|$ is bounded away from 0 and $\infty$, then
\begin{equation}
\Phi'\left(z\right)j\left(z\right) =\mathcal{O}\left( j\left(z\right)\right) =\mathcal{O}\left( \exp\left(- \pi i p_{[m]}z\right)\right).
\end{equation} 
Now under the map $\exp\left(-2 \pi i w\right)$, the vector field $f \frac{d}{dz}$ in the punctured disk is determined by
\begin{equation}
  f\left(\exp\left(-2 \pi i w\right)\right)=-2 \pi i \exp\left(-2 \pi i w\right) g\left(w\right) =\mathcal{O}\left( \exp\left(-2 \pi i ( w + \left(p_{[m]}/2\right)z )\right)\right),
\end{equation}
which gives that $f =o\left( z\right)$ in the punctured disk. Hence, $\left(\eta_{[m]}\right)_{\ast}\left(\xi_{\M}\right)=o\left(z\right) \vf$, where the puncture at 0 in $V_{[m]}^{\ast}$ corresponds to the puncture in $U^{\ast}_{[m]}$ which can be uniquely extended to $U_{[m]}$ by addition of the point $\infty_{[m]}$.  The vector field $\xi_{\M}$ can be extended holomorphically to $\infty_{[m]}$, i.e. such that $\left(\eta_{[m]}\right)_{\ast}\left(\xi_{\M}\right)=0\vf$ at 0 in $V_{[m]}$. \par
\begin{figure}%
   \centering
   \resizebox{!}{6cm}{\input{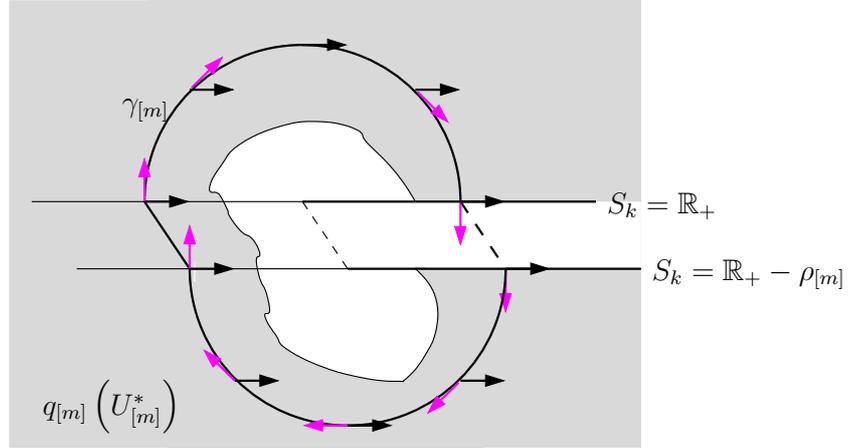}}%
   \caption{Possible depiction of $q_{[m]}\left(U^{\ast}_{[m]}\right)$ in rectifying coordinates  (shaded area) with $p_{[m]}=2$ and $\Im \left(\rho_{[m]}  \right)>0$. The vector field $\vf$ is tangent to the curve $\gamma_{[m]}$ in exactly $p_{[m]}$ places. Compare with Figure \ref{index2}.}
  \label{indexfig1}
\end{figure}
\begin{figure}%
   \centering
   \input{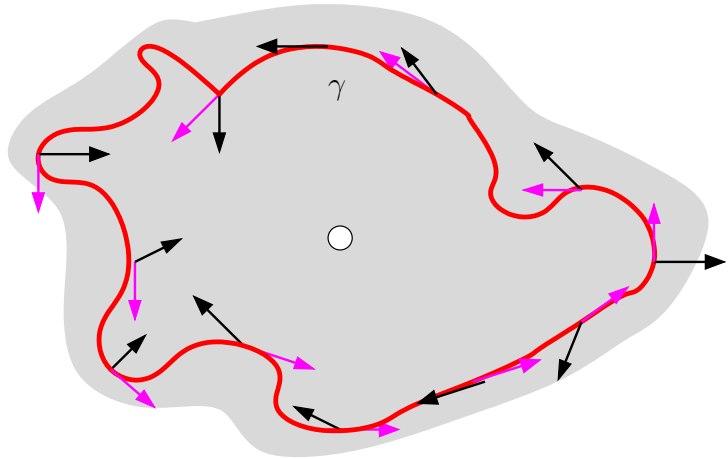}%
   \caption{Possible depiction of the index of $\left(\eta_{[m]}\right)_{\ast}\left(\xi_{\M}\right)$ in a neighborhood of zero where the multiplicity of $\infty_{[m]}$ is 2. The light vectors are tangent vectors to $\gamma$, and the black vectors belong to the vector field at points along $\gamma$. Notice that the vector field is tangent to $\gamma$ exactly twice: once pointing in the same direction as the tangent vector and once against. Compare with Figure \ref{indexfig1}.}
  \label{index2}
\end{figure}

Now we show that the multiplicity of the zero is in fact $\left(p_{[m]}/2+1\right)$. We do this by an index argument. Consider the (piecewise smooth) Jordan curve about $\infty_{[m]}$ in $U_{[m]}$ that corresponds to the curve $\gamma_{[m]}$ in rectifying coordinates with half-circles in each upper and lower half-plane and appropriate line segments in the half-strips  with clockwise orientation so that $\infty_{[m]}$ is to the left of the curve (see Figure \ref{indexfig1}). 
The curve $\gamma_{[m]}$ maps to another Jordan curve $\gamma$ in $V_{[m]}^{\ast}$ under $\eta_{[m]}$ since $\eta_{[m]}$ is univalent.  Moreover,  $\left(\eta_{[m]}\right)_{\ast}\left(\xi_{\M}\right)$ is never 0 along $\gamma$.  The angles between $\left(\eta_{[m]}\right)_{\ast}\left(\xi_{\M}\right)$ on $\gamma$ and the tangent vectors of $\gamma$ are the same as the angles between the vector field $\vf$ and the tangent vectors of $\gamma_{[m]}$ since $\eta_{[m]}$ is conformal (see Figure \ref{index2}).  
This implies that on $\gamma$, $\left(\eta_{[m]}\right)_{\ast}\left(\xi_{\M}\right)$ will be tangent to $\gamma$ exactly $p_{[m]}$ times, and alternating with pointing along the orientation of $\gamma$ and against the orientation of $\gamma$ when we travel along $\gamma$.  Along the orientation of $\gamma$, on the arcs between the tangents going along to the tangents going against, the vector field must be pointing inward, since this corresponds to what happens in rectifying coordinates.  This implies that the vector field must be rotating in the same orientation as $\gamma$, otherwise there would be a place in this arc where there is another tangent.  So with respect to the tangents of $\gamma$, $\left(\eta_{[m]}\right)_{\ast}\left(\xi_{\M}\right)$ has rotated $+p_{[m]}/2$ times.  We must add one more time around, accounting for the index of $\gamma$.  This gives that $\left(\eta_{[m]}\right)_{\ast}\left(\xi_{\M}\right)$ has index $p_{[m]}/2+1$ at 0.\par
We have that $\xi_{\M}$ is holomorphic on $\M \setminus \{ E \}$ since 
it can be 
expressed as $g \vf$ where $g$ is holomorphic in the domain of each chart.

\section{Proof of the Structure Theorem}
\label{finalpf}
We now show that there is a unique conformal isomorphism $\Phi:\M \rightarrow \rs$ that 
induces the 
vector field $\xi_P\in \Xi_d$ having the given invariants.
\begin{theorem}
\label{fundthm}
There exists a unique conformal isomorphism $\Phi:\M \rightarrow \rs$ such that
$\Phi\left(E\right)=\infty$ and such that $\Phi_{\ast}\left(\xi_{\M}\right)=\xi_{P}$ where $\xi_P \in
\Xi_d$ and $\xi_P$ has the given invariants.
\end{theorem}
\begin{proof}
 By the Uniformization Theorem, there exists an isomorphism $\Psi:\M \rightarrow \rs$, which is unique up to post composition by a M\"{o}bius transformation and induces a vector field  $\Psi_{\ast}\left(\xi_{\M}\right)$ defined on $\rs$.  Choose $\Psi$ such that $\Psi\left(E\right)=\infty$. Then  $\Psi_{\ast}\left(\xi_{\M}\right)$ is a holomorphic vector field on $\C$, expressed in canonical coordinates as $g \vf$. We know $g$ has a unique pole of order $d-2$ at $\infty$ since 
$\left(\eta_E\right)_{\ast}\left(\xi_{\M}\right)=\frac{-1}{z^{d-2}}\vf$ for $z \in V_0$.
\begin{lemma}
If $g\vf$ is holomorphic in $\C$ and $g\vf$ has a pole of order $d-2$ at $\infty$, then $g\vf$ is a polynomial vector field of degree $d$.
\end{lemma}
\begin{proof}
Since $g$ is entire, it admits a Taylor series expansion
\begin{equation}
g\left(z\right)=\sum \limits_{n=0}^{\infty}a_n z^n, \quad z\in \C.
\end{equation}
To examine the behavior near $\infty$, we need to make a transformation by $w=\varphi(z)=\frac{1}{z}$ using the transformation rule for vector fields
\begin{equation}
\tilde{g}\left(w\right)=\left(-w^2\right)\sum \limits_{n=0}^{\infty}a_{n} \left(\frac{1}{w} \right)^n, \quad w\in \C^{\ast}.
\end{equation}
So if $\tilde{g} \vfw$ has a pole of order $d-2$ at $w=0$, then $a_n=0$ for $n>d$, $a_d\neq 0$, and  it follows that $g$ is a polynomial of degree $d$.
\end{proof}
If $g$ has leading coefficient $a_d$, then for $A=\sqrt[d-1]{a_d}$, $g \vf$ is conjugate to a monic polynomial vector field whose incoming separatrices have asymptotic directions $\delta_j=\exp\left(2\pi i\frac{ j}{2\left(d-1\right)}\right)$ where $j \in\{0,2,\dots,2d-4  \}$.  We choose $A$ such that $A\Psi(S_0)$ is asymtotic to $\R_+$ near infinity\footnote{See Remark \ref{seplabelremark}.}.  We choose $B$ such that
\begin{equation}
\label{centering}
\sum \limits_{H_{\{[\ell]\}}\subset H}\left(\Psi\left( \infty_{\{[\ell]\}}\right)+B\right)+\sum \limits_{[\ell]\subset L}\left(p_{[\ell]}/2+1\right)\left(\Psi\left( \infty_{[\ell]}\right)+B\right)=0,
\end{equation}
centering the \eqpt s.
Let 
\begin{equation}
\sum \limits_{H_{\{[\ell]\}}\subset H}\Psi\left(\infty_{\{[\ell]\}}\right)+\sum \limits_{[\ell]\subset L}\left(p_{[\ell]}/2+1\right)\Psi\left(\infty_{[\ell]}\right)=\tilde{c}.
\end{equation}
If $\tilde{c}=0$, we are finished.  If $\tilde{c}\neq 0$, then Equation \eqref{centering} gives the equation $\tilde{c}+dB=0$, giving $B=\frac{-\tilde{c}}{d}$. Let $\tilde{A}$ be the unique affine map with $A$ and $B$ as described above.  Then $\Phi=\tilde{A} \circ \Psi$ is the unique isomorphism conjugating $g \vf$  to a 
monic,
centered polynomial vector field $P \vfw$ of degree $d$.\par
The combinatorial invariant is preserved since $\Phi$ is a conformal isomorphism such that $s_0=\Phi(S_0)$ has $\R_+$  as asymptote.\par 
The vector fields $\xi_{\M}$ and $\Phi_{\ast}\left(\xi_{\M}\right)=P \vfw$ must have the same number of zeros with the same multiplicities (and at the images $\Phi\left(\infty_{[\ell]}\right)$ and $\Phi\left(\infty_{\{[\ell]\}}\right)$) and the same dynamical residues (see, for example, Theorem 1 in \cite{BT76}), hence the same analytic invariants.\par
Therefore, $\left(\sim_P,H_P\right)=\left(\sim,H\right)$, and $A\left(\xi_P\right)=A\left(\xi_{\M}\right)$ as desired.
\end{proof}

\bibliographystyle{halpha}
\bibliography{article1}{}
Department of Mathematics \\
Technical University of Denmark\\
Building 303S\\
2800 Kgs. Lyngby, Denmark\\
\\
e-mail: B.Branner@mat.dtu.dk\\
e-mail: K.Dias@mat.dtu.dk
\end{document}